\documentclass[12pt]{article}
\usepackage{amssymb}
\usepackage{amsmath}
\usepackage{color}
\usepackage{graphics}
\usepackage{epsfig}
\usepackage{hyperref}
\usepackage{amsthm,amsmath,amssymb}
\usepackage{mathrsfs}

\newtheorem{claim}{\bf \t}[part]

\newtheorem{Lemma}{Lemma}[part]
\newtheorem{Proposition}{Proposition}[part]
\newtheorem{Remark}{Remark}[part]
\newtheorem{Theorem}{Theorem}[part]

\numberwithin{Assumption}{section} \numberwithin{Corollary}{section}
\numberwithin{Definition}{section} \numberwithin{equation}{section}
\numberwithin{Example}{section} \numberwithin{Lemma}{section}
\numberwithin{Proposition}{section} \numberwithin{Remark}{section}
\numberwithin{Theorem}{section}

\def \Sum{\displaystyle\sum}

\def \ess{\mathrm{esssup}}

\def\t{\theta}

\def\text#1{{\rm #1}}
\begin{document}
\date{}
\title{\Large \bf Local  well-posedness of strong
solutions to the 2D
nonhomogeneous primitive equations with density-dependent viscosity}
\author{\small \textbf{Quansen Jiu},$^{a}$
\thanks{
E-mail: jiuqs@mail.cnu.edu.cn.}\quad
\textbf{Lin Ma}$^{a}$,
  \thanks{
E-mail: malin.cnu@foxmail.com.}\quad
and\quad \textbf{Fengchao Wang}$^{b}$
\thanks{
E-mail: wfcwymm@163.com.}
}
 \maketitle
\small $^a$ School of Mathematical Sciences, Capital Normal University, Beijing
100048, P.~R. China

\small $^b$ College of Mathematics and Physics, Beijing University of Chemical Technology, Beijing 100029, P.~R. China.

 {\bf Abstract:} In this paper, we consider the initial-boundary value problem of the nonhomogeneous primitive equations with density-dependent viscosity. Local  well-posedness of strong
solutions is established for this system with a natural compatibility condition. The initial density does not need to be strictly positive and may contain vacuum.  Meanwhile, we also give the corresponding blow-up criterion if the maximum existence interval with respect to the time is finite.

{\bf Key Words:} Nonhomogeneous primitive equations, local solutions, density-dependent viscosity, blow-up criterion, vacuum

\section{Introduction } \setcounter{equation}{0}
\setcounter{Assumption}{0} \setcounter{Theorem}{0}
\setcounter{Proposition}{0} \setcounter{Corollary}{0}
\setcounter{Lemma}{0}

 \ \ \ \
The primitive equations of the ocean and the atmosphere are considered to be a fundamental model in geophysical fluid dynamics (see, for example,
Pedlosky \cite{ped} and Vallis \cite{val}). They are described by a system of equations which are derived from the Navier-Stokes or Boussinesq equations by assuming that the vertical motion is modeled by the hydrostatic balance.

In this paper we consider the two dimensional   nonhomogeneous primitive equations in $(0,T)\times \Omega$, which are given by
\begin{equation}\label{a1}
	\left\{
		\begin{aligned}
	&\rho_{t}+\partial_{x}(\rho u)+\partial_{y}(\rho w)=0,	\\
		&\partial_{t}(\rho u)+\partial_{x}(\rho u^{2})+\partial_{y}(\rho u w)+\partial_{x}P-\partial_{x}(\mu(\rho)\partial_{x}u)-\partial_{y}(\mu(\rho)\partial_{y}u)=\rho f, \\
	&	\partial_{x}u+\partial_{y}w=0, \\
	&	\partial_{y}P=0.\\
	\end{aligned}
\right.
\end{equation}
Here $(x,y) \in \Omega=M\times (0,1)$, where $M=[0,L]$ and $L>0$. The horizontal velocity $u$, the vertical velocity $w$, the density $\rho$ and the pressure $P$  are the unknowns. And $f$ is a given external force.
The viscosity coefficient $\mu=\mu(\rho)$ is  a function of the density $\rho$ and is assumed to satisfy
\begin{equation}\label{12}
	\mu\in C^{2}[0,\infty) \quad \text{and} \quad \mu\geq\underline{\mu}>0 \quad\text{on}\quad [0,\infty).
\end{equation}

The first mathematical and systematical studies of the (homogeneous) primitive equations were made by Lions, Temam and Wang in  \cite{lions1992-1}-\cite{lions1995}
in which they established the global existence of weak solutions to the system for initial data $a\in L^{2}$. Weak solutions in 2D turn out to be unique, see  Bresch, Kazhikjov and Lemoine \cite{bd} and Kukavica, Pei, Rusin and Ziane \cite{Kuk2014}. However, the uniqueness of weak solutions in 3D case is open up to now. Concerning the strong solutions, Cao and Titi \cite{cao1} proved a breakthrough result for the 3D viscous primitive equations. In \cite{cao1}, they established the global well-posedness of strong solutions
to the 3D incompressible primitive equations by making full use of the hydrostatic balance to exploit the two dimensional structure of the key part of the pressure and decomposing the velocity into barotropic
and baroclinic components, see also Kobelkov \cite{kob} and Kukavica and Ziane \cite{kuk7}-\cite{kuk8}. Modifications of the primitive equations dealing with either only horizontal viscosity and diffusion or with horizontal or vertical eddy diffusivity were recently investigated by Cao and Titi in \cite{cao2}, by Cao, Li, and  Titi in \cite{cao3}-\cite{cao6}. Here, global well-posedness results are established for initial data belonging to $H^{2}$.
For the well-posedness results in the
$L^{p}$ type spaces based on the maximal regularity technique, one can see the works by Hieber and Kashiwabara \cite{hie} and Giga, Gries, Hieber, Hussein and Kashiwabara \cite{giga20} or \cite{giga21}.

The mathematical analysis of the compressible primitive equations model, meaning CPEs, has recently been attracting the attention of many mathematicians. Gatapov and Kazhikhov \cite{gag} proved the global existence of weak solutions in 2D case. The uniqueness of such weak solutions was proved by Jiu, Li and Wang in \cite{jiu}. Ersoy, Ngom and Sy \cite{ers2011} showed the stability of weak solutions with the viscosity coefficients depending on the density for 3D CPEs, see also Tang and Gao \cite{tang} for 2D case. The stability implies that approximate solutions will converge to a weak solution if it satisfies some uniform bounds. Wang, Dou and Jiu \cite{wang} proved the global existence of weak solutions to 3D CPEs with gravity field, see also Liu and Titi \cite{liu2019}. Concerning strong solutions, Liu and Titi \cite{liu2021} proved the local existence and uniqueness of strong solutions to the 3D CPEs. Very recently, Wang \cite{wang2021} gave local existence and uniqueness of strong solutions to the 2D CPEs with density-dependent viscosity.  Meanwhile,  Liu and Titi consider the zero Mach number limit of CPEs under well-prepared initial data or ill-prepared initial data, see \cite{liu2020} or \cite{liu2023}, respectively.

 As intermediate equations between (homogeneous) primitive equations and compressible primitive equations, the nonhomogeneous primitive equations were recently studied. When initial vacuum is taken into account and $\mu$ is still a constant, Jiu and Wang $\cite{wang2}$ first proved the local existence of the unique strong solution to 2D nonhomogeneous  primitive equations for sufficiently regular data. Furthermore, Wang, Jiu and Xu $\cite{wang1}$ established a global existence result on strong solutions with nonnegative density in case that $\mu$ is a constant for initial data with small $H^{\frac{1}{2}}$-norm. However, there have been  few studies on the nonhomogeneous primitive equations with density dependent viscosity.

In this paper, we consider the initial-boundary value problem of the nonhomogeneous primitive equations (\ref{a1}) with the viscosity $\mu(\rho)$ which depends on density $\rho$.

The boundary conditions on the boundary of $\Omega$ are imposed as
\begin{equation}\label{a2}
	\begin{aligned}
		&\rho,~ u,~ w,~ P~ \text {are~ periodic~ in~the ~direction}~ x,\\
		& u|_{y=0}=u|_{y=1}=0,\\
		&w|_{y=0}=w|_{y=1}=0,
	\end{aligned}
\end{equation}
and the initial data of  velocity and density are
\begin{equation}\label{14}
	\begin{aligned}
		& u|_{t=0}=u_{0}(x,y),\\
		& w|_{t=0}=w_{0}(x,y),\\
		&\rho|_{t=0}=\rho_{0}(x,y).\\
	\end{aligned}
\end{equation}
Also the following compatibility condition is imposed:\begin{equation}\label{1.5}\begin{aligned}
	&	-\rho_{0}u_{0}\partial_{x}u_{0}+\rho_{0}\int_{0}^{y}\partial_{x}u_{0
		}(s)ds\partial_{y }u_{0}-\partial_{x}P_{0}+\partial_{x}(\mu(\rho_{0})\partial_{x}u_{0})+\partial_{y }(\mu(\rho_{0})\partial_{y}u_{0})+\rho_{0}f_{0}:=\rho_{0} V_{1}, \\& \partial_{x}\int_{0}^{1}u_{0}dy=0.
	\end{aligned}
\end{equation}
By the divergence free condition, we can write $$w(x,y,t)=-\int_{0}^{y}\partial_{x}u(x,s,t)ds.$$

 We define the following function spaces for primitive equations  :\begin{equation*}
 	H^{1}_{0,\,per}(\Omega)=\{u\in H^1(\Omega); \partial_{x}\int_{0}^{1}u(x,y)dz=0,\,u|_{y=0,1}=0, u~ \text{is\ periodic \ in\ the\ direction}\ x\}.
 \end{equation*}

Now we state the main result as follows:
\begin{Theorem}\label{the-1}
 Suppose that the initial data $(\rho_0,u_0)$ satisfies the regularity condition
\begin{equation}\label{1711}
0\le \rho_{0}\in H^{2},\nabla \rho_{0}\in L^{\infty}, u_{0}\in H^{1}_{0,per}\cap H^{2}
\end{equation} and $ f \in H^{1}((0,T)\times\Omega)$, and the compatibility condition (\ref{1.5})
for some $(P_{0}, f_{0},V_{1})\in H^{1}\times L^{2}\times L^{2}$. Then there exists a small time $T_{*}\in(0,T)$ and a unique strong solution $(\rho, u, w ,P)$ to the initial boundary value problem $(\ref{a1})$-$(\ref{14})$ such that
\begin{equation}\label{17}
\begin{aligned}
&\rho\geq0,\quad \rho \in C([0,T_{*}];W^{1,\infty})\cap C([0,T_{*}];W^{2,2}),\quad\rho_{t}\in C([0,T_{*}];L^{q}),\quad \sqrt{\rho}u_{t}\in L^{\infty}(0,T_{*};L^{2}),\\
&P\in C([0,T_{*}];H^{1})\cap L^{2}(0,T_{*};H^{2}),\quad
\nabla u\in C([0,T_{*}];H^{1})\cap L^{2}(0,T_{*};H^{2}),\\
&u_{t} \in L^{2}(0,T_{*};H^{1}),\quad
w\in C([0,T_{*}];H^{1})\cap L^{2}(0,T_{*};H^{2}),\\
\end{aligned}
\end{equation}
where $1\leqslant q<\infty$. Furthermore, if $T^{*}$ is the maximal existence time of the local strong solution $(\rho,u,w,P)$ and $T^{*}<T$, then
\begin{equation}\label{16}
\sup_{0\leq t<T^{*}}(\left\|\nabla \rho(t)\right\|_{L^{\infty}}+\left\|\nabla^{2}\rho(t)\right\|_{L^{2}}+\left\|\nabla u(t)\right\|_{L^{2}})=\infty.
\end{equation}
\end{Theorem}
We intend to use the linearization and iteration method to obtain the existence and regularity results for the original system (\ref{a1}).
Due to the loss of the derivation of vertical velocity, the nonhomogeneous primitive equations  are slightly more complicated than the incompressible Navier-Stokes equations. One key issue is to choose a suitable working space in which we could close the estimates and prove the existence and uniqueness of the solutions. On the other hand, the final convergence could not be obtained easily due to the loss of the derivative order of the vertical velocity. Indeed, we deduce
\begin{equation*}
	\begin{aligned}
\frac{d}{dt}\int_{\Omega}\rho^{k+1}|\eta^{k+1}|^{2}+\int_{\Omega}|\nabla \eta^{k+1}|^{2}=C\|\sigma^{k+1}\|_{L^{2}}^{2}B_{k}(t)+C\|\sqrt{\rho^{k}}\eta^{k}\|_{L^{2}}^{2}+C\|\partial_{x }\eta^{k}\|_{L^{2}}^{2}
	\end{aligned}
\end{equation*}and
$$\|\sigma^{k+1}\|_{L^{2}}^{2}\leq C(\int_{0}^{t}\|\sqrt{\rho^{k}}\eta^{k}\|_{L^{2}}+\|\partial_{x }\eta^{k}\|_{L^{2}}^{2}ds),$$
where $\eta^{k+1}=u^{k+1}-u^{k}$, $\sigma^{k+1}=\rho^{k+1}-\rho^{k}$ and $B_{k}(t)\in L^{1}([0,T])$.
However, we except such result (see (\ref{441}) and (\ref{4422})): \begin{equation*}
	\begin{aligned}
		\frac{d}{dt}\int_{\Omega}\rho^{k+1}|\eta^{k+1}|^{2}\leq C\|\sigma^{k+1}\|_{L^{2}}^{2}B_{k}(t)+C\|\sqrt{\rho^{k}}\eta^{k}\|_{L^{2}}^{2}
	\end{aligned}
\end{equation*}and
$$\|\sigma^{k+1}\|_{L^{2}}^{2}\leq C\int_{0}^{t}\|\sqrt{\rho^{k}}\eta^{k}\|_{L^{2}}^{2}.$$ Thus, instead  we consider regularization problems (\ref{3.1}) and set the solutions as the approximate solutions to original equations (\ref{a1}). Later, we use standard methods to get the well-posedness of the original equations (\ref{a1}) and  the initial vacuum is allowed. To get higher order a priori estimates for the approximate solutions, we give the regularity results of the hydrostatic  Stokes equations for the general case of density-dependent viscosities, which is introduced precisely in Section 2. It turns out that the estimates in Lemma $\ref{Lemma-41}$ hold when $\rho\in W^{2,2}$ and $\nabla \rho \in L^{\infty}$. Thus, to expect the existence of strong solutions of $(\ref{a1})$-$(\ref{14})$, the initial  density $\rho_{0}$ must have at least the regularity
$\rho_{0}\in W^{2,2} $ and $\nabla \rho_{0} \in L^{\infty}$.

The strategy of our approach may be described as follows. In Section 3,  we mainly provide existence and regularity results for the regularization problems  (\ref{3.1}) which are solved by the linearization and iteration method. In a first step, we construct approximate solutions of (\ref{3.1}) by solving iteratively linearized problems for (\ref{3.1}). We may consider a problem for higher regularity of solutions. To settle this, we need to control the time  derivatives  of $u$ at initial time and hence to impose the compatibility condition (\ref{1.5}). Second, we derive (local in time) uniform bounds with respect to the norms of the approximate solutions of (\ref{3.1}). We prove it by introducing functionals
\begin{equation}
	\Phi_{K}(t)=\max _{1 \leqslant k \leqslant K}\left(1+\left\|\nabla \rho^{k}(t)\right\|_{W^{1,\infty}}+\left\|\nabla
	^{2}\rho^{k}(t)\right\|_{W^{1,2}}+\left\|\nabla u^{k}(t)\right\|_{L^{2}}+\left\|\nabla u^{k}_{t}(t)\right\|_{L^{2}}\right),
\end{equation}
and showing that each $\Phi_{K}(t)$ satisfies the integral inequality $$\Phi_{K}(t) \leqslant C \exp \left[C \exp \left(C \int_{0}^{t} \Phi_{K}(s)^{N} \mathrm{~d} s\right)\right],$$ for some positive constant $C,N>0$. And then  we can conclude the local bound of $\Phi_{K}(t)$ by applying a Gronwall-type argument (see Lemma \ref{Lemma-2}). In the processing of closing this integral inequality, it is more difficult to deal with the term  of $\left\|\nabla^{2}\rho\right\|_{W^{1,2}}$ and $\left\|\nabla\rho\right\|_{W^{1,\infty}}$. To overcome this difficulty, using the transport equation $(\ref{a1})_1$, we obtain	\begin{equation}\label{1.10}
	\|\nabla \rho(t)\|_{L^{\infty}}\leq \|\nabla \rho_{0}\|_{L^{\infty}}\exp(C\int_{0}^{t}\|u(s)\|_{H^{3}}ds),
\end{equation}
\begin{equation}\label{1.11}
	\|\nabla ^{2}\rho(t)\|_{L^{2}}\leq (\|\nabla \rho_{0}\|_{L^{\infty}}+\|\nabla ^{2}\rho_{0}\|_{L^{2}})\exp(C\int_{0}^{t}\|u(s)\|_{H^{3}}ds),
\end{equation}

\begin{equation}\label{1.12}	\| \nabla ^{2}\rho(t)\|_{L^{\infty}}\leq \|\nabla^{2}\rho_{0}\|_{L^{\infty}}\exp(C\int_{0}^{t}\|u(s)\|_{H^{4}}ds),
\end{equation}
and
\begin{equation}\label{1.13}
	\| \nabla ^{3}\rho\|_{L^{2}}\leq(\| \nabla ^{3}\rho_{0}\|_{L^{2}}+\|\nabla^{2}\rho_{0}\|_{L^{\infty}}+\|\nabla\rho_{0}\|_{L^{\infty}}) \exp(C\int_{0}^{t}\|u(s)\|_{H^{4}}ds).
\end{equation}	
It is easily seen that both the $W^{1,\infty}$-norm and the $W^{2,2}$-norm of the density $\rho$ need the  $H^{3}$-regularity of the velocity $u$ and both the $W^{2,\infty}$-norm and the $W^{3,2}$-norm of the density $\rho$ need the  $H^{4}$-regularity of the velocity $u$. This requires more regularity than the solutions to the Navier-Stokes equations.
According to the regularity results of the hydrostatic Stokes equations (see Lemma \ref{Lemma-41}), one has\begin{equation}\label{1.14}
\|u\|_{H^{3}}+\|P\|_{H^{2}} \leqslant C(1+\|\nabla\rho\|_{L^{ \infty}}+\|\nabla\rho\|_{L^{ \infty}}^{2})(1+\|\nabla^{2}\rho\|_{L^{2}})\|f\|_{H^{1}},
\end{equation}
and
	\begin{equation}\label{1.15}
	\|u\|_{H^{4}}+\|P\|_{H^{3}} \leqslant \widetilde{C}(1+\|\nabla\rho\|_{L^{ \infty}}+\|\nabla\rho\|_{L^{ \infty}}^{2}+\|\nabla^{2}\rho\|_{L^{\infty}}+\|\nabla ^{3}\rho\|_{L^{2}})^{2}(1+\|\nabla^{2}\rho\|_{L^{2}})\|f\|_{H^{2}}.
\end{equation}
Combining (\ref{1.10})-(\ref{1.13}) with (\ref{1.14})-(\ref{1.15}), we can obtain the estimate of $\Phi_{K}(t)$ and furthermore
 the uniform bound of the approximate solutions to (\ref{3.1}) (see Lemmas $\ref{Pro-1}$-$\ref{pro-8}$ for details). We remark that it is sufficient to solve the problems (\ref{a1}) when we set $$\Phi_{K}(t)=\max _{1 \leqslant k \leqslant K}\left(1+\left\|\nabla \rho^{k}(t)\right\|_{L^{\infty}}+\left\|\nabla
 ^{2}\rho^{k}(t)\right\|_{L^{2}}+\left\|\nabla u^{k}(t)\right\|_{L^{2}}\right).$$
 That is, we use the terms of  $\|\nabla^{2} \rho^{k}\|_{L^{\infty}},\|\nabla^{3}\rho^{k}\|_{L^{2}}$ and $\|\nabla u^{k}_{t}\|_{L^{2}}$ to deal with the difficulty caused by the term of $-\lambda\rho_{xx}$. Based on the  obtained uniform estimates, an usual approach to prove the existence of the solution is by proving that the approximate solutions to (\ref{3.1}) are Cauchy sequences in suitable spaces. However, it seems difficult to prove that the approximate solutions we construct are Cauchy sequences. A method in our  proof is that we apply the Aubin-Lions lemma to choose subsequences which converge strongly in suitable spaces and we prove that the neighbor subsequences $\{u^{k_j}\}_{j=1}^{\infty}$ and $\{u^{k_j-1}\}_{j=1}^{\infty}$ have same limit function which is in fact the strong solution we expect. A novelty in our proof of solving the difficulty caused by the loss of the derivation of vertical velocity is that we take the derivates of vertical velocity to the density and add the regular term $-\lambda\rho_{xx}$ on the $(\ref{a1})_{1}$ in the process of convergence.
 More precisely, multiplying by $\sigma^{k+1}$ and taking by parts in x-directions, we obtain $$\frac{d}{dt}\int_{\Omega}(\sigma^{k+1})^{2}\leq\|\nabla \rho^{k}\|_{W^{1,\infty}}\|\eta^{k}\|_{L^{2}}\|\sigma^{k+1}\|_{L^{2}}+\|\nabla \rho^{k}\|_{L^{\infty}}\| \eta^{k}\|_{L^{2}}\|\sigma_{x}^{k+1}\|_{L^{2}},$$
 where $ \sigma^{k+1}=\rho^{k+1}-\rho^{k}, \eta^{k+1}=u^{k+1}-u^{k}$. As in the above inequality, there is a term $\|\sigma_{x}^{k+1}\|_{L^{2}}$ which brings difficulty to prove that $\|\sigma^{k+1}\|_{L^{2}}\rightarrow0, k\rightarrow \infty$. To overcome this difficulty, we add a new term $-\lambda \rho_{xx}$. Then, we have \begin{equation}
 \frac{d}{dt}	\|\sigma^{k+1}\|_{L^{2}}^{2}\leq  C \|\sqrt{\rho^{k}}\eta^{k}\|_{L^{2}}^{2}+C\|\sigma^{k+1}\|_{L^{2}}^{2}.
 \end{equation}
  Next, we have (see (\ref{d}))  $$	\|\sqrt{\rho^{k+1}}\eta^{k+1}(t)\|_{L^{2}}^{2}\leq C\int_{0}^{t}\|\sqrt{\rho^{k}}\eta^{k}(s)\|_{L^{2}}^{2}ds,$$
  where use the fact that $\rho^{k}\geq \delta.$
 Thus, we can prove that $u^{k+1}$ and $u^{k}$ have the same limit by a simple recursive argument. Moreover, we can use Aubin-Lions lemma to find some convergence subsequences, denoted again by $(u^{k+1}, u^{k}, \rho^{k+1}, P^{k})$, to show the convergence successively. Lastly, we establish the uniqueness of strong solutions to problems $(\ref{3.1})$. We remark that the estimates in Section 3 are independent on $\lambda$. In Section 4, We regularize the initial data in Theorem \ref{the-1} so that it satisfies the assumptions in the Proposition \ref{pro-1}. We then get some uniform estimates that have nothing to do with the lower bound of density $\delta$. Recall that the estimates in Section 3 are independent on $\lambda$, so we can finish the proof of the Theorem \ref{the-1} by applying standard regularizing techniques and compactness arguments, that is, doing a two-level approximation.
 We still give the strong continuity of time and a blow-up criterion to the system $(\ref{a1})$-$(\ref{14})$. It is remarked that the vertical velocity losses of one derivative but the integral in vertical direction balances out one derivative by using Poincar\'e inequality, which makes a help to prove the uniqueness.

The rest of this article is organised as follows: In Section 2, we present some preliminary lemmas which will be used frequently throughout the paper.  In Section 3, We give the higher regularity of $u$ which is deduced from the hydrostatic Stokes system and prove the local existence of unique strong solution to the regularization problems $(\ref{3.1})$ for the positive case. Based on these results, we construct the approximate solutions to problems $(\ref{a1})$-$(\ref{14})$. In Section 4, We complete the proof of Theorem \ref{the-1}.

{\bf Notations and conventions} The following notations will be used:\\
(1) $\|f\|_{L_{y}^{2}}=(\int_{0}^{1}|f|^{2}dy)^{\frac{1}{2}}$,\\
(2) $\|f\|_{L_{y}^{2}L_{x}^{\infty}}=(\int_{0}^{1}\|f\|^{2}_{L_{x}^{\infty}}dy)^{\frac{1}{2}}$,
$\|f\|_{L_{y}^{\infty}L_{x}^{\infty}}=\ess_{(x,y)\in \Omega} |f| $,\\
(3)  For $ 1<p<\infty,\ \|f\|_{L^{p}}=\|f\|_{L_{y}^{p}L_{x}^{p}}=\|f\|_{L_{x}^{p}L_{y}^{p}}=(\int_{0}^{1}\int_{0}^{L}|f|^{p}dxdy)^{\frac{1}{p}}$,\\
(4)$\int_{\Omega}=\int_{\Omega} \quad dxdy$.

 The positive constant $C$, which may be different from line to line, depends on $\|\mu\|_{C^{2}}$, $\Omega$, $\|\nabla^{2} \rho_{0}\|_{L^{2}}$, $\| \rho_{0}\|_{W^{1,\infty}}$, $\|u_{0}\|_{H^{1}}$, $\|f\|_{H^{1}}$ and $T$.

\section{Preliminaries }
\setcounter{equation}{0}
\ \ \ \ \
In this section, we present some basic facts needed later.

We have the following lemma about the regularity result of the hydrostatic Stokes system in $L^{2}$-setting:
\begin{Lemma}\label{Lemma-1}
(see $\cite{ziane}$) Let $\mathcal{O}=\prod_{\alpha=1}^{n-1}\left(0, L_{\alpha}\right)$ with $L_{\alpha}>0, \alpha=1,2, \cdots, n-1$, and let\\ $\Omega=\left\{\left(x^{\prime}, z\right) \in \mathcal{O} \times R ; 0<z<h\left(x^{\prime}\right)\right\}$, where $x^{\prime}=\left(x_{1}, \cdots, x_{n-1}\right)$ and $0<h_{1}<h\left(x^{\prime}\right)<h_{2}$ with some positive constants $h_{1}, h_{2}$. Assume that $h \in C^{k+3}(\mathcal{O})$ and that $h$ is periodic with respect to $x_{1}, x_{2}, \cdots, x_{n-1}$ with periods $L_{1}, L_{2}, \cdots, L_{n-1}$ respectively. If $f \in\left(H_{p}^{k}(\Omega)\right)^{n-1},~ k \geq 0$, then the hydrostatic Stokes system
$$
\left\{\begin{array}{l}
-\nu_{h} \triangle_{x^{\prime}} u\left(x^{\prime}, z\right)-\nu_{v} \partial_{z z} u\left(x^{\prime}, z\right)+\nabla_{x^{\prime}} P\left(x^{\prime}\right)=f\left(x^{\prime}, z\right), \quad\left(x^{\prime}, z\right) \in \Omega ,\\
\operatorname{div}_{\mathrm{x}^{\prime}} \int_{0}^{h\left(x^{\prime}\right)} u\left(x^{\prime}, \tau\right) \mathrm{d} \tau=0, \\
u \text { ~is~ periodic  ~with ~ respect ~to~} x_{1}, x_{2}, \cdots, x_{n-1}, \\
u\left(x^{\prime}, 0\right)=u\left(x^{\prime}, h\left(x^{\prime}\right)\right)=0,
\end{array}\right.
$$
possesses a unique solution $(u, P)$ which belongs to $\left(H_{p}^{k+2}(\Omega)\right)^{n-1} \times H_{p}^{k+1}(\Omega)$, where the Sobolev space $H_{p}^{m}(\Omega)$ defined as $H_{p}^{m}(\Omega)=\left\{u \in H^{m}(\Omega) ; \partial^{l} u / \partial x_{i}^{l}, i=1,2, \cdots, n-1, l=0,1, \cdots, m,\right.\\\left.\text{ is}~ \mathcal{O} \text{-periodic} \right\}.$ Moreover, we have the continuous dependence, that is, there exists a constant $C=C\left(\Omega, \nu_{h}, \nu_{v}\right)>0$ such that
$$
\|u\|_{H_{p}^{k+2}(\Omega)}+\|P\|_{H_{p}^{k+1}(\Omega)} \leq C\|f\|_{H_{p}^{k}(\Omega)}.
$$
\end{Lemma}
We state the corresponding results for the $L^{p}$ case.
\begin{Lemma}\label{the31} Let $\mathcal{O}=\prod_{\alpha=1}^{n-1}\left(0, L_{\alpha}\right)$ with $L_{\alpha}>0, \alpha=1,2, \cdots, n-1$, and let\\ $\Omega=\left\{\left(x^{\prime}, z\right) \in \mathcal{O} \times R ; 0<z<1\right\}$, where $x^{\prime}=\left(x_{1}, \cdots, x_{n-1}\right).$
  If $p \in (1,\infty)$ and $f \in L^{p}(\Omega),$ then the system
$$
\left\{\begin{array}{l}
	- \triangle_{x^{\prime}} u\left(x^{\prime}, z\right)-\partial_{z z} u\left(x^{\prime}, z\right)+\nabla_{x^{\prime}} P\left(x^{\prime}\right)=f\left(x^{\prime}, z\right), \quad\left(x^{\prime}, z\right) \in \Omega ,\\
	\operatorname{div}_{\mathrm{x}^{\prime}} \int_{0}^{1} u\left(x^{\prime}, \tau\right) \mathrm{d} \tau=0, \\
	u \text { ~is~ periodic  ~with ~ respect ~to~} x_{1}, x_{2}, \cdots, x_{n-1}, \\
	u\left(x^{\prime}, 0\right)=u\left(x^{\prime}, 1\right)=0,
\end{array}\right.
$$
	admits a unique solution $(u,P) \in (W^{2,p}_{per}(\Omega))^{n-1}\times W^{1,p}_{per}(\mathcal{O}) \cap L^{p}_{0} (\mathcal{O})$, where the Sobolev space $W_{per}^{2,p}(\Omega)$ defined as $W_{per}^{2,p}(\Omega)=\left\{u \in W^{2,p}(\Omega) ;  \partial^{l} u / \partial x_{i}^{l}, i=1,2, \cdots, n-1,  l=0,1,2,\right.\\\left. \text{ is}~\mathcal{O}\text{-periodic} \right\}.$ Moreover, there exists a constant $C(\Omega)>0$ such that
	\begin{equation}\label{377}
		\|u\|_{W^{2,p}(\Omega)}+\|P\|_{W^{1,p}(\mathcal{O})}\leq C(\Omega) \|f\|_{L^{p}(\Omega)}.
	\end{equation}
\end{Lemma}
{\bf Proof.}
The proof of this lemma is similar to the method of \cite{ziane} (see also \cite{hie}), and we omit the details here.
$\hfill\Box$

Next, we give the  Gagliardo-Nirenberg interpolation inequality
\begin{Lemma}
(see $\cite{nir}$)(Gagliardo-Nirenberg interpolation inequality). For a function $u: \Omega \rightarrow \mathbb{R}$ defined on a bounded Lipschitz domain $\Omega \subset \mathbb{R}^{n}, \forall~  1 \leq q, r \leq \infty$, and a natural number m. Suppose also that a real number $\theta$ and a natural number $j$ are such that
$$
\frac{1}{p}=\frac{j}{n}+\left(\frac{1}{r}-\frac{m}{n}\right) \theta+\frac{1-\theta}{q}
,$$
and
$$
\frac{j}{m} \leq \theta \leq 1.
$$
Then, we have
$$
\left\|D^{j} u\right\|_{L^{p}} \leq C_{1}\left\|D^{m} u\right\|_{L^{r}}^{\theta}\|u\|_{L^{q}}^{1-\theta}+C_{2}\|u\|_{L^{s}},
$$
where s is a positive constant. The constants $C_{1}$ and $C_{2}$ depend on the domain $\Omega$ and $m, n, j, r, q, \theta$.
\end{Lemma}

Finally, we state a Gronwall-type inequality which will be used to guarantee the local existence of solutions.

\begin{Lemma}\label{Lemma-2}(see \cite{simon})
	Let $g \in W^{1,1}(0,T)$ and $ k \in L^{1}(0,T)$ satisfy
	
	$$\frac{dg}{dt}\leq F(g)+k
	\quad in \quad [0,T],\quad
	g(0)\leq g_{0},$$
	where F is bounded on bounded sets from R into R, that is ,
	
	$$ \forall~ a>0 ~ \exists~ A>0 ~ such ~
	that~ |x|\leq a \rightarrow|F(x)|\leq A.$$
	Then for every $\epsilon>0$, there exists $T_{\epsilon}>0$ independent of g, but depend on
	$ g_{0}$ and $ k $ such that
	$$g(t)\leq g_{0}+\epsilon, ~ \forall~ t\leq T_{\epsilon}.$$
\end{Lemma}


\section{Well-posedness on the  regularization problems}
\ \ \ \

In this section, in order to prove the local-time existence and uniqueness of strong solutions to the nonhomogeneous primitive  equations, we consider the  following regularization of systems (\ref{a1}):
\begin{equation}\label{3.1}
	\left\{	\begin{aligned}
		&\rho_{t}+u\partial_{x}\rho-\int_{0}^{y}\partial_{x}u(s)ds\partial_{y}\rho-\lambda\rho_{xx}=0 ,\\
		&\rho u_{t}+\rho u\partial_{x}u-\rho\int_{0}^{y}\partial_{x}u(s)ds\partial_{y}u+\partial_{x}P-\partial_{x}(\mu(\rho)\partial_{x}u)-\partial_{y}(\mu(\rho)\partial_{y}u)=\rho f ,\\
		&\partial_{y}P=0, \quad \partial_{x}\int_{0}^{1}u(x,y)dy=0,
	\end{aligned}\right.
\end{equation}

We remark that we  get these estimates below which are independent on $\lambda$. Moreover, if we let $\lambda\rightarrow0$, we can get the well-posedness of the original problems (\ref{a1}) with some higher regularity with $\rho\geq\delta$.

\begin{Proposition}\label{pro-1}
Suppose that the initial data $(\rho_0,u_0)$ satisfies the regularity condition
\begin{equation}\label{1711}
	\delta \le \rho_{0}\in H^{3},\nabla \rho_{0}\in W^{1,\infty}, u_{0}\in H^{1}_{0,per}\cap H^{2}\cap H^{3}
\end{equation} and $f \in H^{2}((0,T)\times\Omega)$, and the compatibility condition (\ref{1.5}) for
 $(P_{0},f_{0},V_{1})\in H^{2}\times H^{1} \times H^{1}$. Then there exists a small time $T_{*}\in(0,T)$ and a unique strong solution $(\rho, u, w ,P)$ to the system $(\ref{3.1})$ with the initial conditions $(\ref{a2})-(\ref{14})$ such that
\begin{equation}
	\begin{aligned}
		&\rho \in L^{\infty}(0, T ; H^{3}) \cap C\left([0, T] ; W^{2, q}\right), \quad \rho_{t} \in L^{\infty}\left(0, T ; W^{1, 2}\right),
		\\	&u \in L^{\infty}\left([0, T] ;  H^{3}\right) \cap L^{2}\left(0, T ;  H^{4}\right),~ u_{t} \in L^{\infty}\left(0, T ; H_{0,per}^{1}\right)\cap L^{2}(0,T;H^{2}),\\& P\in L^{\infty}\left(0, T ; H^{1}\right) \cap L^{2}\left(0, T ; H^{2}\right),
	\end{aligned}
\end{equation}
where $1\leqslant q<\infty$.
\end{Proposition}
Our proof of the existence of the system (\ref{3.1}) is based on linearization and iteration method. Next, in the subsection we will show the existence of solutions to the linearized system of (\ref{3.1}) with a known vector field $v$ by the Galerkin method.
In this section,  the positive constant $\widetilde{C}$, which may be different from line to line, depends  on $\|\partial^{3}\mu/\partial \rho^{3}\|_{C}$, $\|V_{1}\|_{H^{1}}$, $\|\nabla^{3}\rho_{0}\|_{L^{2}}$, $\|\nabla ^{3}u_{0}\|_{L^{2}}$, $\|f\|_{H^{2}}$ and the parameters of C.
\subsection{Regularized results on  regularization problems $(\ref{3.1})$}
Our derivation of higher order a prior estimates is based on the following regularity results on the  hydrostatic Stokes system.
\begin{Lemma}\label{Lemma-41}
	Assume that $\nabla^{2}\rho\in W^{1,2}$ and $\nabla \rho \in W^{1,\infty}$. Let $(u,P)\in H^{1}_{0,per}\times L^{2}$ be the unique weak solution to the problem
	\begin{equation}\label{l33}
		-\partial_{x}(\mu(\rho)\partial_{x}u)-\partial_{y}(\mu(\rho)\partial_{y}u)+\partial_{x}P=f, ~ \partial_{y}P=0,~\int_{\Omega}P=0.
	\end{equation}
	Then it holds that\\
	$	(1)$If $f \in L^{2}$, then $(u, P) \in H^{2} \times H^{1}$ and
	\begin{equation}\label{3222}
		\|u\|_{H^{2}}+\|P\|_{H^{1}} \leqslant C\left(1+\|\nabla \rho\|_{L^{\infty}}\right
		)\|f\|_{L^{2}}.
	\end{equation}
	$(2)$If $f \in L^{r}$, then $(u, P) \in W^{2,r} \times W^{1,r}$ and\begin{equation}\label{3333}
		\| u\|_{W^{2,r}}+\| P\|_{W^{1,r}}\leq C(1+\|\nabla\rho\|_{L^{\infty}}+\|\nabla \rho\|_{L^{\infty}}^{2})\|f\|_{L^{r}} ,
	\end{equation}
	where $r>2$.\\
	$(3)$ If $f \in H^{1}$, then $(u, P) \in H^{3} \times H^{2}$ and
	\begin{equation}
		\|u\|_{H^{3}}+\|P\|_{H^{2}} \leqslant C(1+\|\nabla\rho\|_{L^{ \infty}}+\|\nabla\rho\|_{L^{ \infty}}^{2})(1+\|\nabla^{2}\rho\|_{L^{2}})\|f\|_{H^{1}}.
	\end{equation}
	$	(4)$ If $f \in H^{2}$, then $(u, P) \in H^{4} \times H^{3}$ and
	\begin{equation}
		\|u\|_{H^{4}}+\|P\|_{H^{3}} \leqslant \widetilde{C}(1+\|\nabla\rho\|_{L^{ \infty}}+\|\nabla\rho\|_{L^{ \infty}}^{2}+\|\nabla^{2}\rho\|_{L^{\infty}}+\|\nabla ^{3}\rho\|_{L^{2}})^{2}(1+\|\nabla^{2}\rho\|_{L^{2}})\|f\|_{H^{2}}.
	\end{equation}

\end{Lemma}
{\bf Proof.} The existence and uniqueness of the solution can be found in
\cite{ziane}. Here we give the a priori estimates. First, we will show that
\begin{equation}\label{32}
	\left\|\nabla u\right\|_{L_{2}}+\left\|P\right\|_{L{2}}\leq C\left\|f\right\|_{L{2}}.
\end{equation}

Multiplying ($\ref{l33}$) by $u$ and integrating over $\Omega$ it follows that $$\int_{\Omega} \mu(\rho)|\nabla u|^{2} \leq \int_{\Omega} fu \leq C \left\|f\right\|_{L_{2}}\left\|u\right\|_{L_{2}}\leq C \left\|f\right\|_{L_{2}}^{2}+\frac{1}{2C}\left\|\nabla u\right\|_{L_{2}}^{2},$$
and  since $C^{-1}\leq\mu\leq C$, we obtain the energy inequality $\left\|\nabla u\right\|_{L_{2}}\leq C\left\|f\right\|_{L_{2}}$.

According to Bovosgii's theory, we can find a function  $v=(v_{1},v_{2})\in H^{1}_{0}(\Omega)$ such that $P$=div $v$ and $\left\|v\right\|_{H^{1}}\leq C \left\|P\right\|_{L_{2}}$ (see $\cite{galdi}$), then
\begin{equation}\label{33333}
	\begin{aligned}
		&\int_{\Omega} P^{2}dydx=\int_{\Omega} P(\partial_{x}v_{1}+\partial_{y}v_{2})dydx=-\int_{\Omega}\partial_{x}P v_{1}dydx\\&
		=-\int_{\Omega} fv_{1} +\partial_{x}(\mu(\rho)\partial_{x}u)v_{1}+\partial_{y}(\mu(\rho)\partial_{y}u)v_{1}dydx\leq C\left\|f\right\|_{L_{2}}\left\|v\right\|_{H_{1}}\leq C\left\|f\right\|_{L_{2}}\left\|P\right\|_{L_{2}},
	\end{aligned}
\end{equation}
so that $\left\|P\right\|_{L^{2}}\leq C\left\|f\right\|_{L^{2}}$.

Next, We introduce $\tilde{P}=P/\mu(\rho)$  to rewrite the equation $(\ref{l33})$ as
$$-\Delta u+\partial_{x}\tilde{P}=\mu(\rho)^{-1}(f+ \nabla \mu(\rho) \nabla u+\tilde{P} \partial_{x} \mu(\rho)).$$

Then, it follows from the regularity results on the hydrostatic Stokes system (see Lemma \ref{Lemma-1}) that
\begin{equation}\label{34}
	\left\|u\right\|_{H^{2}}+\left\|\tilde{P}\right\|_{H^{1}}\leq C\left\|\mu(\rho)^{-1}(f+ \nabla \mu(\rho) \nabla u+\tilde{P} \partial_{x} \mu(\rho))\right\|_{L^{2}}\leq C (\left\|f\right\|_{L^{2}}+\left\|\nabla \rho \nabla u\right\|_{L^{2}}+\left\|\tilde{P}\partial_{x} \rho\right\|_{L^{2}}).
\end{equation}

Using the H$\ddot{o}$lder inequality and (\ref{32}), we have
\begin{equation}\label{35}
	\begin{aligned}
		\left\|\nabla \rho \nabla u\right\|_{L^{2}}+\left\|\tilde{P}\partial_{x} \rho\right\|_{L^{2}}
		&\leq\left\|\nabla\rho\right\|_{L^{\infty}}\left\|\nabla u\right\|_{L^{2}}+\left\|\partial_{x}\rho\right\|_{L^{\infty}}\left\| \tilde{P}\right\|_{L^{2}}\\
		&\leq C \left\|\nabla\rho\right\|_{L^{\infty}}\left\|f\right\|_{L^{2}}.
	\end{aligned}
\end{equation}

Substituting $(\ref{35})$ into $(\ref{34})$, we proved the first part $(1)$ of the lemma.

Applying Lemma \ref{the31} we deduce that
\begin{equation}
	\begin{aligned}
		&	\left\|u\right\|_{W^{2,r}}+\left\|\tilde{P}\right\|_{W^{1,r}}\\&\leq C\left\|\mu(\rho)^{-1}(f+ \nabla \mu(\rho) \nabla u+\tilde{P} \partial_{x} \mu(\rho))\right\|_{L^{r}}\\&\leq C (\left\|f\right\|_{L^{r}}+\left\|\nabla \rho \nabla u\right\|_{L^{r}}+\left\|\tilde{P}\partial_{x} \rho\right\|_{L^{r}})\\
		&\leq C(\left\|\nabla\rho\right\|_{L^{\infty}}\left\|\nabla u\right\|_{L^{r}}+\left\|\partial_{x}\rho\right\|_{L^{\infty}}\left\| \tilde{P}\right\|_{L^{r}}+\left\|f\right\|_{L^{r}})\\
		&\leq C(\left\|\nabla\rho\right\|_{L^{\infty}}\left\|\nabla u\right\|_{H^{1}}+\left\|\partial_{x}\rho\right\|_{L^{\infty}}\left\|\tilde {P}\right\|_{H^{1}}+\left\|f\right\|_{L^{r}})\\
		&\leq C [\|f\|_{L^{r}}+\|\nabla \rho\|_{L^{\infty}}(1+\|\nabla\rho\|_{L^{\infty}})\|f\|_{L^{2}}]\\
		&\leq C\|f\|_{L^{r}}(1+\|\nabla \rho\|_{\infty}+\|\nabla \rho\|_{L^{\infty}}^{2}).
	\end{aligned}
\end{equation}

The corresponding result $(3)$ can be proved similarly by recalling
\begin{equation}
	\left\|u\right\|_{H^{3}}+\left\|\tilde{P}\right\|_{H^{2}}\leq C\left\|\mu(\rho)^{-1}(f+ \nabla \mu(\rho) \nabla u+\tilde{P} \partial_{x} \mu(\rho))\right\|_{H^{1}},
\end{equation}
and using the previous estimates ($\ref{3222}$) and ($\ref{3333}$), we can get$$
\|u\|_{H^{3}}+\|P\|_{H^{2}} \leqslant C(1+\|\nabla\rho\|_{L^{ \infty}}+\|\nabla\rho\|_{L^{ \infty}}^{2})(1+\|\nabla^{2}\rho\|_{L^{2}})\|f\|_{H^{1}}.$$

In fact,
\begin{equation}
\begin{aligned}
&\|u\|_{H^{3}}+\|\tilde{P}\|_{H^{2}}\\&\leq C (\|\frac{f}{\mu(\rho)}\|_{H^{1}}+\|\frac{\nabla\mu(\rho)\nabla u}{\mu(\rho)}\|_{H^{1}}+\|\frac{P}{\mu(\rho)}\partial_{x}\mu(\rho)\|_{H^{1}})\\
&\leq C(\|f\|_{H^{1}}+\|\nabla \rho\nabla u\|_{L^{2}}+\|P \partial_{x}\rho\|_{L^{2}}\\
&\quad+\|\nabla^{2}\rho \nabla u\|_{L^{2}}+\|\nabla \rho \nabla^{2}u\|_{L^{2}}+\|\nabla P\partial_{x}\rho\|_{L^{2}}+\|P\nabla \partial_{x}\rho\|_{L^{2}})\\
&\leq C (\|f\|_{H^{1}}+\|\nabla \rho\|_{L^{\infty}}\|\nabla u\|_{L^{2}}+\|P\|_{L^{2}}\|\partial_{x}\rho\|_{L^{\infty}}+\|\nabla^{2}\rho\|_{L^{2}}\|\nabla u\|_{L^{\infty}}\\
&\quad+\|\nabla \rho\|_{L^{\infty}}\|\nabla^{2}u\|_{L^{2}}+\|\nabla P\|_{L^{2}}\|\partial_{x}\rho\|_{L^{\infty}}+\|P\|_{L^{\infty}}\|\nabla \partial_{x}\rho\|_{L^{2}})\\
&\leq C (\|f\|_{H^{1}}+\|\nabla \rho\|_{L^{\infty}}\|\nabla u\|_{L^{2}}+\|P\|_{L^{2}}\|\partial_{x}\rho\|_{L^{\infty}}+\|\nabla^{2}\rho\|_{L^{2}}\|\nabla u\|_{W^{1,r}}\\
&\quad+\|\nabla \rho\|_{L^{\infty}}\|\nabla^{2}u\|_{L^{2}}+\|\nabla P\|_{L^{2}}\|\partial_{x}\rho\|_{L^{\infty}}+\|P\|_{W^{1,r}}\|\nabla \partial_{x}\rho\|_{L^{2}})\\
&\leq C [\|f\|_{H^{1}}+\|\nabla \rho\|_{L^{\infty}}\|f\|_{L^{2}}+\|\partial_{x}\rho\|_{L^{\infty}}\|f\|_{L^{2}}\\
&\quad+\|\nabla^{2}\rho\|_{L^{2}}(1+\|\nabla \rho\|_{L^{\infty}}+\|\nabla \rho\|_{L^{\infty}}^{2})\|f\|_{L^{r}}+\|\nabla\rho\|_{L^{\infty}}(1+\|\nabla \rho\|_{L^{\infty}})\|f\|_{L^{2}}]\\
&\leq C\|f\|_{H^{1}}(1+\|\nabla \rho\|_{L^{\infty}}+\|\nabla \rho\|_{L^{\infty}}^{2})(\|\nabla^{2}\rho\|_{L^{2}}+1).
\end{aligned}
\end{equation}
Finally, $(4)$ can be proved similarly by recalling
\begin{equation}
	\begin{aligned}
		\|u\|_{H^{4}}+\|\tilde{P}\|_{H^{3}}&\leq \widetilde{C} (\|\frac{f}{\mu(\rho)}\|_{H^{2}}+\|\frac{\nabla\mu(\rho)\nabla u}{\mu(\rho)}\|_{H^{2}}+\|\frac{P}{\mu(\rho)}\partial_{x}\mu(\rho)\|_{H^{2}})\\
	\end{aligned}
\end{equation}and using the previous estimates $(1)$, $(2)$ and $(3)$.
We complete the proof of the Lemma \ref{Lemma-41}.
$\hfill\Box$

 Next, we state some regularity estimates on density. Due to the loss of derivation of velocity,  this estimates on density require more regularity of horizontal velocity  than  ones of the Navier-Stokes equations. 
\begin{Lemma} \label{Lemma-2.3}
	Assume that $\nabla \rho_{0}\in W^{1,\infty}(\Omega)$ and $\nabla^{2}\rho_{0}\in W^{1,2}(\Omega)$, let $\rho$ be the unique solution to the problem \begin{equation}\label{2.2}
		\partial_{t}\rho +\partial_{x}\rho v-\partial_{y}\rho\int_{0}^{y}\partial_{x}v(s)ds-\lambda\partial_{ xx}\rho=0,~\rho|_{t=0}=\rho_{0},
	\end{equation} where $\lambda>0$ is a small number. Then we have the following regularity results:

	(1)If $v\in L^{1}([0,T],H^{3}(\Omega))$, then
\begin{equation}\label{2.4}
\|\nabla \rho(t)\|_{L^{\infty}}\leq \|\nabla \rho_{0}\|_{L^{\infty}}\exp(C\int_{0}^{t}\|v(s)\|_{H^{3}}ds),
\end{equation}
	and
\begin{equation}\label{2.5}
\|\nabla ^{2}\rho(t)\|_{L^{2}}\leq (\|\nabla \rho_{0}\|_{L^{\infty}}+\|\nabla ^{2}\rho_{0}\|_{L^{2}})\exp(C\int_{0}^{t}\|
v(s)\|_{H^{3}}ds).
\end{equation}

(2)If $v\in L^{1}([0,T],H^{4}(\Omega))$, then
	\begin{equation}\label{2.6}
		\| \nabla ^{2}\rho(t)\|_{L^{\infty}}\leq \|\nabla^{2}\rho_{0}\|_{L^{\infty}}\exp(C\int_{0}^{t}\|v(s)\|_{H^{4}}ds),
	\end{equation}
	and
	\begin{equation}\label{2.7}
		\| \nabla ^{3}\rho(t)\|_{L^{2}}\leq(\| \nabla ^{3}\rho_{0}\|_{L^{2}}+\|\nabla^{2}\rho_{0}\|_{L^{\infty}}+\|\nabla\rho_{0}\|_{L^{\infty}}) \exp(C\int_{0}^{t}\|v(s)\|_{H^{4}}ds).
	\end{equation}	
\end{Lemma}
{\bf Proof.} Taking the gradient operator $\nabla$ to $(\ref{2.2})$ and multiplying by $p|\nabla\rho|^{p-2}\nabla\rho$, we have
\begin{equation}
\begin{aligned}	&(|\nabla\rho|^{p})_{t}+v\cdot\partial_{x}|\nabla\rho|^{p}-\int_{0}^{y}\partial_{x}v(s)ds\partial_{y}|\nabla\rho|^{p}+p\partial_{x}\rho\nabla v|\nabla\rho|^{p-2}\nabla\rho+p\partial_{y}\rho U |\nabla\rho|^{p-2}\nabla\rho\\&-\lambda\nabla\partial_{xx}\rho p|\nabla \rho|^{p-2}\nabla \rho=0,
	\end{aligned}
\end{equation}
where $U=(-\int_{0}^{y}\partial_{xx}v(s)ds, -\partial_{x}v(x))$.
Integrating over $\Omega$, we obtain
\begin{equation*}
	\begin{aligned} \frac{d}{dt}\int_{\Omega}|\nabla\rho|^{p}+2\lambda\int_{\Omega}p|\nabla \partial_{x}\rho|^{2}|\nabla \rho|^{p-2}&\leq p\int_{\Omega}|\nabla \rho|^{p}|\nabla v|+p\int_{\Omega}|\nabla \rho|^{p}|\int_{0}^{y}\partial_{xx}v(s)ds|\\
		&\leq p\left\|\nabla v\right\|_{L^{\infty}}\left\|\nabla\rho\right\|_{L^{p}}^{p}+p\left\|\int_{0}^{y}\partial_{xx}v(s)ds\right\|_{L^{\infty}_{x}L^{\infty}_{y}}
		\left\|\nabla\rho\right\|_{L^{p}}^{p}.\\
	\end{aligned}
\end{equation*}
Then after applying Gronwall's inequality and Sobolev inequality, we derive that
$$\|\nabla \rho(t)\|_{L^{p}} \leqslant\left\|\nabla \rho_{0}\right\|_{L^{p}} \exp \left(C \int_{0}^{t}\|v(s)\|_{H^{3}} \mathrm{~d} s\right),$$
where $p\geq1$,
which implies $(\ref{2.4})$.

To prove $(\ref{2.5})$, we first give the estimate $\|\partial_{x x}\rho\|_{L^{2}}$. The other terms can be dealt with in a similar way. Taking $\partial_{x x}$ on both sides of $(\ref{2.2})$ then multiplying by $\partial_{x x}\rho$ and integrating over $\Omega$, we have
\begin{equation*}\label{315} \frac{d}{dt}\int_{\Omega}|\partial_{xx}\rho|^{2}+\int_{\Omega}[\partial_{x x}v\partial_{x}\rho\partial_{x x}\rho-\int_{0}^{y}\partial_{x xx}v(s)ds\partial_{y}\rho\partial_{x x}\rho+2\partial_{x}v|\partial_{x x}\rho|^{2}-2\int_{0}^{y}\partial_{x x}v(s)ds\partial_{x x}\rho\partial_{x y}\rho]\end{equation*}$$+\int_{\Omega}\lambda|\partial_{x xx}\rho|^{2}=0.$$\\
By H$\ddot{o}$lder inequality and Young inequality, we have
\begin{equation*}\label{2.9}
\begin{aligned}
&\frac{d}{dt}\|\partial_{x x}\rho\|_{L^{2}}^{2}+\int_{\Omega}\lambda|\partial_{x xx}\rho|^{2}\\
&\leq \|\partial_{x x}v\|_{L^{2}}\|\partial_{x}\rho\|_{L^{\infty}}\|\partial_{x x}\rho\|_{L^{2}}
+\|\int_{0}^{y}\partial_{x xx}v(s)ds\|_{L^{2}}\|\partial_{y}\rho\|_{L^{\infty}}\|\partial_{x x}\rho\|_{L^{2}}\\
&\quad+\|\partial_{x }v\|_{L^{\infty}}\|\partial_{x x}\rho\|_{L^{2}}^{2}+\|\int_{0}^{y}\partial_{x x}v(s)ds\|_{L^{\infty }}\|\partial_{x x}\rho\|_{L^{2}}\|\partial_{x y}\rho\|_{L^{2}}\\
&\leq\|\partial_{x}\rho\|_{L^{\infty}}\|v\|_{H^{2}}\|\partial_{x x}\rho\|_{L^{2}}+\|\partial_{x x}\rho\|_{L^{2}}\|\partial_{x xx}v\|_{L^{2}}\|\partial_{y} \rho\|_{L^{\infty}}+ \|\partial_{xx}\rho\|_{L^{2}}^{2}\|v\|_{H^{3}}\\
&\quad+\|v\|_{H^{3}}\|\partial_{x x}\rho\|_{L^{2}}\|\partial_{x y}\rho\|_{L^{2}}.
\end{aligned}
\end{equation*}
Therefore, we have\begin{equation*}2\frac{d}{dt}\|\partial_{x x}\rho\|_{L^{2}}\leq\|v\|_{H^{3}}\|\nabla \rho \|_{L^{\infty}}+C\|v\|_{H^{3}}(\|\partial_{x x}\rho\|_{L^{2}}+\|\partial_{x y}\rho\|_{L^{2}}).\end{equation*}
Similar inequality also hold for $\partial_{x y}$$\rho$ and  $\partial_{y y}\rho$. Hence,
\begin{equation*}\label{2.10}
	2\frac{d}{dt}\|\partial_{x y}\rho\|_{L^{2}}
	\leq
	\|v\|_{H^{2}}\|\nabla \rho \|_{L^{\infty}}+C\|v\|_{H^{3}}(\|\partial_{x x}\rho\|_{L^{2}}+\|\partial_{x y}\rho\|_{L^{2}}+\|\partial_{y y}\rho\|_{L^{2}})
\end{equation*}
and \begin{equation*}\label{2.11}
	2\frac{d}{dt}\|\partial_{y y}\rho\|_{L^{2}}
	\leq
	\|v\|_{H^{2}}\|\nabla \rho \|_{L^{\infty}}+C\|v\|_{H^{3}}(\|\partial_{yy}\rho\|_{L^{2}}+\|\partial_{x y}\rho\|_{L^{2}}).
\end{equation*}
Summing up the above inequalities and using the Gronwall's inequality, one gets
\begin{equation*}\label{317}
	\begin{aligned}
		\|\nabla^{2}\rho\|_{L^{2}}&\leq[\|\nabla^{2} \rho_{0}\|_{L^{2}}+\int_{0}^{t}
		\|\nabla \rho_{0}\|_{L^{\infty}}\exp(\int_{0}^{\tau}C\|v(s)\|_{H^{3}}ds)\|v\|_{H^{2}}d\tau] exp(\int_{0}^{t}C\|v(\tau)\|_{H^{3}}d\tau)\\
		&\leq (\|\nabla^{2} \rho_{0}\|_{L^{2}}+\|\nabla\rho_{0}\|_{L^{\infty}})\exp( C{\int_{0}^{t}\|v(s)\|_{H^{3}}}ds).
	\end{aligned}
\end{equation*}
Next, we shall perform some estimates to show that (3.20).
Taking $\nabla \partial_{x}$ on both sides of $(\ref{2.2})$ 
and multiplying by $ p|\nabla \partial_{x}\rho|^{p-2}\nabla \partial_{x}\rho$, we have
\begin{equation*}
\begin{aligned}
&\partial_{t}|\nabla \partial_{x }\rho|^{p}+\nabla \partial_{x }v\partial_{x }\rho p|\nabla\partial_{x }\rho|^{p-2}\nabla\partial_{x }\rho+ \partial_{x }v\nabla\partial_{x }\rho p|\nabla\partial_{x }\rho|^{p-2}\nabla\partial_{x }\rho\\
&+\nabla v\partial_{x x}\rho p|\nabla\partial_{x }\rho|^{p-2}\nabla\partial_{x }\rho+ v\nabla\partial_{x x}\rho p|\nabla\partial_{x }\rho|^{p-2}\nabla\partial_{x }\rho-\nabla \int_{0}^{y}\partial_{xx }vds\partial_{y }\rho p|\nabla \partial_{x }\rho|^{p-2}\nabla \partial_{x }\rho\\
&- \int_{0}^{y} \partial_{xx }vds\nabla\partial_{y }\rho p|\nabla \partial_{x }\rho|^{p-2}\nabla \partial_{x }\rho-\nabla\int_{0}^{y}\partial_{x } vds\partial_{xy }\rho p|\nabla \partial_{x }\rho|^{p-2}\nabla \partial_{x }\rho\\&-\int_{0}^{y} \partial_{x }vds\nabla\partial_{xy }\rho p|\nabla \partial_{x }\rho|^{p-2}\nabla \partial_{x }\rho-\lambda\nabla \partial_{x xx}\rho p|\nabla \partial_{x}\rho|^{p-2}\nabla\partial_{x } \rho=0.
	\end{aligned}
\end{equation*}
Integrating over $\Omega$, we have
\begin{equation*}
\begin{aligned}
&\frac{d}{dt}\|\nabla \partial_{x}\rho\|_{L^{p}}^{p}+\int_{\Omega}\lambda|\nabla\partial_{x x}\rho|^{2}p(p-1)|\nabla\partial_{x}\rho|^{p-2}\\
&\leq\|\nabla \partial_{x }v\|_{L^{\infty}}\|\partial_{x}\rho\|_{L^{p}}p
\|\nabla \partial_{x}\rho\|_{L^{p}}^{p-1}+\|\partial_{x }v\|_{L^{\infty}}p\|\nabla \partial_{x}\rho\|_{L^{p}}^{p}+\|\nabla v\|_{L^{\infty}}p\|\nabla \partial_{x}\rho\|_{L^{p}}^{p}\\
&\quad+\|v\|_{H^{4}}p\|\nabla \partial_{x}\rho\|_{L^{p}}^{p-1}\|\partial_{y}\rho\|_{L^{p}}+\|v\|_{H^{3}}p\|\nabla \partial_{x}\rho\|_{L^{p}}^{p-1}\|\nabla\partial_{y}\rho\|_{L^{p}}
+\|v\|_{H^{3}}p\|\nabla\partial_{x}\rho\|_{L^{p}}^{p},
\end{aligned}
\end{equation*}
which implies that
 \begin{equation*}
	\begin{aligned}
		\frac{d}{dt}\|\nabla \partial_{x}\rho\|_{L^{p}}	\leq
		\|v\|_{H^{4}}\|\nabla \rho \|_{L^{p}}+C\|v\|_{H^{3}}(\|\nabla \partial_{y}\rho\|_{L^{p}}+\|\nabla \partial_{x}\rho\|_{L^{p}}).
	\end{aligned}
\end{equation*}
Similarly, one has 
\begin{equation*}
\begin{aligned}
\frac{d}{dt}\|\nabla \partial_{y}\rho\|_{L^{p}}	\leq
\|v\|_{H^{4}}\|\nabla \rho \|_{L^{p}}+C\|v\|_{H^{3}}(\|\nabla \partial_{y}\rho\|_{L^{p}}+\|\nabla \partial_{x}\rho\|_{L^{p}}).
\end{aligned}
\end{equation*}
Therefore
\begin{equation*} \frac{d}{dt}\|\nabla^{2}\rho\|_{L^{p}}\leq C\|v\|_{H^{3}}\|\nabla^{2}\rho\|_{L^{p}}+\|v\|_{H^{4}}\|\nabla \rho\|_{L^{p}}.
\end{equation*}
 Applying the Gronwall's inequality yields that
\begin{equation*}\begin{aligned}
		\|\nabla^{2}\rho\|_{L^{p}}&\leq\exp{\int_{0}^{t}C\|v(s)\|_{H^{3}}}ds(\|\nabla ^{2}\rho_{0}\|_{L^{p}}+\int_{0}^{t}\|v\|_{H^{4}}(s)ds\|\nabla\rho\|_{L^{p}})\\&\leq(\|\nabla^{2}\rho_{0}\|_{L^{p}}+\|\nabla\rho_{0}\|_{L^{p}})\exp(\int_{0}^{t}C\|v\|_{H^{4}}(s)ds).
	\end{aligned}
\end{equation*}
Let $p\rightarrow\infty$, one has$$\|\nabla^{2}\rho\|_{L^{\infty}}\leq(\|\nabla^{2}\rho_{0}\|_{L^{\infty}}+\|\nabla\rho_{0}\|_{L^{\infty}})\exp(\int_{0}^{t}C\|v\|_{H^{4}}(s)ds).$$
Similar arguments also hold for the higher order derivatives.
Therefore, direct estimates imply that
\begin{equation*}
\begin{aligned}	\frac{d}{dt}\|\partial_{xxx}\rho\|_{L^{2}}\leq\|v\|_{H^{4}}(\|\nabla ^{2}\rho\|_{L^{\infty}}+\|\nabla \rho\|_{L^{\infty}})+C\|v\|_{H^{3}}(\|\partial_{xxx}\rho\|_{L^{2}}+\|\partial_{xxy}\rho\|_{L^{2}}),
	\end{aligned}
\end{equation*}
\begin{equation*}
\begin{aligned} \frac{d}{dt}\|\partial_{xxy}\rho\|_{L^{2}}\leq\|v\|_{H^{3}}(\|\nabla^{2}\rho\|_{L^{\infty}}+\|\nabla\rho\|_{L^{\infty}})+C\|v\|_{H^{3}}(\|\partial_{xxx}\rho\|_{L^{2}}+\|\partial_{xyy}\rho\|_{L^{2}}+\|\partial_{xxy}\rho\|_{L^{2}}),
\end{aligned}
\end{equation*}
\begin{equation*}
\begin{aligned} \frac{d}{dt}\|\partial_{xyy}\rho\|_{L^{2}}\leq\|v\|_{H^{3}}(\|\nabla^{2}\rho\|_{L^{\infty}}+\|\nabla\rho\|_{L^{\infty}})+C\|v\|_{H^{3}}(\|\partial_{yyy}\rho\|_{L^{2}}+\|\partial_{xyy}\rho\|_{L^{2}}+\|\partial_{xxy}\rho\|_{L^{2}}),
\end{aligned}
\end{equation*}
and
\begin{equation*}\label{}	
\begin{aligned}
\frac{d}{dt}\| \partial_{yyy}\rho\|_{L^{2}}&\leq\|v\|_{H^{3}}(\|\nabla^{2}\rho\|_{L^{\infty}}+\|\nabla\rho\|_{L^{\infty}})+C\|v\|_{H^{3}}(\|\partial_{yyy}\rho\|_{L^{2}}+\|\partial_{xyy}\rho\|_{L^{2}}).
\end{aligned}
\end{equation*}
Combining the above estimates, one has
 $$\frac{d}{dt}\|\nabla^{3}\rho\|_{L^{2}}\leq\|v\|_{H^{4}}(\|\nabla^{2}\rho\|_{L^{\infty}}+\|\nabla\rho\|_{L^{\infty}})+C\|v\|_{H^{3}}\|\nabla ^{3}\rho\|_{L^{2}}.$$
This will imply after applying the Gronwall's inequality
$$\|\nabla^{3}\rho\|_{L^{2}}\leq(\|\nabla^{3}\rho_{0}\|_{L^{2}}+\|\nabla^{2}\rho_{0}\|_{L^{\infty}}+\|\nabla \rho_{0}\|_{L^{\infty}})\exp(\int_{0}^{t}C\|v\|_{H^{4}}ds).$$
Here, we complete the proof of this lemma.$\hfill\Box$
\subsection{Well-posedness on linearized equations }
In this subsection, we will consider the following linearized problem of \eqref{3.1} in $(0,T)\times\Omega$: \begin{equation}\label{37}
\left\{	\begin{aligned}
		&\rho_{t}+v\partial_{x}\rho-\int_{0}^{y}\partial_{x}v(s)ds\partial_{y}\rho-\lambda\rho_{xx}=0 ,\\
		&\rho u_{t}+\rho v\partial_{x}u-\rho\int_{0}^{y}\partial_{x}v(s)ds\partial_{y}u+\partial_{x}P-\partial_{x}(\mu(\rho)\partial_{x}u)-\partial_{y}(\mu(\rho)\partial_{y}u)=\rho f ,\\
		&\partial_{y}P=0, \quad \partial_{x}\int_{0}^{1}u(x,y)dy=0,
	\end{aligned}\right.
\end{equation}
where $v$ is a known divergence-free vector field and $\lambda$ is a small parameter. The initial and boundary conditions are formulated as
\begin{equation}\label{3122}
\begin{aligned}
		&\rho,~ u, ~ P~\text {are~periodic~ in~ the ~ direction}~ x,\\
	& u|_{y=0}=u|_{y=1}=0,\\
& u|_{t=0}=u_{0}(x,y),\\
			&\rho|_{t=0}=\rho_{0}(x,y).\\	
\end{aligned}
\end{equation}  We will prove the global existence and uniqueness of strong solutions to $(\ref{37})$-$(\ref{3122})$ for the positive initial data so that we can construct a sequence of the approximate solutions to the initial boundary value problem $(\ref{3.1})$ by iteration method.

Now we give the following lemma:
\begin{Lemma}
Assume that the hypotheses of Proposition $\ref{pro-1}$ are satisfied by ($\rho_{0},u_{0}, f)$. If $v_{t} \in L^{\infty}\left(0, T ;  H_{0,per}^{1}\right)\cap  L^{2}(0,T;H^{2})$, $v \in L^{\infty}\left(0, T ; H_{0,per}^{1}
	\cap H^{2}\cap H^{3}\right) \cap L^{2}\left(0, T ;  H^{4}\right)$, $v|_{t=0}=u_{0}$ and $ \partial_{x}\int_{0}^{1}v(x,y)dy=0$ in $\Omega$, then for any $T >0$, there exists a unique strong solution $(\rho, u, P)$ to the initial boundary value problem $(\ref{37})$-$(\ref{3122})$ such that \begin{equation}\label{311}
\begin{aligned}
&\rho \in L^{\infty}(0, T ; W^{3, 2}) \cap L^{\infty}(0, T ; W^{2, \infty})\cap C\left([0, T] ; W^{2, q}\right), \quad \rho_{t} \in L^{\infty}\left(0, T ; W^{1, 2}\cap L^{\infty}\right), \\
&u \in C\left([0, T] ; H_{0,per}^{1}\cap H^{2} \cap H^{3}\right) \cap L^{2}\left(0, T ; H^{4}\right),~u_{t}\in L^{\infty}(0,T;H^{1})\cap L^{2}(0,T;H^{2}), \\
&\sqrt{\rho} u_{t} \in L^{\infty}\left(0, T ; L^{2}\right), \quad P\in L^{\infty}\left(0, T ; H^{1}\right) \cap L^{2}\left(0, T ; H^{2}\right),
\end{aligned}
\end{equation}
where $1\leq q< \infty$.
\end{Lemma}

{\bf Proof.} Equations $(\ref{37})_{1}$ and
$(\ref{37})_{2}$ are  linear parabolic equations and can be solved by standard Galerkin approximation.  The existence and regularity of solutions to $(\ref{37})_{1}$ follow with similar arguments as $(\ref{37})_{2}$ below.  Now we only give some estimates. Together with  Lemma \ref{Lemma-2.3},  $(\ref{37})_1$ and Aubin Lions lemma, we deduce that $\rho\in L^{\infty}(0,T;W^{3,2}\cap W^{2,\infty})$,  $\rho_{t}\in L^{\infty}(0,T;W^{1,2}\cap L^{\infty})$ and $\rho \in C([0,T];W^{2,q})$ for any $1\leq q<\infty$. 


Now we  turn to the non-stationary hydrostatic Stokes equations. Since the uniqueness of strong solutions can be easily proved, we only prove the existence. First, we construct the approximate solutions by the semi-discrete-Galerkin method. To do this, we define a finite-dimensional space $X_{n}$=$span\big\{\psi_{1},\cdots,\psi_{n}\big\}$, $n\in N$, where the function $\psi_{n}$ is the $n$-$th$ eigenfunction to the following eigenvalue problem of the Dirichlet-Periodic problem for the hydrostatic Stokes equations in $\Omega$:

\begin{equation}\left\{\begin{aligned}
		&-\partial_{x x} \psi_{n}-\partial_{yy} \psi_{n}+\partial_{x} P_{n}=\lambda_{n} \psi_{n}, \\
		&\partial_{x} \int_{0}^{1} \psi_{n}(x, y) \mathrm{d} y=0, \quad \partial_{y} P_{n}=0,\\
		&\left.\psi_{n}\right|_{y=0}=\left.\psi_{n}\right|_{y=1}=0,\\
		&\psi_{n}~\text{ is~ periodic~ in~ the~ direction}~ x .
	\end{aligned}\right.
\end{equation}

Here $\lambda_{n}$ is the $n$-$th$ eigenvalue of the hydrostatic Stokes operator. The sequence $\big\{\psi_{i}\big\}_{i=1}^{\infty}$ can be renormalized in such a way that it is an orthonormal basis of $L^{2}$, which is also an orthogonal basis of $H^{1}_{0,per}(\Omega)=\Big\{\psi\in H^{1}(\Omega);\psi|_{y=0,1}=0, \psi(0,y)=\psi(L,y)\Big\}$ and a basis in $H^{1}_{0,per}(\Omega)\cap H^{2}$=$X$.
Introduce the function spaces
$$X=H^{1}_{0,per}(\Omega)\cap H^{2}=\Big\{\psi\in H^{1}\cap H^{2}(\Omega); \psi|_{y=0,1}=0, \psi(0,y)=\psi(L,y)\Big\}$$
and $X^{m}=span\big\{\psi^{1},\psi^{2},\cdots,\psi^{m}\big\},$ $ \psi=\Sum_{m=1}^{\infty}\left(\psi, \psi^{m}\right)_{L^{2}} \psi^{m}$  in $ H^{2}$, where we denote $\left(\psi, \psi^{m}\right)_{L^{2}}=\int_{\Omega} \psi \psi^{m} $.


 Then applying the standard Galerkin method, we construct an approximate solutions $u^{m}\in C^{1}([0,T];X^{m})$ such that for all $\phi\in X^{m}$,
\begin{equation}\label{314}
\int_{\Omega}[\rho u^{m}_{t}+\rho v\partial_{x}u^{m}-\rho\int_{0}^{y}\partial_{x}v(s)ds\partial_{y}u^{m}-\partial_{x}(\mu(\rho)\partial_{x}u^{m})
-\partial_{y}(\mu(\rho)\partial_{y}u^{m})]\cdot\phi=\int_{\Omega}\rho f\cdot \phi,
\end{equation}
\begin{equation}
u^{m}(0)=\sum_{k=1}^{m}(u_{0},\psi^{k})_{L^{2}}\psi^{k}.
\end{equation}

In fact, since $\rho\geq\delta>0$ in $\Omega $ and $\rho\in W ^{1,\infty}(0,T ;W^{1,2})$, we can rewrite $(\ref{314})$ with $\phi=\psi^{k} (k=1,\cdots,m)$ as a system of linear ODEs with regular coefficients. From the theory of linear ODEs, we obtain  a unique Galerkin-solution $u^{m}$ for every fixed m.
Next, we will show that a subsequence of the approximate solutions $u^{m}$ converges to a solution of the original  problem (\ref{37}). To show this, we need to derive some uniform bounds.

Note that $$C^{-1} \leqslant \mu \leqslant C,$$$$\left|\nabla \mu\right| \leqslant C|\nabla \rho|,  ~\left|\nabla^{2} \mu\right| \leqslant C\left(|\nabla \rho|^{2}+\left|\nabla^{2} \rho\right|\right)~\text { and }~|\nabla^{3}\mu|\leq \widetilde{C}(|\nabla\rho|^{3}+|\nabla \rho||\nabla^{2}\rho|+|\nabla^{3}\rho|),
$$
on $[0,T]\times\bar \Omega$.

Taking $\phi=u^{m}_{t}$ in $(\ref{314})$ and integrating by parts, we obtain
\begin{equation}
\begin{aligned}
&\int_{\Omega} \rho\left|u_{t}^{m}\right|^{2}+\frac{1}{2}\frac{\mathrm{d}}{\mathrm{d} t} \int_{\Omega} \mu\left|\nabla u^{m}\right|^{2}
\\&=  \int_{\Omega}\rho f u_{t}^{m}-\int \rho v \partial_{x}u^{m} u_{t}^{m}
+\int_{\Omega} \rho \int_{0}^{y} \partial_{x}v(s)d s \partial_{y} u^{m} u_{t}^{m}+\frac{1}{2} \int_{\Omega} \mu^{\prime}\left(\rho\right) \partial_{t} \rho\left|\nabla u^{m}\right|^{2}.
\end{aligned}
\end{equation}

Then H$\ddot{o}$lder inequality, the linear transport equation $(\ref{37})_1$,  (\ref{2.4}) and   (\ref{2.6}) yield
\begin{equation}
\begin{aligned}
&\int_{\Omega} \rho\left|u_{t}^{m}\right|^{2} +\frac{1}{2}\frac{\mathrm{d}}{\mathrm{d} t} \int_{\Omega} \mu\left|\nabla u^{m}\right|^{2} \\
&\leq\left\|\sqrt{\rho}u^{m}_{t}\right\|_{L^{2}}\left\|\sqrt{\rho}\right\|_{L^{\infty}}\left\|f\right\|_{L^{2}}
+\left\|\partial_{x}u^{m}\right\|_{L^{2}}\left\|\sqrt{\rho}u^{m}_{t}\right\|_{L^{2}}\left\|\sqrt{\rho}\right\|_{L^{\infty}}\left\|v\right\|_{L^{\infty}}\\
&\quad+\left\|\partial_{y}u^{m}\right\|_{L^{2}}\left\|\sqrt{\rho}u^{m}_{t}\right\|_{L^{2}}\left\|\sqrt{\rho}
\right\|_{L^{\infty}}\left\|\int_{0}^{y}\partial_{x}v(s)ds\right\|_{L^{\infty}}+C\left\|v\right\|_{L^{\infty}}\left\|\partial_{x}\rho\right\|_{L^{\infty}}\left\|\nabla u^{m}\right\|_{L^{2}}\\
&\quad+C\left\|\int_{0}^{y}\partial_{x}v(x,s)ds\right\|_{L^{\infty}}\left\|\partial_{y}\rho\right\|_{L^{\infty}}\left\|\nabla u^{m}\right\|_{L^{2}}+\lambda C\|\rho_{xx}\|_{L^{\infty}}\|\nabla u^{m}\|_{L^{2}}\\
&\leq C +\epsilon\left\|\sqrt{\rho}u^{m}_{t}\right\|_{L^{2}}^{2}+C\left\|\nabla u^{m}\right\|_{L^{2}}^{2},
\end{aligned}
\end{equation}
where we have used the following estimate $$\left\|\int_{0}^{y}\partial_{x}v(s)ds\right\|_{L^{\infty}}=\left\|\int_{0}^{y}\partial_{x}v(s)ds\right\|_{W^{1,2}_{x}W^{1,2}_{y}}
\leq\left\|\partial_{x}v\right\|_{W^{1,2}_{x}L^{2}_{y}}\leq\left\|v\right\|_{H^{2}}.$$

Therefore, applying Gronwall's inequality, we obtain
\begin{equation}\label{318}
\int_{0}^{T}\left\|\sqrt{\rho} u_{t}^{m}(t)\right\|_{L^{2}}^{2} \mathrm{~d} t+\sup _{0 \leqslant t \leqslant T}\left\|\nabla u^{m}(t)\right\|_{L^{2}}^{2} \leqslant C.
\end{equation}

Next, we differentiate $(\ref{314})$ with respect to $t$ and take $\phi=u^{m}_{t}$.
We can deduce
\begin{equation}\label{319}
\begin{aligned}
&\frac{1}{2}\frac{\mathrm{d}}{\mathrm{d} t}\int_{\Omega}\rho\left|u_{t}^{m}\right|^{2}+\int_{\Omega}\mu(\rho)\left|\nabla u_{t}^{m}\right|^{2}\\
&=-\frac{1}{2}\int_{\Omega}\rho_{t}(u^{m}_{t})^{2}-\int_{\Omega}\rho_{t}v\cdot\partial_{x} u^{m}u^{m}_{t}+\int_{\Omega}\rho_{t}\int_{0}^{y}\partial_{x}v(s)ds\cdot\partial_{y} u^{m}u^{m}_{t}\\
&\quad-\int_{\Omega}\mu_{t}(\rho)\nabla u^{m} \nabla u_{t}^{m}-\int_{\Omega}\rho v_{t}\partial_{x} u^{m} u^{m}_{t}-\int_{\Omega}\rho v\partial_{x}u_{t}^{m}u^{m}_{t}+\int_{\Omega}\rho\int_{0}^{y}\partial_{x}v_{t}(s)ds\partial_{y}u^{m}u^{m}_{t}\\
&\quad+\int_{\Omega}\rho\int_{0}^{y}\partial_{x}v(s)ds\partial_{y}u_{t}^{m}u^{m}_{t}+\int_{\Omega}\rho f_{t}u^{m}_{t}+\int_{\Omega}\rho_{t}fu^{m}_{t}=\sum_{i=1}^{10}\mathbb{I}_{i}.
\end{aligned}
\end{equation}

Then using $(\ref{37})_1$, we can deduce
\begin{equation}
\begin{aligned}
\mathbb{I}_{1}&=-\frac{1}{2}\int_{\Omega}(v\partial_{x}\rho-\int_{0}^{y}\partial_{x}v(s)ds\partial_{y}\rho-\lambda\rho_{xx}u^{m}_{t}u^{m}_{t}\\
&=\frac{1}{2}\int_{\Omega}( v \rho\partial_{x}u^{m}_{t} u^{m}_{t}-\int_{0}^{y}\partial _
{x}v(s)ds\rho\partial_{y}u^{m}_{t}u^{m}_{t}-\lambda C_{\delta} \rho_{xx}u^{m}_{t}u_{t}^{m})\\
&\leq \left\|v\right\|_{H^{2}}\left\|\sqrt{\rho}u^{m}_{t}\right\|_{L^{2}}\left\|\nabla u^{m}_{t}\right\|_{L^{2}}\left\|\sqrt{\rho}\right\|_{L_{\infty}}+\lambda C_{\delta}\|\rho_{xx}\|_{L^{\infty}}\|\sqrt{\rho}u^{m}_{t}\|_{L^{2}}\left\|\nabla u^{m}_{t}\right\|_{L^{2}}\\
&\leq C_{\delta}\left\|\sqrt{\rho}u^{m}_{t}\right\|_{L^{2}}^{2}+\epsilon\left\|\nabla u^{m}_{t}\right\|_{L_{2}}^{2}
\end{aligned}
\end{equation}where using Pioncare inequality in y-direction
and
\begin{equation}
\begin{aligned}
\mathbb{I}_{2}&=-\int_{\Omega}(v\partial_{x}\rho-\int_{0}^{y}\partial_{x}v(s)ds\partial_{y}\rho-\lambda\rho_{xx})v \partial_{x}u^{m}u^{m}_{t}\\
&\leq C(\left\|v\right\|_{L^{\infty}}\left\|\partial_{x}\rho\right\|_{L^{\infty}}\left\|\partial_{x}u^{m}\right\|_{L^{2}}
\left\|u^{m}_{t}\right\|_{L^{\infty}_{y}L^{2}_{x}}\left\|v\right\|_{L^{2}_{y}L^{\infty}_{x}}\\
&\quad+\left\|v\right\|_{L^{\infty}}
\left\|\partial_{y}\rho\right\|_{L^{\infty}}\left\|\partial_{x}u^{m}\right\|_{L^{2}}
\left\|u^{m}_{t}\right\|_{L^{2}_{y}L^{\infty}_{x}}\left\|\int_{0}^{y}\partial_{x}v(s)ds\right\|_{L^{2}_{x}L^{\infty}_{y}}\\
&\quad+\|\rho_{xx}\|_{L^{\infty}}\|\partial_{x}u^{m}\|_{L^{2}}\left\|u^{m}_{t}\right\|_{L^{\infty}_{y}L^{2}_{x}}\left\|v\right\|_{L^{2}_{y}L^{\infty}_{x}})\\
&\leq C+\epsilon\left\|\nabla u^{m}_{t}\right\|_{L^{2}}^{2}.
\end{aligned}
\end{equation}

Similarly, we can get $\mathbb{I}_{3}\leq C+\epsilon\left\|\nabla u^{m}_{t}\right\|_{L^{2}}^{2}$, and $\mathbb{I}_{4}\leq
 C+\epsilon\left\|\nabla u^{m}_{t}\right\|_{L^{2}}^{2}$, where $\epsilon$ is some small number.

For the term $\mathbb{I}_{5}$, we have
\begin{equation}
\begin{aligned}
\mathbb{I}_{5}&=\int_{\Omega}\rho v_{t}\partial_{x}u^{m}u^{m}_{t}\leq  \left\|\rho\right\|_{L_{\infty}}\left\|u^{m}_{t}\right\|_{L^{6}_{y}L^{2}_{x}}\left\|v_{t}\right\|_{L^{3}_{y}L^{\infty}_{x}}
\left\|\nabla u^{m}\right\|_{L^{2}}\\
&\leq C \left\|\rho\right\|_{L_{\infty}}\left\|\partial_{y}u^{m}_{t}\right\|_{L^{2}_{y}L^{2}_{x}}\left\|v_{t}\right\|_{W^{1,2}_{y}W^{1,2}_{x}}
\left\|\nabla u^{m}\right\|_{L^{2}}\\
&\leq C\left\|v_{t}\right\|_{H^{1}}^{2}+\epsilon\left\|\nabla u^{m}_{t}\right\|_{L^{2}}^{2},
\end{aligned}
\end{equation}
where we use Poincare in y-direction.

 \begin{equation}
 \begin{aligned}
 \mathbb{I}_{8}+\mathbb{I}_{6}
\leq C\left\|\sqrt{\rho}u^{m}_{t}\right\|_{L^{2}}^{2}+\epsilon\left\|\nabla u^{m}_{t}\right\|_{L_{2}}^{2}.
 \end{aligned}
 \end{equation}

For the term  $\mathbb{I}_{7}$, one has
\begin{equation}
\begin{aligned}
\mathbb{I}_{7}&\leq C \left\|\partial_{y}u^{m}\right\|_{L^{2}}
\left\|u^{m}_{t}\right\|_{L^{2}_{y}L^{\infty}_{x}}\left\|\int_{0}^{y}\partial_{x}v_{t}(s)ds\right\|_{L^{2}_{x}L^{\infty}_{y}}\\
&\leq C\left\|v_{t}\right\|_{H^{1}}^{2}+\epsilon\left\|\nabla u^{m}_{t}\right\|_{L^{2}}^{2}.
\end{aligned}
\end{equation}

A similiar argument to $\mathbb{I}_{9}$ yields that
\begin{equation}
\mathbb{I}_{9}\leq C (\left\|f_{t}\right\|_{L^{2}}^{2}+\left\|\sqrt{\rho}u^{m}_{t}\right\|_{L^{2}}^{2}).
\end{equation}

For the last term $\mathbb{I}_{10}$, we get that
\begin{equation}
\begin{aligned}
\mathbb{I}_{10}&\leq\left\|v\right\|_{L^{\infty}}\left\|\partial_{x}\rho\right\|_{L^{\infty}}\left\|f\right\|_{L^{\infty}_{y}L^{2}_{x}}
\left\|u^{m}_{t}\right\|_{L^{\infty}_{x}L^{2}_{y}}+\left\|f\right\|_{L^{2}}\left\|\partial_{x}\rho\right\|_{L^{\infty}}
\left\|\int_{0}^{y}\partial_{x}v(s)ds\right\|_{L^{\infty}_{y}L^{2}_{x}}
\left\|u^{m}_{t}\right\|_{L^{\infty}_{x}L^{2}_{y}}\\
&\quad+\lambda\|\rho_{xx}\|_{L^{\infty}}\left\|f\right\|_{L^{\infty}_{y}L^{2}_{x}}
\left\|u^{m}_{t}\right\|_{L^{\infty}_{x}L^{2}_{y}}
\leq C\left\|f\right\|_{H^{1}}^{2}+\epsilon\left\|\nabla u^{m}_{t}\right\|_{L^{2}}^{2}.
\end{aligned}
\end{equation}
where we use  Poincare inequality in y-direction.

Substituting all the estimates into $(\ref{319})$ with suitably small $\epsilon>0$, we obtain$$
\frac{\mathrm{d}}{\mathrm{d} t}\left\|\sqrt{\rho} u_{t}^{m}\right\|_{L^{2}}^{2}+\left\|u_{t}^{m}\right\|_{H^{1}}^{2} \leqslant C\left(1+\|f\|_{H^{1}}^{2}+\left\|f_{t}\right\|_{L^{2}}^{2}+\left\|v_{t}\right\|_{H^{1}}^{2}\right)+C\left\|\sqrt{\rho} u_{t}^{m}\right\|_{L^{2}}^{2}.$$


By Gronwall's inequality, we deduce that 
\begin{equation}\label{326}
\sup _{0 \leqslant t \leqslant T}\left\|\sqrt{\rho} u_{t}^{m}(t)\right\|_{L^{2}}^{2}+\int_{0}^{T}\left\|u_{t}^{m}(t)\right\|_{H^{1}}^{2} \mathrm{~d} t \leqslant C\left(1+\left\|\sqrt{\rho} u_{t}^{m}(0)\right\|_{L^{2}}^{2}\right).
\end{equation}

Recall that $V_{1}$ is given in ($\ref{1.5}$), essentially $V_{1}=u_{t}|_{t=0}$, and $\rho_{0}\geq
 \delta$.
Using the facts that   $u^{m}_{t}(0)=\Sum_{k=1}^{m}(u_{t}(0),\psi^{k})_{L^{2}}\psi^{k}\rightarrow u_{t}(0)$ in $L^{2}$ as $m\rightarrow\infty$, we deduce from $(\ref{318})$ and $(\ref{326})$ that\begin{equation}\label{330}
	\sup _{0 \leqslant t \leqslant T}\left(\left\|\nabla u^{m}(t)\right\|_{L^{2}}^{2}+\left\|\sqrt{\rho}u_{t}^{m}(t)\right\|_{L^{2}}^{2}\right)+\int_{0}^{T}\left\|u_{t}^{m}(t)\right\|_{H^{1}}^{2} \mathrm{~d} t \leqslant  \widetilde{C}.
\end{equation}

Using the uniform bound on $m$, there is a sequence $\left\{m_{k}\right\}$ such that $\left(\nabla u^{m_{k}}, u_{t}^{m_{k}}\right) \stackrel{*}{\rightarrow}$ $\left(\nabla u, u_{t}\right)$ in $L^{\infty}\left(0, T_{*} ; L^{2}\right)$ and $u_{t}^{m_{k}} \stackrel{w}{\rightarrow} u_{t}$ in $L^{2}\left(0, T_{*} ; H^{1}\right)$ for some limit function $u .$ 

By Aubin compactness theorem, we have also $$u_{m}\rightarrow u~\text{in}~ L^{2}(0,T;L^{2}(\Omega)) ~\text{strongly}.$$It is easy to see that $u$ is a weak solution to the original non-stationary hydrostatic Stokes equations (with the initial data $u_{0}$).
 Moreover, observe that for a.e. $t$ in $(0,T)$,$$\int_{\Omega} \mu(\rho)\nabla u\cdot\nabla \phi=\int_{\Omega}\rho(f-u_{t}-v\partial_{x}u+\int_{0}^{y}\partial_{x}v(s)ds\partial_{y}u)\cdot \phi$$for all $\phi\in X$, which follows from $(\ref{314})$ by passing to the limits. Hence it follows from the classical results in $\cite{abe}$ that for a.e. $t$ in $(0,T)$, there exists $P\in L^{2}$ such that $(u,P)$ is a weak solution of the stationary hydrostatic Stokes equations $(\ref{l33})$  with $F= \rho(f-u_{t}-v\partial_{x}u+\int_{0}^{y}\partial_{x}v(s)ds\partial_{y}u)$. From Lemma $\ref{Lemma-41}$ and estimates $(\ref{330})$, we easily deduce that $\nabla u, P\in L^{\infty}(0,T;H^{1})\cap L^{2}(0,T;H^{2})$. The time-continuity of $u$ in $H^{2}$ follows from a standard embedding result. This verifies $(u,\rho,P)$ is the solution to (\ref{37}) and \begin{equation}\label{3.38}
 	\begin{aligned}
 	&	\rho\in L^{\infty}(0,T;W^{2,\infty}\cap H^{3}),~\rho_{t}\in L^{\infty}(0,T;W^{1,2}\cap L^{\infty})\\&
 	u\in L^{\infty}(0,T;H^{2})\cap L^{2}(0,T;H^{3}),~u_{t}\in L^{\infty}(0,T;L^{2})\cap L^{2}(0,T;H^{1}).\end{aligned}
 \end{equation}

Next, we prove that $u\in L^{\infty}(0,T;H^{3})\cap L^{2}(0,T;H^{4})$, under the initial condition $u_{0} \in H^{3}$ and the compatibility condition (\ref{1.5}). Additionaly, supposed $u_{0}\in H^{3}$, from the equation $(\ref{3.1})_{2}$, one has $u_{t}|_{t=0}\in H^{1}(\Omega).$

We differentiate $(\ref{3.1})_{2}$ with respect to $t$ and multiply $u_{tt}$, one has
\begin{equation}\label{3.24}
\begin{aligned}		&\int_{\Omega}\rho|u_{tt}|^{2}+\frac{1}{2}\frac{d}{dt}\int_{\Omega}\mu(\rho)|\nabla u_{t}|^{2}\\&=-\int_{\Omega}\rho_{t}u_{t}u_{tt}-\int_{\Omega}\rho_{t}v\partial_{x}uu_{tt}-\int_{\Omega}\rho v_{t}\partial_{x}uu_{tt}-\int_{\Omega}\rho v\partial_{x}u_{t}u_{tt}\\
&\quad+\int_{\Omega}\rho_{t}\int_{0}^{y}\partial_{x}v(s)ds\partial_{y }uu_{tt}+\int_{\Omega}\rho \int_{0}^{y}\partial_{x}v_{t}(s)ds\partial_{y }u+\int_{\Omega}\rho\int_{0}^{y}\partial_{x}v(s)ds\partial_{y }u_{t}u_{tt}\\
&\quad+\int_{\Omega}div(\mu(\rho)_{t}\nabla u)u_{tt}+\int_{\Omega}\rho_{t}fu_{tt}+\int_{\Omega}\rho
f_{t}u_{tt}+\frac{1}{2}\int_{\Omega}\mu(\rho)_{t}|\nabla u_{t}|^{2}=\sum_{i=1}^{11}J_{i}.
\end{aligned}
\end{equation}
For $\mathbb{J}_{1}$, we have 
\begin{equation*}
	\mathbb{J}_{1}\leq C_{\delta}\|\rho_{t}\|_{L^{\infty}}\|\sqrt{\rho}u_{tt}\|_{L^{2}}\|u_{t}\|_{L^{2}}\leq C_{\delta}\|\nabla u_{t}\|_{L^{2}}^{2}+\epsilon\|\sqrt{\rho}u_{tt}\|_{L^{2}}^{2},
\end{equation*}
where we use  Poincare inequality in y-direction and $C_{\delta}$ is some positive constant depending on $\delta$.
To control 	$\mathbb{J}_{2}-	\mathbb{J}_{4}$, by Holder inequality, we have 
\begin{equation*}
	\begin{aligned}
\mathbb{J}_{2}+\mathbb{J}_{3}+\mathbb{J}_{4}&\leq
	C_{\delta}\|\rho_{t}\|_{L^{\infty}}\|v\|_{L^{\infty}}\|\partial_{x}u\|_{L^{2}}\|\sqrt{\rho}u_{tt}\|_{L^{2}}\\&\quad+\|\sqrt{\rho}\|_{L^{\infty}}\|v_{t}\|_{L^{\infty}}\|\partial_{x}u\|_{L^{2}}\|\sqrt{\rho}u_{tt}\|_{L^{2}}\\&\quad+\|\sqrt{\rho}\|_{L^{\infty}}\|v\|_{L^{\infty}}\|\sqrt{\rho}u_{tt}\|_{L^{2}}\|\partial_{x}u_{t}\|_{L^{2}}\\&\leq C_{\delta} +\epsilon\|\sqrt{\rho}u_{tt}\|_{L^{2}}^{2}+C\|v_{t}\|_{H^{2}}^{2}+C\|\nabla u_{t}\|_{L^{2}}^{2}.
	\end{aligned}
\end{equation*}
As the same manner of the derivation above, we can see 
\begin{equation*}
	\begin{aligned}
		\mathbb{J}_{5}+\mathbb{J}_{6}&\leq C_{\delta}\|\rho_{t}\|_{L^{\infty}}\|\int_{0}^{y}\partial_{x}vds\|_{L^{\infty}}\|\partial_{y }u\|_{L^{2}}\|\sqrt{\rho}u_{tt}\|_{L^{2}}\\&\quad+\|\sqrt{\rho}\|_{L^{\infty}}\|\int_{0}^{y}\partial_{x}v_{t}ds\|_{L^{\infty}}\|\partial_{y }u\|_{L^{2}}\|\sqrt{\rho}u_{tt}\|_{L^{2}}\\&\leq C_{\delta} +\epsilon\|\sqrt{\rho}u_{tt}\|_{L^{2}}^{2}+C\|v_{t}\|_{H^{2}}^{2},
	\end{aligned}
\end{equation*}
and 
\begin{equation*}
	\begin{aligned}
		\mathbb{J}_{7}+\mathbb{J}_{8}&\leq\|\sqrt{\rho}\|_{L^{\infty}}\|\int_{0}^{y}\partial_{x}v(s)ds\|_{L^{\infty}}\|\partial_{y }u_{t}\|_{L^{2}}\|\sqrt{\rho}u_{tt}\|_{L^{2}}\\&\quad+C_{\delta}\|\nabla \rho_{t}\|_{L^{2}}\|\nabla u\|_{L^{\infty}}\|\sqrt{\rho}u_{tt}\|_{L^{2}}\\&\quad+C_{\delta}\|\rho_{t}\|_{L^{\infty}}\|\Delta u\|_{L^{2}}\|\sqrt{\rho}u_{tt}\|_{L^{2}} \\&\leq \epsilon\|\sqrt{\rho}u_{tt}\|_{L^{2}}^{2}+C_{\delta}\|v_{t}\|_{H^{2}}^{2}.	\end{aligned}
\end{equation*}
Furthermore, we have
\begin{equation*}
	\begin{aligned}
		\mathbb{J}_{9}+\mathbb{J}_{10}+\mathbb{J}_{11}&\leq \|\rho_{t}\|_{\infty}\|f\|_{L^{2}}\|\sqrt{\rho}u_{tt}\|_{L^{2}}+\|\sqrt{\rho}\|_{\infty}\|f_{t}\|_{L^{2}}\|\sqrt{\rho}u_{tt}\|_{L^{2}}+\|\rho_{t}\|_{L^{\infty}}\|\nabla u_{t}\|^{2}\\&\leq C_{\delta} +\epsilon\|\sqrt{\rho}u_{tt}\|_{L^{2}}^{2}+C\|\nabla u_{t}\|_{L^{2}}^{2}.	\end{aligned}
\end{equation*}
Therefore, collecting the estimates above, after choosing $\epsilon$ small enough, we get \begin{equation*}
	\begin{aligned}
		\|\sqrt{\rho}u_{tt}\|^{2}+\frac{d}{dt}\|\nabla u_{t}\|^{2}_{L^{2}}\leq C_{\delta}\|\nabla
		u_{t}\|_{L^{2}}^{2}+C_{\delta}+C_{\delta}\|v_{t}\|_{H^{2}}^{2}.
	\end{aligned}
\end{equation*}

Then Gronwall inequality yields\begin{equation}\label{3.25}
	\int_{0}^{		t}\|\sqrt{\rho}u_{tt}(s)\|_{L^{2}}^{2}ds+\|\nabla u_{t}(t)\|_{L^{2}}^{2}\leq C_{\delta}.
\end{equation}

From Lemma
$\ref{Lemma-41}$ and $(\ref{37})_{2}$, we have \begin{equation}
	\|u\|_{H^{3}}\leq C \|F\|_{H^{1}}(1+\|\nabla\rho\|_{L^{ \infty}}+\|\nabla\rho\|_{L^{ \infty}}^{2})(1+\|\nabla^{2}\rho\|_{L^{2}}),
\end{equation}
where $$F=\rho(f-u_{t}-v\partial_{x}u+\int_{0}^{y}\partial_{x}v(s)ds\partial_{y }u).$$

	Therefore, using the above estimates (\ref{3.38}) and  (\ref{3.25}), we have  \begin{equation*}
F\in L^{\infty}(0,T;L^{2}),\nabla F\in L^{\infty}(0,T;L^{2})~\text{and}~u \in L^{\infty}(0,T;H^{3}).
	\end{equation*}

We differentiate  $(\ref{37})_{2}$ with respect to $t$, consider the hydrostatic Stokes problem:
$$-\mathrm{div}(\mu(\rho)\nabla u_{t})+\partial_{x}P_{t}=F_{t},$$ 
where 
\begin{eqnarray*}
  &&F_{t}=-\rho_{t}u_{t}-\rho u_{tt}-\rho_{t}v\partial_{x}u-\rho v \partial_{x}u_{t}-\rho v_{t}\partial_{x}u+\rho_{t}\int_{0}^{y}\partial_{x }v(s)ds\partial_{y}u+\rho\int_{0}^{y}\partial_{x }v_{t}(s)ds\partial_{y}u\\
  &&\quad\quad+\rho\int_{0}^{y}\partial_{x }v(s)ds\partial_{y}u_{t}+\mu(\rho)_{t}\nabla u+\mu(\rho)_{t}\Delta u+\rho_{t} f+\rho f_{t}.
\end{eqnarray*}
From Lemma
$\ref{Lemma-41}$, we have \begin{equation}\label{{3.44}}
	\|u_{t}\|_{H^{2}}\leq C\|F_{t}\|_{L^{2}}(1+\|\nabla\rho\|_{L^{ \infty}}).
\end{equation}
Combining with estimates (\ref{3.38}) and  (\ref{3.25}), we have $F_{t}\in L^{2}(0,T;L^{2})$ and $u_{t}\in L^{2}(0,T;H^{2})
$.

Following with similar arguments as above, we easily deduce that $u\in L^{2}(0,T;H^{4})$. The time-continuity of u in $H^{3}$ follows from a standard embedding result.
This completes the proof of the lemma.
$\hfill\Box$

\begin{Remark}\label{remark-1}
If $(\rho, u, P)$ is a strong solution to the linearized problem $(\ref{37})$-$(\ref{3122})$ satisfying $(\ref{311})$, then it satisfies the following identity analogous to $(\ref{319})$: for a.e. t in $(0,T)$,
\begin{equation}\label{444}
\begin{aligned}
&\frac{1}{2}\frac{\mathrm{d}}{\mathrm{d} t}\int_{\Omega}\rho\left|u_{t}\right|^{2}dx+\int_{\Omega}\mu(\rho)\left|\nabla u_{t}\right|^{2}\\
&=-\int_{\Omega}\rho_{t}(u_{t})^{2}-\int_{\Omega}\rho_{t}v\cdot\partial_{x}u u_{t}+\int_{\Omega}\rho_{t}\int_{0}^{y}\partial_{x}v(s)ds\cdot\partial_{y} uu_{t}-\int_{\Omega}\mu_{t}(\rho)\nabla u\nabla u_{t}\\
&\quad-\int_{\Omega}\rho v_{t}\partial_{x} u u_{t}+\int_{\Omega}\rho\int_{0}^{y}\partial_{x}v_{t}(s)ds\partial_{y}u u_{t}+\int_{\Omega}\rho f_{t}u_{t}+\int_{\Omega}\rho_{t}fu_{t}.
\end{aligned}
\end{equation}
\end{Remark}


\subsection{Proof of Proposition \ref{3.1}}

\ \ \
In this subsection, we mainly  prove the existence and uniqueness of strong solutions to the original system $(\ref{3.1})$ for the case of a positive initial density.

To see this, assume that $\mu\in C^{3}[0, \infty)$, $\rho_{0}\in W^{3,2}$, $\nabla \rho_{0}\in W^{1,\infty}$ and $\rho_{0}\geq\delta$ in $\bar{\Omega}$ for some constant $\delta>0$. We construct approximate solutions, inductively, as follows:
~\\
$(1)$ Define $u^{0}=0$, and
~\\
$(2)$ Assume that $u^{k-1}$ is defined for $k\geq1$. Then define $(\rho^{k},u^{k},P^{k})$ to be the unique global strong solution to the linearized problem $(\ref{37})$-$(\ref{3122})$ with $v$ replaced by $u^{k-1}$.
~

Next, we will give the uniform bounds of the approximate solutions and prove that these solutions converge to a local strong solution to the system $(\ref{3.1})$. Moreover,  we show that the local strong solution is unique.

\subsubsection{Uniform bounds}
To prove uniform bounds for the approximate solutions, we introduce a functional $\Phi_{K}(t)$ defined by
\begin{equation}
	\Phi_{K}(t)=\max _{1 \leqslant k \leqslant K}\left(1+\left\|\nabla \rho^{k}(t)\right\|_{W^{1,\infty}}+\|\nabla u^{k}_{t}(t)\|_{L^{2}}+\left\|\nabla
	^{2}\rho^{k}(t)\right\|_{W^{1,2}}+\left\|\nabla u^{k}(t)\right\|_{L^{2}}\right)
\end{equation}
where  $K\ge1$ an integer.  Second, we will estimate each term of $\Phi_{K}$ in terms of some integrals of $\Phi_{K}$, apply arguments of Gronwall-type inequality and thus prove that $\Phi_{K}$ is locally bounded. Finally, it is easy to see that the boundness of $\Phi_{K}$ is sufficient to give uniform bounds for the approximate solutions.

To simplify  the presentation, we use the following notations
$$\mu^{k}=\mu(\rho^{k}),~ k=1,2,\cdots.$$
Using the facts that $C\geq \rho_{0} \geq\delta$, according to the parabolic maximum principle, we deduce that
$$\delta \leqslant \rho^{k} \leqslant C~  \text { on }~  {[0, T] \times \bar{\Omega}}.$$
Furthermore, we have  $$ \left|\nabla \mu^{k}\right| \leqslant C\left|\nabla \rho^{k}\right|, $$
$$ \left|\nabla^{2} \mu^{k}\right| \leqslant C\left(\left|\nabla \rho^{k}\right|^{2}+\left|\nabla^{2} \rho^{k}\right|\right),\\~\text { and }~|\nabla^{3}\mu^{k}|\leq \widetilde{C}(|\nabla\rho^{k}|^{3}+|\nabla \rho^{k}||\nabla^{2}\rho^{k}|+|\nabla^{3}\rho^{k}|)~  \text { on }~  {[0, T] \times \bar{\Omega}}.
$$
\begin{Lemma}\label{Pro-1}
Assume that $(\rho^{k},u^{k},P^{k})$ is the unique global strong solution to the linearized problem $(\ref{37})$-$(\ref{3122})$ with $v$ replaced by $u^{k-1}$. Then, for any $t \in [0,T]$, it holds that
\begin{equation}\label{42}
\int_{0}^{t}\left\|\sqrt{\rho^{k}} u_{t}^{k}(s)\right\|_{L^{2}}^{2} \mathrm{~d} s+\left\|\nabla u^{k}(t)\right\|_{L^{2}}^{2} \leqslant C+C\int_{0}^{t} \Phi_{K}(s)^{8} \mathrm{~d} s
\end{equation}
for all $k, 1 \leqslant k \leqslant K$.
\end{Lemma}
{\bf Proof.}
Multiplying the linearized momentum equation by $u_{t}^{k}$ and integrating over $\Omega$, we derive
\begin{equation}\label{44444}
\begin{aligned}
&\int_{\Omega} \rho^{k}\left|u_{t}^{k}\right|^{2}+\frac{1}{2}\frac{\mathrm{d}}{\mathrm{d} t} \int_{\Omega} \mu^{k}\left|\nabla u^{k}\right|^{2} \\
& =  \int_{\Omega}\rho^{k}f u_{t}^{k}-\int_{\Omega} \rho^{k} u^{k-1} \partial_{x}u^{k} u_{t}^{k}
+\int_{\Omega} \rho^{k} \int_{0}^{y} \partial_{x}u^{k-1}(s)d s \partial_{y}u^{k}u_{t}^{k}+\frac{1}{2} \int_{\Omega} \mu^{\prime}\left(\rho^{k}\right) \rho^{k}_{t}\left|\nabla u^{k}\right|^{2}\\&=\int_{\Omega}\rho^{k}f u_{t}^{k}-\int_{\Omega} \rho^{k} u^{k-1} \partial_{x}u^{k} u_{t}^{k}
+\int_{\Omega} \rho^{k} \int_{0}^{y} \partial_{x}u^{k-1}(s)d s \partial_{y}u^{k}u_{t}^{k}\\&\quad+\int_{\Omega} u^{k-1}\partial_{x}\rho^{k} \left|\nabla u^{k}\right|^{2} -\int_{\Omega} \int_{0}^{y}\partial_{x} u^{k-1}(s)ds \partial_{y}\rho^{k}\left|\nabla u^{k}\right|^{2}-\lambda \rho^{k}_{xx}|\nabla u^{k}|^{2}=\sum_{i=1}^{6}I_{i}.
\end{aligned}
\end{equation}
where we use 
the linear transport equation  $(\ref{37})_1$ .

To estimate the right hand side of $(\ref{44444})$, the higher-order a priori estimates of the approximate solutions are required. Note that $\left(u^{k}, P^{k}\right) \in H_{0,per}^{1} \times L^{2}$ is a solution of the hydrostatic Stokes equations
$$-\partial_{x}\left(\mu(\rho^{k})\partial_{x}u^{k}\right)-\partial_{y}\left(\mu(\rho^{k})\partial_{y}u^{k}\right)+\partial_{x}P^{k}=F^{k} \quad \text{and}\quad \partial_{y} P^{k}=0 \quad \text  { in}
\quad \Omega,$$
where $F^{k}=\rho^{k}f-\rho^{k}u^{k}_{t}+\rho^{k}u^{k-1}\partial_{x}u^{k}+\rho^{k}\int_{0}^{y}\partial_{x}u^{k-1}(s)ds\partial_{y}u^{k} .$

Then, it follows from Lemma $\ref{Lemma-41}$ that
\begin{equation}\label{1}
\left\|u^{k}\right\|_{H^{2}}+\left\|P^{k}\right\|_{H^{1}}\leq C\left\|F^{k}\right\|_{L^{2}}\Phi_{K}.
\end{equation}

Using H$\ddot{o}$lder inequalities, Gagliardo-Nirenberg inequalities, and anisotropy Sobolev inequalities, one has
\begin{equation}\label{2}
\begin{aligned}
\left\|F^{k}\right\|_{L^{2}}&\leq \left\|\rho^{k}f-\rho^{k}u^{k}_{t}+\rho^{k}u^{k-1}\partial_{x}u^{k}+\rho^{k}\int_{0}^{y}\partial_{x}u^{k-1}(s)ds\partial_{y}u^{k}\right\|_{L^{2}}\\
&\leq  C\left\|\rho^{k}f\right\|_{L^{2}}+\left\|\sqrt{\rho^{k}}u_{t}^{k}\right\|_{L^{2}}+\left\|u^{k-1}\right\|_{L^{6}}\left\|\partial_{x} u^{k}\right\|_{L^{3}}\\
&\quad+\left\|\int_{0}^{y}\partial_{x}u^{k-1}(s)ds\right\|_{L_{y}^{\infty}L_{x}^{2}}\left\|\partial_{y}u^{k}\right\|_{L_{y}^{2}L_{x}^{\infty}}\\
&\leq C(1+\left\|\sqrt{\rho^{k}}u_{t}^{k}\right\|_{L^{2}}+\left\|u^{k}\right\|_{H^{2}}^{\frac{1}{3}}\Phi_{K}^{\frac{5}{3}}+\left\|u^{k}\right\|_{H^{2}}^{\frac{1}{2}}\Phi_{K}^{\frac{3}{2}}+\Phi_{K}^{2}),
\end{aligned}
\end{equation}
where we used the following facts:
$$
\left\|u^{k-1}\right\|_{L^{6}}\leq C_{1}\left\|\nabla u^{k-1}\right\|_{L^{2}}^{\frac{2}{3}}\left\| u^{k-1}\right\|_{L^{2}}^{\frac{1}
{3}}+C_{2}\left\| u^{k-1}\right\|_{L^{2}}\leq C\left\|\nabla u^{k-1}\right\|_{L^{2}},
$$
$$
\left\|\partial _{x}u^{k}\right\|_{L^{3}}\leq C_{1}\left\|\nabla\partial_{x} u^{k}\right\|_{L^{2}}^{\frac{1}{3}}\left\| \partial_{x} u^{k}\right\|_{L^{2}}^{\frac{2}
{3}}+C_{2}\left\|\partial_{x} u^{k}\right\|_{L^{2}}\\
\leq C_{1}\left\|\nabla^{2} u^{k}\right\|_{L^{2}}^{\frac{1}{3}}\left\| \nabla u^{k}\right\|_{L^{2}}^{\frac{2}
{3}}+C_{2}\left\|\nabla u^{k}\right\|_{L^{2}},
$$
and
$$
\left\|\partial_{y}u^{k}\right\|_{L_{y}^{2}L_{x}^{\infty}} \leq C_{1}\left\|\partial_{y}\partial_{x} u^{k}\right\|_{L^{2}}^{\frac{1}{2}}\left\| \partial_{y} u^{k}\right\|_{L^{2}}^{\frac{1}
{2}}+C_{2}\left\|\partial_{y} u^{k}\right\|_{L^{2}}\\
\leq C_{1}\left\|\nabla^{2} u^{k}\right\|_{L^{2}}^{\frac{1}{2}}\left\| \nabla u^{k}\right\|_{L^{2}}^{\frac{1}
{2}}+C_{2}\left\|\nabla u^{k}\right\|_{L^{2}}.
$$
Substituting (\ref{2}) into (\ref{1}) and using Young's inequality, we obtain
\begin{equation}\label{22}
\left\|u^{k}\right\|_{H^{2}}+\left\|P^{k}\right\|_{H^{1}}\leq C(\phi_{K}^{5}+\left\|\sqrt{\rho^{k}}u_{t}^{k}\right\|_{L^{2}}\phi_{K}).
\end{equation}

Applying ($\ref{22}$), we obtain
\begin{equation}\label{48888}
\begin{aligned}
I_{2}&\leq \epsilon\|\sqrt{\rho^{k}}u^{k}_{t}\|_{L^{2}}^{2}+C\int_{\Omega}\rho^{k}|u^{k-1}|^{2}|\partial_{x} u^{k}|^{2}\\&\leq\epsilon\|\sqrt{\rho^{k}}u^{k}_{t}\|_{L^{2}}^{2}+C\|\nabla u^{k-1}\|_{L^{2}}^{2}(\|\nabla u^{k}\|_{L^{2}}\|u^{k}\|_{H^{2}}+\|\nabla u^{k}\|_{L^{2}}^{2})\\&\leq
\epsilon\|\sqrt{\rho^{k}}u^{k}_{t}\|_{L^{2}}^{2}+C\phi_{K}^{8},
\end{aligned}
\end{equation}
 where we use\begin{equation*}
 	\begin{aligned}
 		\|\partial_{x}u^{k}\|_{L^{4}}\leq C_{1} \|\partial_{x}u^{k}\|_{L^{2}}^{\frac{1}{2}}\|\nabla \partial_{x}u^{k}\|_{L^{2}}^{\frac{1}{2}}+C_{2}\|\partial_{x}u^{k}\|_{L^{2}}.
 	\end{aligned}
 \end{equation*}
 Using H$\ddot{o}$lder inequalities, Gagliardo-Nirenberg inequalities, and anisotropy Sobolev inequalities, one has
\begin{equation}\label{49999}
\begin{aligned}
I_{3}&\leq C\|\sqrt{\rho^{k}}u^{k}_{t}\|_{L^{2}}\|\int_{0}^{y}\partial_{x}u^{k-1}(s)ds\|_{L^{\infty}L^{2}}\|\partial_{y}u^{k}\|_{L^{2}L^{\infty}}\\&\leq C\phi_{K}^{4}+\epsilon\|\sqrt{\rho^{k}}u^{k}_{t}\|_{L^{2}}^{2},
\end{aligned}
\end{equation}
and
\begin{equation}\label{40000}
\begin{aligned}
I_{4}=\int_{\Omega} u^{k-1}\partial_{x}\rho^{k} \left|\nabla u^{k}\right|^{2}
&\leq \left\|u^{k-1}\right\|_{L_{y}^{\infty}L_{x}^{2}}\left\|\partial_{x}\rho^{k}\right\|_{L^{\infty}}\left\|\nabla u^{k}\right\|_{L^{2}}\left\|\nabla u^{k}\right\|_{L_{y}^{2}L_{x}^{\infty}}\\
&\leq \left\|\nabla u^{k-1}\right\|_{L^{2}}\left\|\partial_{x} \rho^{k}\right\|_{L^{\infty}}\left\|\nabla u^{k}\right\|_{L^{2}}(C_{1}\left\|\nabla^{2} u^{k}\right\|_{L^{2}}^{\frac{1}{2}}\left\| \nabla u^{k}\right\|_{L^{2}}^{\frac{1}
	{2}}\\
&\quad+C_{2}\left\|\nabla u^{k}\right\|_{L^{2}})\leq C \Phi_{K}^{6}+\epsilon\left\|\sqrt{\rho^{k}}u_{t}^{k}\right\|_{L^{2}}^{2}.
\end{aligned}
\end{equation}
Similary, it yields
\begin{equation}
\begin{aligned}
I_{5}&=-\int_{\Omega} \int_{0}^{y}\partial_{x} u^{k-1}(s)ds\partial_{y}\rho^{k} \left|\nabla u^{k}\right|^{2}\\&\leq\left\|\int_{0}^{y}\partial_{x}u^{k-1}(s)ds\right\|_{L_{y}^{\infty}L_{x}^{2}}\left\|\partial_{x}\rho^{k}\right\|_{L^{\infty}}\left\|\nabla u^{k}\right\|_{L^{2}}\left\|\nabla u^{k}\right\|_{L_{y}^{2}L_{x}^{\infty}}\\
&\leq C \Phi_{K}^{6}+\epsilon\left\|\sqrt{\rho^{k}}u_{t}^{k}\right\|_{L^{2}}^{2},
\end{aligned}
\end{equation}and
\begin{equation}\label{41111}
	I_{6}=\int_{\Omega}-\lambda\rho^{k}_{xx}|\nabla u^{k}|_{L^{2}}\leq \lambda \|\rho^{k}_{xx}\|_{L^{\infty}}\Phi_{K}^{2}\leq C\Phi_{K}^{3}.
\end{equation}
Finally, substituting  ($\ref{48888}$)-($\ref{41111}$) into ($\ref{44444}$), with small $\epsilon$, we obtain$$\int_{\Omega} \rho^{k}\left|u_{t}^{k}\right|^{2} \mathrm{d} x+\frac{\mathrm{d}}{\mathrm{d} t} \int_{\Omega} \mu^{k}\left|\nabla u^{k}\right|^{2} \mathrm{d} x\leq C\Phi_{K}^{8}.$$

Integrating the above inequality over $(0,t)$, we complete the proof of the lemma.
$\hfill\Box$
\begin{Lemma}\label{Pro-2}
Assume that $(\rho^{k},u^{k},P^{k})$ is the unique global strong solution to the linearized problem $(\ref{37})$-$(\ref{3122})$ with $v$ replaced by $u^{k-1}$. Then, for any $t \in [0,T]$, it holds that
\begin{equation}
\left\|\sqrt{\rho^{k}} u_{t}^{k}(t)\right\|_{L^{2}}^{2}+\int_{0}^{t}\left\|\nabla u_{t}^{k}(s)\right\|_{L^{2}}^{2} \mathrm{~d} s \leqslant C(1+\|V_{1}\|_{L^{2}})\exp \left[C\int_{0}^{t} \Phi_{K}(s)^{16} \mathrm{~d} s\right]
\end{equation}
for all $k, 1 \leqslant k \leqslant K$.
\end{Lemma}
{\bf Proof.} Recall from (\ref{444})  that
\begin{equation}
\begin{aligned}
&\frac{1}{2}\frac{\mathrm{d}}{\mathrm{d} t}\int_{\Omega}\rho^{k}\left|u_{t}^{k}\right|^{2}+\int_{\Omega}\mu(\rho^{k})\left|\nabla u_{t}^{k}\right|^{2}+\mu_{t}(\rho^{k})\nabla u^{k} \nabla u_{t}^{k}\\&=
\int_{\Omega}\left[\rho^{k}(f_{t}-u^{k-1}\partial_{x}u_{t}^{k}+\int_{0}^{y}\partial_{x}u^{k-1}(s)ds\partial_{y}u_{t}^{k}-u_{t}^{k-1}\partial_{x}u^{k}
+\int_{0}^{y}\partial_{x}u_{t}^{k-1}(s)ds\partial_{y}u^{k})\right.\\
&\left.\quad\quad+\rho_{t}^{k}(f-\frac{1}{2}u_{t}^{k}-u^{k-1}\partial_{x}u^{k}+
\int_{0}^{y}\partial_{x}u^{k-1}(s)ds\partial_{y}u^{k})\right]u_{t}^{k}.
\end{aligned}
\end{equation}

Hence using the linearized continuity equation  $(\ref{37})_1$, we have
\begin{equation}\label{33}
\begin{aligned}
&\frac{1}{2}\frac{\mathrm{d}}{\mathrm{d} t}\int_{\Omega}\rho^{k}\left|u_{t}^{k}\right|^{2}+\int_{\Omega}\mu(\rho^{k})\left|\nabla u_{t}^{k}\right|^{2}\\&\leq
\int_{\Omega}\rho^{k}u_{t}^{k}f_{t}+\int_{\Omega}\rho^{k}u_{t}^{k}u^{k-1}\partial_{x}u_{t}^{k}+\int_{\Omega}[\rho^{k}u_{t}^{k}\int_{0}^{y}\partial_{x}u^{k-1}(s)ds\partial_{y}u_{t}^{k}]
+\int_{\Omega}\rho^{k}u_{t}^{k}u_{t}^{k-1}\partial_{x}u^{k}
\\&\quad+\int_{\Omega}[\rho^{k}u_{t}^{k}\int_{0}^{y}\partial_{x}u_{t}^{k-1}(s)ds\partial_{y}u^{k}]
+\int_{\Omega}(u^{k-1}\partial_{x}\rho^{k}-\int_{0}^{y}\partial_{x}u^{k-1}(s)ds\partial_{y}\rho^{k})fu_{t}^{k} \\ &\quad+\int_{\Omega}(u^{k-1}\partial_{x}\rho^{k}-\int_{0}^{y}\partial_{x}u^{k-1}(s)ds\partial_{y}\rho^{k})\frac{1}{2}u_{t}^{k}u_{t}^{k}\\
&\quad+\int_{\Omega}(u^{k-1}\partial_{x}\rho^{k}-\int_{0}^{y}\partial_{x}u^{k-1}(s)ds\partial_{y}\rho^{k})u^{k-1}\partial_{x}u^{k}u_{t}^{k}\\
&\quad+\int_{\Omega}(u^{k-1}\partial_{x}\rho^{k}-\int_{0}^{y}\partial_{x}u^{k-1}(s)ds\partial_{y}\rho^{k})\int_{0}^{y}\partial_{x}u^{k-1}(s)ds\partial_{y}u^{k}u_{t}^{k}\\
&\quad-\int_{\Omega}(u^{k-1}\partial_{x}\rho^{k}-\int_{0}^{y}\partial_{x}u^{k-1}(s)ds\partial_{y}\rho^{k})\nabla u^{k} \nabla u_{t}^{k}\\
&=C \sum_{i=1}^{10} J_{i}.
\end{aligned}
\end{equation}

 By H$\ddot{o}$lder inequalities, Gagliardo-Nirenberg inequalities and anisotropy Sobolev inequalities, it follows that
\begin{equation}
\begin{aligned}
\end{aligned}
\end{equation}
\begin{equation}
\begin{aligned}
J_{1}&\leq \left\|\rho^{k}\right\|_{L_{\infty}}^{\frac{1}{2}}\left\|f_{t}\right\|_{L_{2}}\left\|\sqrt{\rho^{k}}u_{t}^{k}\right\|_{L_{2}}\leq C ( \left\|f_{t}\right\|_{L_{2}}^{2}+\epsilon\left\|\sqrt{\rho^{k}}u_{t}^{k}\right\|_{L_{2}}^{2}).\\
J_{2}
&\leq\left\|\rho^{k}u_{t}^{k}\right\|_{L_{4}}\left\|\partial_{x}u_{t}^{k}\right\|_{L_{2}}
\left\|u^{k-1}\right\|_{L_{4}}\\
&\leq C(C_{1}\left\|\nabla(\rho^{k}u_{t}^{k})\right\|_{L_{2}}^{\frac{1}{2}}\left\|\rho^{k}u_{t}^{k}\right\|_{L_{2}}^{\frac{1}{2}}+
C_{2}\left\|\rho^{k}u_{t}^{k}\right\|_{L_{2}})(C_{1}\left\|\nabla u^{k-1}\right\|_{L_{2}}^{\frac{1}{2}}\left\|u^{k-1}\right\|_{L_{2}}^{\frac{1}{2}}\\&\quad+
C_{2}\left\|u^{k-1}\right\|_{L_{2}})\left\|\partial_{x}u_{t}^{k}\right\|_{L_{2}}\\
&\leq C(C_{1}\left\|\nabla\rho^{k}u_{t}^{k}+\rho^{k}\nabla u_{t}^{k}\right\|_{L_{2}}^{\frac{1}{2}}\left\|\rho^{k}u_{t}^{k}\right\|_{L_{2}}^{\frac{1}{2}}+
C_{2}\left\|\rho^{k}u_{t}^{k}\right\|_{L_{2}})\left\|\nabla u^{k-1}\right\|_{L_{2}}\left\|\partial_{x}u_{t}^{k}\right\|_{L_{2}}\\
&\leq C(C_{1}\Phi_{K}^{\frac{1}{2}}\left\|\nabla u_{t}^{k}\right\|_{L_{2}}^{\frac{1}{2}}\left\|\sqrt{\rho^{k}}u_{t}^{k}\right\|_{L_{2}}^{\frac{1}{2}}+C_{2}\left\|\sqrt{\rho^{k}}u_{t}^{k}\right\|_{L_{2}})
\left\|\nabla u^{k-1}\right\|_{L_{2}}\left\|\partial_{x}u_{t}^{k}\right\|_{L_{2}}\\
&\leq C\Phi_{K}^{\frac{3}{2}}\left\|\sqrt{\rho^{k}}u_{t}^{k}\right\|_{L_{2}}^{\frac{1}{2}}\left\|\nabla u_{t}^{k}\right\|_{L_{2}}^{\frac{3}{2}}+C\Phi_{K}^{\frac{3}{2}}\left\|\sqrt{\rho^{k}}u_{t}^{k}\right\|_{L_{2}}\left\|\nabla u_{t}^{k}\right\|_{L_{2}}\\
&\leq C\Phi_{K}^{6}\left\|\sqrt{\rho^{k}}u_{t}^{k}\right\|_{L_{2}}^{2}+\epsilon\left\|\nabla u_{t}^{k}\right\|_{L_{2}}^{2},\\
\end{aligned}
\end{equation}
where we have used
$$
\left\|\nabla\rho^{k}u_{t}^{k}+\rho^{k}\nabla u_{t}^{k}\right\|_{L_{2}}^{\frac{1}{2}}\leq C(\left\|\nabla\rho^{k}u_{t}^{k}\right\|_{L_{2}}^{\frac{1}{2}}+\left\|\rho^{k}\nabla u_{t}^{k}\right\|_{L_{2}}^{\frac{1}{2}})
\leq C \Phi_{K}^{\frac{1}{2}}\left\|\nabla u^{k}_{t}\right\|_{L_{2}}^{\frac{1}{2}}.
$$ Here and below, $\epsilon$ is some small number.

For the term $J_{3}$, we deduce
\begin{equation}
\begin{aligned}
J_{3}&\leq\left\|\int_{0}^{y}\partial_{x} u^{k-1}(s)ds\right\|_{L_{y}^{\infty}L_{x}^{2}}\left\|\rho^{k}u_{t}^{k}\right\|_{L_{y}^{2}L_{x}^{\infty}}
\left\|\partial_{y}u_{t}^{k}\right\|_{L_{2}}\\
&\leq C\left\|\partial_{x} u^{k-1}\right\|_{L_{2}}(\Phi_{K}^{\frac{1}{2}}\left\|\nabla u^{k}_{t}\right\|_{L_{2}}^{\frac{1}{2}}\left\|\sqrt{\rho^{k}}u_{t}^{k}\right\|_{L_{2}}^{\frac{1}{2}}+
\left\|\sqrt{\rho^{k}}u_{t}^{k}\right\|_{L_{2}})\|\nabla u^{k}_{t}\|_{L^{2}}\\
&\leq C\left\|\sqrt{\rho^{k}}u_{t}^{k}\right\|_{L_{2}}^{2}\Phi_{K}^{6}+\epsilon\left\|\nabla u_{t}^{k}\right\|_{L_{2}}^{2},\\
\end{aligned}
\end{equation}
where we have used
$$
\left\|\rho^{k}u_{t}^{k}\right\|_{L_{y}^{2}L_{x}^{\infty}}\leq
C\left\|\partial_{x}\rho^{k}u_{t}^{k}+\rho^{k}\partial_{x} u_{t}^{k}\right\|_{L_{2}}^{\frac{1}{2}}\left\|\rho^{k}u_{t}^{k}\right\|_{L_{2}}^{\frac{1}{2}}+
C\left\|\rho^{k}u_{t}^{k}\right\|_{L_{2}}.
$$

Similarly, it yields
\begin{equation*}
\begin{aligned}
J_{4}+J_{5}&\leq\left\|\rho^{k}u_{t}^{k}\right\|_{L_{3}}\left\|u_{t}^{k-1}\right\|_{L_{6}}
\left\|\partial_{x}u^{k}\right\|_{L_{2}}+\left\|\int_{0}^{y}\partial_{x} u_{t}^{k-1}(s)ds\right\|_{L_{y}^{\infty}L_{x}^{2}}\left\|\rho^{k}u_{t}^{k}\right\|_{L_{y}^{2}L_{x}^{\infty}}
\left\|\partial_{y}u^{k}\right\|_{L_{2}}\\
&\leq C(\Phi_{K}^{\frac{1}{3}}\left\|\nabla u^{k}_{t}\right\|_{L_{2}}^{\frac{1}{3}}\left\|\rho^{k}u_{t}^{k}\right\|_{L_{2}}^{\frac{2}{3}}+
\left\|\rho^{k}u_{t}^{k}\right\|_{L_{2}})\left\|\nabla u_{t}^{k-1}\right\|_{L_{2}}\left\|\partial_{x}u^{k}\right\|_{L_{2}}\\&
\quad+C(\Phi_{K}^{\frac{1}{2}}\left\|\nabla u^{k}_{t}\right\|_{L_{2}}^{\frac{1}{2}}\left\|\rho^{k}u_{t}^{k}\right\|_{L_{2}}^{\frac{1}{2}}+
\left\|\rho^{k}u_{t}^{k}\right\|_{L_{2}})\left\|\nabla u_{t}^{k-1}\right\|_{L_{2}}\left\|\partial_{y}u^{k}\right\|_{L_{2}}\\
&\leq C\Phi_{K}^{\frac{8}{3}}\left\|\nabla u^{k}_{t}\right\|_{L_{2}}^{\frac{2}{3}}\left\|\sqrt{\rho^{k}}u_{t}^{k}\right\|_{L_{2}}^{\frac{4}{3}}+\frac{\theta}{2}\left\|\nabla u_{t}^{k-1}\right\|_{L_{2}}^{2}+C\left\|\sqrt{\rho^{k}}u_{t}^{k}\right\|_{L_{2}}^{2}\Phi_{K}^{2}\\&\quad+C\Phi_{K}^{3}\left\|\nabla u^{k}_{t}\right\|_{L_{2}}\left\|\sqrt{\rho^{k}}u_{t}^{k}\right\|_{L_{2}}+\frac{\theta}{2}\left\|\nabla u_{t}^{k-1}\right\|_{L_{2}}^{2}+C\left\|\sqrt{\rho^{k}}u_{t}^{k}\right\|_{L_{2}}^{2}\Phi_{K}^{2}\\&\leq C\Phi_{K}^{6}\left\|\sqrt{\rho^{k}}u_{t}^{k}\right\|_{L_{2}}^{2}+\theta\left\|\nabla u_{t}^{k-1}\right\|_{L_{2}}^{2}+\epsilon\left\|\nabla u_{t}^{k}\right\|_{L_{2}}^{2},
\end{aligned}
\end{equation*}
where $\epsilon, \theta$ is some small number.

For the term $J_{6}$, we deduce
\begin{equation}
\begin{aligned}
J_{6}&=\int_{\Omega}(u^{k-1}\partial_{x}\rho^{k}-\int_{0}^{y}\partial_{x}u^{k-1}(s)ds\partial_{y}\rho^{k}-\lambda\rho^{k}_{xx})fu_{t}^{k} \\
&\leq\left\|u^{k-1}\right\|_{L_{4}}\left\|u^{k}_{t}\right\|_{L_{4}}\left\|\partial_{x}f\right\|_{L_{2}}\left\|\rho^{k}\right\|_{L_{\infty}}+
\left\|\int_{0}^{y}\partial_{x} u^{k-1}(s)ds\right\|_{L_{x}^{2}L_{y}^{\infty}}\left\|u^{k}_{t}\right\|_{L_{y}^{2}L_{x}^{\infty}}\left\|\partial_{y}f\right\|_{L_{2}}\left\|\rho^{k}\right\|_{L_{\infty}}\\
&\quad+\left\|u^{k-1}\right\|_{L_{4}}\left\|\partial_{x}u^{k}_{t}\right\|_{L_{2}}\left\|f\right\|_{L_{2}}\left\|\rho^{k}\right\|_{L_{\infty}}+
\left\|\int_{0}^{y}\partial_{x} u^{k-1}(s)ds\right\|_{L_{x}^{2}L_{y}^{\infty}}\left\|\partial _{y}u^{k}_{t}\right\|_{L_{2}}\left\|f\right\|_{L_{y}^{2}L_{x}^{\infty}}\left\|\rho^{k}\right\|_{L_{\infty}}\\
&\quad+\lambda\|\rho^{k}_{xx}\|_{L^{2}}\|u^{k}_{t}\|_{L^{\infty}L^{2}}\|f\|_{L^{\infty}L^{2}}\leq C\Phi_{K}^{2}+\epsilon\left\|\nabla u^{k}_{t}\right\|_{L_{2}}^{2}.
\end{aligned}
\end{equation}
Similarly, it yields
\begin{equation}
\begin{aligned}
J_{7}&=\int_{\Omega}(u^{k-1}\partial_{x}\rho^{k}-\int_{0}^{y}\partial_{x}u^{k-1}(s)ds\partial_{y}\rho^{k}-\lambda\rho^{k}_{xx})\frac{1}{2}u_{t}^{k}u_{t}^{k}\leq C(\epsilon,\delta)\left\|\sqrt{\rho^{k}}u_{t}^{k}\right\|_{L_{2}}^{2}\Phi_{K}^{6}+\epsilon\left\|\nabla u_{t}^{k}\right\|_{L_{2}}^{2}.\\
\end{aligned}
\end{equation}

To estimate the term $J_{8}$, we use $(\ref{22})$ to get
\begin{equation}
\begin{aligned}
J_{8}
&=\int_{\Omega}(u^{k-1}\partial_{x}\rho^{k}-\int_{0}^{y}\partial_{x}u^{k-1}(s)ds\partial_{y}\rho^{k}-\lambda\rho^{k}_{xx})u^{k-1}\partial_{x}u^{k}u_{t}^{k}\\
&\leq \left\|u^{k-1}\right\|_{L^{2}_{x}L^{\infty}_{y}}\left\|\partial_{x}\rho^{k}\right\|_{L^{\infty}}\left\|u^{k-1}\right\|_{L^{\infty}_{y}L^{2}_{x}}\left\|\partial_{x}u^{k}
\right\|_{L^{2}_{y}L^{\infty}_{x}}\left\|u^{k}_{t}\right\|_{L^{2}_{y}L^{\infty}_{x}}\\
&\quad+
\left\|\int_{0}^{y}\partial_{x}u^{k-1}(s)ds\right\|_{L^{\infty}_{y}L^{2}_{x}}\left\|\partial_{y}\rho^{k}\right\|_{L_{\infty}}\left\|u^{k-1}\right\|_{L_{x}^{2}L_{y}^{\infty}}
\left\|\partial_{x}u^{k}\right\|_{L_{y}^{2}L^{\infty}_{x}}\left\|u_{t}^{k}\right\|_{L_{y}^{2}L_{x}^{\infty}}\\&\quad+\lambda\|\rho^{k}_{xx}\|_{L^{\infty}}\|u^{k-1}\|_{L_{y}^{\infty}L_{x}^{2}}\|\partial_{x}u^{k}\|_{L^{2}}\|u^{k}_{t}\|_{L_{x}^{\infty}L_{y}^{2}}\\
&\leq \Phi_{K}^{3}\left\|u^{k}\right\|_{H^{2}}\left\|\nabla u^{k}_{t}\right\|_{L^{2}}+\Phi_{K}^{3}\left\|\nabla u^{k}_{t}\right\|_{L^{2}}\\
&\leq C\Phi_{K}^{16}+C\Phi_{K}^{8}\left\|\sqrt{\rho^{k}}u_{t}^{k}\right\|_{L^{2}}^{2}+\epsilon\left\|\nabla u^{k}_{t}\right\|_{L^{2}}^{2}.\\
\end{aligned}
\end{equation}

The terms $J_{9}$ and $J_{10}$ are estimated in a similar way:
\begin{equation}
\begin{aligned}
J_{9}&
=\int_{\Omega}(u^{k-1}\partial_{x}\rho^{k}-\int_{0}^{y}\partial_{x}u^{k-1}(s)ds\partial_{y}\rho^{k}-\lambda\rho^{k}_{xx})\int_{0}^{y}\partial_{x}u^{k-1}(s)ds\partial_{y}u^{k}u_{t}^{k}\\
&\leq \left\|u^{k-1}\right\|_{L^{2}_{y}L^{\infty}_{x}}\left\|\partial_{x}\rho^{k}\right\|_{L^{\infty}}\left\|\int_{0}^{y}\partial_{x}u^{k-1}(s)ds \right\|_{L^{2}_{x}L^{\infty}_{y}}\left\|\partial_{y}u^{k}\right\|_{L^{\infty}_{y}L^{2}_{x}}\left\| u^{k}_{t}\right\|_{L^{2}_{y}L^{\infty}_{x}}\\
&\quad+\left\|\int_{0}^{y}\partial_{x}u^{k-1}(s)ds\right\|_{L^{2}_{x}L^{\infty}_{y}}\left\|\partial_{y}\rho^{k}\right\|_{L^{\infty}}
\left\|\int_{0}^{y}\partial_{x}u^{k-1}(s)ds \partial_{y}u^{k}\right\|_{L^{2}_{y}L^{\infty}_{x}}\left\|\partial_{y}u^{k}\right\|_{L^{2}_{x}L^{\infty}_{x}}\left\| u^{k}_{t}\right\|_{L^{2}_{y}L^{\infty}_{x}}\\&\quad+\|\rho^{k}_{xx}\|_{L^{\infty}}\|\int_{0}^{y}\partial_{x}u^{k-1}\|_{L^{\infty}L^{2}}\|\partial_{y}u^{k}\|_{L^{2}}\|u^{k}_{t}\|_{L^{\infty}L^{2}}\\
&\leq \Phi_{K}^{3}\left\|u^{k}\right\|_{H^{2}}\left\|\nabla u^{k}_{t}\right\|_{L^{2}}+\Phi_{K}^{3}\left\|\nabla u^{k}_{t}\right\|_{L^{2}}\\
&\leq C\phi_{K}^{16}+C\phi_{K}^{8}\left\|\sqrt{\rho^{k}}u^{k}_{t}\right\|_{L^{2}}^{2}+\epsilon\left\|\nabla u^{k}_{t}\right\|_{L^{2}}^{2},
\end{aligned}
\end{equation}
and
\begin{equation}
\begin{aligned}
J_{10}
&=-\int_{\Omega}(u^{k-1}\partial_{x}\rho^{k}-\int_{0}^{y}\partial_{x}u^{k-1}(s)ds\partial_{y}\rho^{k}-\lambda\rho^{k}_{xx})\nabla u^{k} \nabla u_{t}^{k}\\
&\leq\left\|u^{k-1}\right\|_{L^{\infty}_{y}L^{2}_{x}}\left\|\partial_{x} \rho^{k}\right\|_{L_{\infty}}\left\|\nabla u^{k}\right\|_{L^{\infty}_{x}L^{2}_{y}}\left\|\nabla u_{t}^{k}\right\|_{L_{2}}\\&\quad+\left\|\int_{0}^{y}\partial_{x}u^{k-1}(s)ds\right\|_{L^{\infty}_{y}L^{2}_{x}}\left\|\partial_{y}\rho^{k}\right\|_{L_{\infty}}
\left\|\nabla u^{k}\right\|_{L^{2}_{y}L^{\infty}_{x}}
\left\|\nabla u_{t}^{k}\right\|_{L_{2}}\\&\quad+\lambda\|\rho^{k}_{xx}\|_{L^{\infty}}\|\nabla u^{k}\|_{L^{2}}\|\nabla u^{k}_{t}\|_{L^{2}}\\
&\leq C \Phi_{K}^{14}+C\Phi_{K}^{6}\left\|\sqrt{\rho^{k}}u_{t}^{k}\right\|_{L_{2}}^{2}+\epsilon\left\|\nabla u_{t}^{k}\right\|_{L_{2}}^{2}.
\end{aligned}
\end{equation}

Substituting these estimates of $J_{1}$-$J_{10}$ into $(\ref{33})$ with small $\theta, \epsilon>0$, we obtain
\begin{equation}
\begin{aligned}
&\frac{\mathrm{d}}{\mathrm{d} t}\int_{\Omega}\rho^{k}\left|u_{t}^{k}\right|^{2}dx+\int_{\Omega}\left|\nabla u_{t}^{k}\right|^{2}dx\\&\leq
C(1+\left\|f_{t}\right\|_{L_{2}}^{2})+C\Phi_{K}^{16}(1+\left\|\sqrt{\rho^{k}}u_{t}^{k}\right\|_{L_{2}}^{2})+\theta\left\|\nabla u_{t}^{k-1}\right\|_{L_{2}}^{2}.
\end{aligned}
\end{equation}

Thus, the desired estimates can be obtained exactly as in [$\cite{cho2004}$, Lemma $8$]. We complete the proof of Lemma \ref{Pro-2}.$\hfill\Box$

Now, we give estimates of $\left\|\nabla \rho^{k}\right\|_{L^{\infty}}
$ and $\|\nabla^{2}\rho^{k}\|_{L^{2}}$, which are
\begin{Lemma}\label{Pro-3}
Assume that $(\rho^{k},u^{k},P^{k})$ is the unique global strong solution to the linearized problem $(\ref{37})$-$(\ref{3122})$ with $v$ replaced by $u^{k-1}$. Then, for any $t \in [0,T]$, it holds that
\begin{equation}\label{433l}
\left\|\nabla \rho^{k}(t)\right\|_{L^{\infty}} \leqslant C \exp \left[C\left(1+\|V_{1}\|_{L^{2}}\right) \exp \left(C \int_{0}^{t} \Phi_{K}(s)^{20} \mathrm{~d} s\right)\right],
\end{equation}
and
\begin{equation}\label{431}
\left\|\nabla^{2} \rho^{k}(t)\right\|_{L^{2}} \leqslant C\exp \left[C\left(1+\|V_{1}\|_{L^{2}}\right) \exp \left(C\int_{0}^{t} \Phi_{K}(s)^{20} \mathrm{~d} s\right)\right]
\end{equation}
for all $k, 1 \leqslant k \leqslant K$.
\end{Lemma}
Before we give the proof of Lemma $\ref{Pro-3}$, we first prove
\begin{Lemma}\label{Pro-4}
Assume that $(\rho^{k},u^{k},P^{k})$ is the unique global strong solution to the linearized problem $(\ref{37})$-$(\ref{3122})$ with $v$ replaced by $u^{k-1}$. Then, for any $t \in [0,T]$, it holds that
\begin{equation}\label{433}
\left\|u^{k}\right\|_{H^{3}}+\left\|P^{k}\right\|_{H^{2}}\leq C (\left\|f\right\|_{H^{1}}^{2}+\left\|\sqrt{\rho^{k}} u^{k}_{t}\right\| _{L^{2}}^{2}+\left\|\sqrt{\rho^{k-1}}u^{k-1}_{t}\right\| _{L^{2}}^{2}+\|\nabla u^{k}_{t}\|_{L^{2}}\Phi_{K}^{4}+\Phi_{K}^{20})
\end{equation}
for all $k, 1 \leqslant k \leqslant K$.
\end{Lemma}
{\bf Proof.}
  It follows from  Lemma $\ref{Lemma-41}$ that $$
\|u^{k}\|_{H^{3}}+\|P^{k}\|_{H^{2}} \leqslant C\|F^{k}\|_{H^{1}}(1+\|\nabla\rho^{k}\|_{L^{ \infty}}+\|\nabla\rho^{k}\|_{L^{ \infty}}^{2})(1+\|\nabla^{2}\rho^{k}\|_{L^{2}}),$$
where $F^{k}=\rho^{k}f-\rho^{k}u^{k}_{t}+\rho^{k}u^{k-1}\partial_{x}u^{k}+\rho^{k}\int_{0}^{y}\partial_{x}u^{k-1}(s)ds\partial_{y}u^{k}$.

Therefore, we have
\begin{equation}\label{43111}
\begin{aligned}\|u^{k}\|_{H^{3}}+\|P^{k}\|_{H^{2}} \leq C\left\|F^{k}\right\|_{H^{1}}\Phi_{K}^{3}&\leq C (\|F^{k}\|_{L^{2}}\Phi_{K}^{3}+\left\|\nabla F^{k}\right\|_{L^{2}}\Phi_{K}^{3}).\\
\end{aligned}
\end{equation}

The first term on the right hand side of (\ref{43111}) can be estimated in a similar way as in ($\ref{2}$).

To estimate the second term on the right hand side of (\ref{43111}), one has
\begin{equation}\label{434}
\begin{aligned}
\left\|\nabla F^{k}\right\|_{L^{2}}\Phi_{K}^{3}
&\leq( \left\|\nabla \rho^{k}\right\|_{L^{\infty
}}\left\|f\right\|_{L^{2}}+\left\|\rho^{k}\right\|_{L^{\infty}}\left\|\nabla f\right\|_{L_{2}}+\left\|\nabla\rho^{k}\right\|_{L^{\infty}}\left\|\partial_{t}u^{k}\right\|_{L^{2}}\\
&\quad+ \left\|\rho^{k}\right\|_{L^{\infty}}\left\|\nabla\partial_{t}u^{k}\right\|_{L^{2}}+\left\|\nabla \rho^{k}\right\|_{L^{\infty}}\left\|u^{k-1}\right\|_{L^{\infty}_{x}L^{2}_{y}}\left\|\partial_{x}u^{k}\right\|_{L^{\infty}_{y}L^{2}_{x}}\\
&\quad+ \left\|\nabla u^{k-1}\right\|_{L^{\infty}_{y}L^{2}_{x}} \left\|\partial_{x} u^{k}\right\|_{L^{2}_{y}L^{\infty}_{x}}\left\|\rho^{k}\right\|_{L^{\infty}}+\left\| u^{k-1}\right\|_{L^{\infty}_{y}L^{2}_{x}} \left\|\nabla \partial_{x} u^{k}\right\|_{L^{2}_{y}L^{\infty}_{x}}\left\|\rho^{k}\right\|_{L^{\infty}}\\
&\quad+\left\|\nabla \rho^{k}\right\|_{L^{\infty}}\left\|\int_{0}^{y}\partial_{x}u^{k-1}(s)ds\right\|_{L_{y}^{\infty}L_{x}^{2}}\left\|\partial_{y}u^{k}\right\|_{L_{x}^{\infty}L_{y}^{2}}+
\left\|\rho^{k}\int_{0}^{y}\partial_{xx}u^{k-1}(s)ds\partial_{y}u^{k}\right\|_{L^{2}}\\
&\quad+\left\|\rho^{k}\partial_{x}u^{k-1}\partial_{y}u^{k}\right\|_{L^{2}}
+\left\|\rho^{k}\right\|_{L^{\infty}}\left\|\int_{0}^{y}\partial_{x}u^{k-1}(s)ds\right\|_{L^{\infty}_{y}L^{2}_{x}}\left\|\nabla \partial_{y}u^{k}\right\|_{L^{\infty}_{x}L^{2}_{y}})\Phi_{K}^{3}\\
&=\sum_{i=1}^{11}L_{i}.
\end{aligned}
\end{equation}
The estimates for $L_{1}$-$L_{4}$ are given as follows.

$$L_{1}+L_{2}\leq \Phi_{K}^{4}\left\|f\right\|_{H^{1}},$$
and
$$L_{3}+L_{4}\leq \Phi_{K}^{4}\left\|\nabla u^{k}_{t}\right\|_{L_{2}}.$$

To estimate the term $L_{5}$, thanks to the Gagliardo-Nirenberg inequality, one has
$$L_{5}\leq \|\nabla \rho^{k}\|_{L^{\infty}}\|\nabla u^{k-1}\|_{L^{2}}\|u^{k}\|_{H^{2}}\Phi_{K}^{3}\leq C\|\sqrt{\rho} u^{k}_{t}\|_{L^{2}}\phi_{K}^{6}+C\phi_{K}^{10}.$$

A similar argument yields
\begin{equation}
\begin{aligned}
L_{6}&\leq(C_{1}\left\|\partial_{y}\nabla u^{k-1}\right\|_{L^{2}}^{\frac{1}{2}}\left\|\nabla u^{k-1}\right\|_{L^{2}}^{\frac{1}{2}}+C_{2}\left\|\nabla u^{k-1}\right\|_{L^{2}})(C_{1}\left\|\partial_{xx} u^{k}\right\|_{L^{2}}^{\frac{1}{2}}\left\|\partial_{x} u^{k}\right\|_{L^{2}}^{\frac{1}{2}}\\&\quad+C_{2}\left\|\partial_{x} u^{k}\right\|_{L^{2}})\Phi_{K}^{3}\\
&\leq C\left\|\sqrt{\rho} u^{k-1}_{t}\right\|_{L^{2}}\Phi_{K}^{8}+\Phi_{K}^{9}+\epsilon\|u^{k}\|_{H^{2}}.
\end{aligned}
\end{equation}

Using $(\ref{22})$ and the Gagliardo-Nirenberg inequality, we deduce
\begin{equation}
\begin{aligned}
L_{7}&\leq C(\left\|\nabla u^{k-1}\right\|_{L^{2}}\left\|u^{k}\right\|_{H^{3}}^{\frac{1}{2}}\left\|\nabla^{2}u^{k}\right\|_{L^{2}}^{\frac{1}{2}}+
\left\|\nabla u^{k-1}\right\|_{L^{2}}\left\|\nabla^{2}u^{k}\right\|_{L^{2}})\phi_{K}^{8}\\
&\leq \epsilon\left\|u^{k}\right\|_{H^{3}}+C\left\|\sqrt{\rho} u^{k}_{t}\right\|_{L^{2}}\Phi_{K}^{9}+\Phi_{K}^{13},
\end{aligned}
\end{equation}
where we have used
$$\left\|\nabla\partial_{x}u^{k}\right\|_{L^{2}_{y}L^{\infty}_{x}}\leq C_{1}\left\|\partial_{x}\nabla\partial_{x}u^{k}\right\|_{L^{2}}^{\frac{1}{2}}\left\|\nabla\partial_{x}u^{k}\right\|_{L^{2}}^{\frac{1}{2}}+
C_{2}\left\|\nabla\partial_{x}u^{k}\right\|_{L^{2}}.$$

The estimates for $L_{8}$ and $L_{9}$ are given as follows.
\begin{equation}
\begin{aligned}
L_{8}&\leq \Phi_{K}^{4} \left\|\partial_{x}u^{k-1}\right\|_{L^{2}}(C_{1}\left\|u^{k}\right\|_{H^{2}}^{\frac{1}{2}}\left\|\nabla u^{k}\right\|_{L^{2}}^{\frac{1}{2}}+C_{2}\left\|\nabla u^{k}\right\|_{L^{2}})\\
&\leq  C\Phi_{K}^{11}+\epsilon\left\|u^{k}\right\|_{H^{2}}.
\end{aligned}
\end{equation}

\begin{equation}
\begin{aligned}
L_{9}&\leq  \left\|\int_{0}^{y}\partial_{xx}u^{k-1}(s)ds\right\|_{L^{\infty}_{y}L^{2}_{x}}
\left\|\partial_{y}u^{k}\right\|_{L^{\infty}_{x}L^{2}_{y}}\phi_{K}^{3}\\
&\leq  \left\|u^{k-1}\right\|_{H^{2}}(C_{1}\left\|u^{k}\right\|_{H^{2}}^{\frac{1}{2}}\left\|\nabla u^{k}\right\|_{L^{2}}^{\frac{1}{2}}+C_{2}\left\|\nabla u^{k}\right\|_{L^{2}})\phi_{K}^{3}\\
&\leq C (\left\|\sqrt{\rho} u^{k}_{t}\right\|_{L^{2}}^{2}+\left\|\sqrt{\rho} u^{k-1}_{t}\right\|_{L^{2}}^{2}+\Phi_{K}^{20}).
\end{aligned}
\end{equation}

The estimate of $L_{10}$ is the same as that of $L_{6}$. For the last term $L_{11}$, we get that
\begin{equation}
\begin{aligned}
L_{11}&\leq  \Phi_{K}^{4}(C_{1} \left\|u^{k}\right\|_{H^{3}}^{\frac{1}{2}}\left\|u^{k}\right\|_{H^{2}}^{\frac{1}{2}}+C_{2}\left\|u^{k}\right\|_{H^{2}})\\
&\leq \epsilon\left\|u^{k}\right\|_{H^{3}}+C\left\|\sqrt{\rho} u^{k}_{t}\right\|_{L^{2}} \Phi_{K}^{9}+\phi_{K}^{13}.
\end{aligned}
\end{equation}

Substituting estimates $L_{1}$-$L_{11}$ into $(\ref{434})$ with small $\epsilon>0$, we obtain
$$\left\|u^{k}\right\|_{H^{3}}+\left\|P^{k}\right\|_{H^{2}}\leq C(\left\|f\right\|_{H^{1}}^{2}+\left\|\sqrt{\rho^{k}}u^{k}_{t}\right\| _{L^{2}}^{2}+\left\|\sqrt{\rho^{k-1}} u^{k-1}_{t}\right\| _{L^{2}}^{2}+\|\nabla u^{k}_{t}\|_{L^{2}}\Phi_{K}^{4}+\Phi_{K}^{20})\leq
C(1+\Phi_{K}^{20}).$$

The proof of the lemma is finished.$\hfill\Box$

 Now we prove Lemma $\ref{Pro-3}$.

\textbf{Proof of Lemma $\ref{Pro-3}$}
It follows from $(\ref{2.4})$ and $(\ref{2.5})$ that $$
\|\nabla \rho^{k}(t)\|_{L^{ \infty}} \leqslant\left\|\nabla\rho_{0}\right\|_{L^{ \infty}} \exp \left(C\int_{0}^{t}\|u^{k-1
}(s)\|_{H^{3}} \mathrm{~d} s\right),$$
and$$
\|\nabla^{2}\rho^{k}(t)\|_{L^{ 2}} \leqslant(\left\|\nabla^{2} \rho_{0}\right\|_{L^{2}}+\|\nabla \rho_{0}\|_{L^{\infty}}) \exp \left(C \int_{0}^{t}\|u^{k-1}(s)\|_{H^{3}} \mathrm{~d} s\right).$$

Applying Lemma $\ref{Pro-4}$ and Lemma $\ref{Pro-2}$ leads to
\begin{equation}
\begin{aligned}
\left\|\nabla\rho^{k}\right\|_{L^{\infty}}&\leq C \exp \big[C \int_{0}^{t} \Phi_{K}(s)^{20}ds+C\int_{0}^{t}\left\|\nabla u^{k-1}_{t}\right\|_{L^{2}}^{2}ds+C\int_{0}^{t}\left\|\nabla u^{k-2}_{t}\right\|_{L^{2}}^{2}ds\big]\\
&\leq C \exp [\exp(C\int_{0}^{t}\Phi_{K}(s)^{20}ds)+C (1+\|V_{1}\|_{L^{2}})\exp(C \int_{0}^{t}\Phi_{K}(s)^{16}ds)]\\
&\leq C \exp [C (1+\|V_{1}\|_{L^{2}})\exp \big( \int_{0}^{t}\Phi_{K}(s)^{20}ds)].
\end{aligned}
\end{equation}

Similarly, we can obtain $$\left\|\nabla^{2} \rho^{k}(t)\right\|_{L^{2}} \leqslant C \exp \left[C\left(1+\|V_{1}\|_{L^{2}}\right) \exp \left(C \int_{0}^{t} \Phi_{K}^{20}(s) \mathrm{~d} s\right)\right].$$ The proof of Lemma $\ref{Pro-3}$ is complete.$\hfill\Box$

Now, we give eatimates of $\left\|\nabla^{2} \rho^{k}\right\|_{L^{\infty}}
$ and $\|\nabla^{3}\rho^{k}\|_{L^{2}}$, which are
\begin{Lemma}\label{Pro-3.6}
	Assume that $(\rho^{k},u^{k},P^{k})$ is the unique global strong solution to the linearized problem $(\ref{37})$-$(\ref{3122})$ with $v$ replaced by $u^{k-1}$. Then, for any $t \in [0,T]$, it holds that
	\begin{equation}\label{3.65}
		\left\|\nabla^{2} \rho^{k}(t)\right\|_{L^{\infty}} \leqslant C\exp[\widetilde{C}_{\delta}(1+\|V_{1}\|_{L^{2}})^{2}\exp(\widetilde{C}\int_{0}^{t}\Phi_{K}^{42}(s)ds)],
	\end{equation}
	and
	\begin{equation}\label{3.66}
		\left\|\nabla^{3} \rho^{k}(t)\right\|_{L^{2}}\leq C\exp[\widetilde{C}_{\delta}(1+\|V_{1}\|_{L^{2}})^{2}\exp(\widetilde{C}\int_{0}^{t}\Phi_{K}^{42}(s)ds)].
	\end{equation}
	for all $k, 1 \leqslant k \leqslant K$.
\end{Lemma}
Before we prove Lemma \ref{Pro-3.6}, we still need the following lemma:
\begin{Lemma}	Assume that $(\rho^{k},u^{k},P^{k})$ is the unique global strong solution to the linearized problem $(\ref{37})$-$(\ref{3122})$ with $v$ replaced by $u^{k-1}$. Then, for any $t \in [0,T]$, it holds that
	\begin{equation}\label{3.67}
\int_{0}^{t}\|\sqrt{\rho^{k}}u^{k}_{tt}(s)\|^{2}_{L^{2}}ds+\|\nabla u^{k}_{t}(t)\|_{L^{2}}^{2}\leq C_{\delta}
(1+\|V_{1}\|_{L^{2}})^{2}\exp[C\int_{0}^{t}\Phi_{K}^{42}(s)ds],
\end{equation}
	for all $k, 1 \leqslant k \leqslant K$ and $C_{\delta}$ is a positive constant depending on $\delta$.
\end{Lemma}
{\bf Proof.} We differentiate $(\ref{3.1})_{2}$ with respect to $t$ and multiply $u^{k}_{tt}$, one has
\begin{equation}\label{3.67}
	\begin{aligned}
	&	\int_{\Omega}\rho^{k}|u^{k}_{tt}|^{2}+\frac{1}{2}\frac{d}{dt}\int_{\Omega}\mu(\rho^{k})|\nabla u^{k}_{t}|^{2}\\&=-\int_{\Omega}\rho^{k}_{t}u^{k}_{t}u^{k}_{tt}-\int_{\Omega}\rho^{k}_{t}u^{k-1}\partial_{x}u^{k}u^{k}_{tt}-\int_{\Omega}\rho^{k}u^{k-1}_{t}\partial_{x}u^{k}u^{k}_{tt}-\int_{\Omega}\rho^{k}u^{k-1}\partial_{x}u_{t}^{k}u^{k}_{tt}\\&\quad-\int_{\Omega}\rho^{k}_{t}\int_{0}^{y}\partial_{x}u^{k-1}\partial_{y}u^{k}u^{k}_{tt}-\int_{\Omega}\rho^{k}\int_{0}^{y}\partial_{x}u^{k-1}_{t}\partial_{y}u^{k}u^{k}_{tt}-\int_{\Omega}\rho^{k}\int_{0}^{y}\partial_{x}u^{k-1}\partial_{y}u_{t}^{k}u^{k}_{tt}\\&\quad+\int_{\Omega}\nabla \rho^{k}_{t}\nabla u^{k}u^{k}_{tt}+\int_{\Omega} \rho^{k}_{t}\nabla^{2} u^{k}u^{k}_{tt}+\frac{1}{2}\int_{\Omega}\mu(\rho^{k})_{t}|\nabla u^{k}_{t}|^{2}+\int_{\Omega}\rho^{k}_{t}fu^{k}_{tt}+\int_{\Omega}\rho^{k} f_{t}u^{k}_{tt}=\sum_{i=1}^{12}M_{i}.
	\end{aligned}
\end{equation}
Then one has the following estimates to the terms in the right-hand side of (\ref{3.67}).
\begin{equation*}
	\begin{aligned}
		M_{1}=-\int_{\Omega}\rho^{k}_{t}u^{k}_{t}u^{k}_{tt}\leq\epsilon
	\|\sqrt{\rho^{k}}	u^{k}_{tt}\|^{2}_{L^{2}}+C_{\delta}\|\rho^{k}_{t}\|^{2}_{L^{\infty}}\|\nabla u^{k}_{t}\|^{2}_{L^{2}}\leq \epsilon
	\|\sqrt{\rho^{k}}	u^{k}_{tt}\|^{2}_{L^{2}}+C_{\delta}\Phi_{K}^{4}\|\nabla u^{k}_{t}\|^{2}_{L^{2}}.
	\end{aligned}
\end{equation*}
\begin{equation*}
	\begin{aligned}
	M_{2}=\int_{\Omega}\rho^{k}_{t}u^{k-1}\partial_{x}u^{k}u^{k}_{tt}\leq
C_{\delta}	\|\sqrt{\rho^{k}}	u^{k}_{tt}\|_{L^{2}}\|\rho^{k}_{t}\|_{L^{\infty}}\|\nabla u^{k}\|_{L_{x}^{\infty}L_{y}^{2}}\|u^{k-1}\|_{L_{y}^{\infty}L_{x}^{2}}\leq
	\epsilon
	\|\sqrt{\rho^{k}}	u^{k}_{tt}\|^{2}_{L^{2}}+C_{\delta}\Phi_{K}^{16}.
	\end{aligned}
\end{equation*}
\begin{equation*}
	\begin{aligned}M_{3}&=-\int_{\Omega}\rho^{k}u^{k-1}_{t}\partial_{x}u^{k}u^{k}_{tt}\leq\|\sqrt{\rho^{k}}\|_{L^{\infty}}\|\sqrt{\rho^{k}}	u^{k}_{tt}\|_{L^{2}}\|u^{k-1}_{t}\|_{L_{x}^{\infty}L_{y}^{2}}\|\partial_{x}u^{k}\|_{L_{y}^{2}L_{x}^{\infty}}\\&\leq	\epsilon\|\sqrt{\rho^{k}}	u^{k}_{tt}\|^{2}_{L^{2}}+C_{\delta}\Phi_{K}^{12}.
	\end{aligned}
\end{equation*}
\begin{equation*}
	\begin{aligned}
		M_{4}&=-\int_{\Omega}\rho^{k}u^{k-1}\partial_{x}u_{t}^{k}u^{k}_{tt}\leq\|\sqrt{\rho^{k}}\|_{L^{\infty}}	\|\sqrt{\rho^{k}}u^{k}_{tt}\|_{L^{2}}\|u^{k-1}\|_{L^{\infty}}\|\partial_{x}u_{t}^{k}\|_{L^{2}}\\&\leq	\epsilon\|\sqrt{\rho^{k}}	u^{k}_{tt}\|^{2}_{L^{2}}+C\Phi_{K}^{10}\|\partial_{x}u_{t}^{k}\|_{L^{2}}^{2},
	\end{aligned}
\end{equation*}
 where we use\begin{equation*}
 \|u^{k}\|_{H^{2}}\leq C\Phi_{K}^{5}.
 \end{equation*}  Similarly as above, we have
\begin{equation*}
	M_{5}\leq \epsilon
	\|\sqrt{\rho^{k}}	u^{k}_{tt}\|^{2}_{L^{2}}+C_{\delta}\Phi_{K}^{16},
\end{equation*}
\begin{equation*}
M_{6}\leq \epsilon\|\sqrt{\rho^{k}}	u^{k}_{tt}\|^{2}_{L^{2}}+C_{\delta}\Phi_{K}^{12},
\end{equation*}
 and
\begin{equation*}
M_{7}\leq\epsilon\|\sqrt{\rho^{k}}	u^{k}_{tt}\|^{2}_{L^{2}}+C\Phi_{K}^{10}\|\partial_{x}u_{t}^{k}\|_{L^{2}}^{2}.
\end{equation*}
We get by $(\ref{37})_{1}$ that 
\begin{equation*}
	\begin{aligned}
	M_{8}&=\int_{\Omega}\nabla \rho^{k}_{t}\nabla u^{k}u^{k}_{tt}\\&=\int_{\Omega}	u^{k}_{tt}(\nabla u^{k-1}\partial_{x}\rho+u^{k-1}\nabla \partial_{x}\rho^{k}+\partial_{x }u^{k-1}\partial_{y }\rho^{k}+\int_{0}^{y}\partial_{x}u^{k-1}\nabla \partial_{y}\rho\\&\quad+\int_{0}^{y}\partial_{x x}u^{k-1}\partial_{y}\rho^{k}+\lambda\nabla\partial_{x x}\rho^{k})\nabla u^{k}\\&
	\leq C_{\delta}(\|\sqrt{\rho^{k}}	u^{k}_{tt}\|_{L^{2}}\|u^{k-1}\|_{H^{2}}\Phi_{K}\|u^{k}\|_{H^{2}}+\|\sqrt{\rho^{k}}	u^{k}_{tt}\|_{L^{2}}\Phi_{K}^{2}\|u^{k}\|_{H^{2}}+\|\sqrt{\rho^{k}}	u^{k}_{tt}\|_{L^{2}}\Phi_{K}\|u^{k}\|_{H^{3}})\\&\leq
	\epsilon\|\sqrt{\rho^{k}}	u^{k}_{tt}\|_{L^{2}}^{2}+C_{\delta}\Phi_{K}^{42}, 
	\end{aligned}
\end{equation*}
where we use \begin{equation*}
	\|u^{k}\|_{H^{3}}\leq C(1+\Phi_{K}^{20}).
\end{equation*}
For $M_{9}$ and $M_{10}$, we have 
\begin{equation*}
\begin{aligned}
		M_{9}+M_{10}&=\int_{\Omega} \rho^{k}_{t}\nabla^{2} u^{k}u^{k}_{tt}+\frac{1}{2}\int_{\Omega}\mu(\rho^{k})_{t}|\nabla u^{k}_{t}|^{2}\\&\leq	\epsilon\|\sqrt{\rho^{k}}	u^{k}_{tt}\|^{2}_{L^{2}}+C_{\delta}\Phi_{K}^{6}\|\nabla u_{t}^{k}\|_{L^{2}}^{2}+C_{\delta}\Phi_{K}^{14}.
	\end{aligned}
\end{equation*}
\begin{equation*}
	\begin{aligned}
	M_{11}+M_{12}&=\int_{\Omega}\rho^{k}_{t}fu^{k}_{tt}+\int_{\Omega}\rho^{k} f_{t}u^{k}_{tt}\\&\leq C_{\delta}\|\rho^{k}_{t}\|_{L^{\infty}}\|f\|_{L^{2}}\|\sqrt{\rho^{k}}u^{k}_{tt}\|_{L^{2}}+C_{\delta}\|\rho^{k}\|_{L^{\infty}}\|f_{t}\|_{L^{2}}\|\sqrt{\rho^{k}}u^{k}_{tt}\|_{L^{2}}\\&\leq\epsilon\|\sqrt{\rho^{k}}	u^{k}_{tt}\|^{2}_{L^{2}}+C_{\delta}\Phi_{K}^{12}.
\end{aligned}
\end{equation*}
Summing up the above inequalities with small enough $\epsilon$ yields
\begin{equation*}
	\begin{aligned}
	\int_{\Omega}\rho^{k}|u^{k}_{tt}|^{2}+\frac{1}{2}\frac{d}{dt}\int_{\Omega}\mu(\rho^{k})|\nabla u^{k}_{t}|^{2}\leq C_{\delta}\Phi_{K}^{10}\|\nabla u_{t}^{k}(t)\|_{L^{2}}^{2}+C_{\delta}\Phi_{K}^{42}.
	\end{aligned}
\end{equation*}
By Gronwall inequality, we have
\begin{equation*}
	\begin{aligned}
	\int_{0}^{t}\|\sqrt{\rho^{k}}	u^{k}_{tt}(s)ds\|^{2}_{L^{2}}+\|\nabla u^{k}_{t}(t)\|_{L^{2}}^{2}\leq C_{\delta}(1+\|V_{1}\|_{H^{1}}^{2})\exp(C\int_{0}^{t}\Phi_{K}^{42}(s)ds).
	\end{aligned}
\end{equation*}

\begin{Lemma}\label{pro-8}
	Assume that $(\rho^{k},u^{k},P^{k})$ is the unique global strong solution to the linearized problem $(\ref{37})$-$(\ref{3122})$ with $v$ replaced by $u^{k-1}$. Then, for any $t \in [0,T]$, it holds that	\begin{equation}\label{433}
			\left\|u^{k}\right\|_{H^{4}}+\left\|P^{k}\right\|_{H^{3}}\leq \widetilde{C} (\Phi_{K}^{28}+\|\sqrt{\rho^{k}}u^{k}_{tt}\|_{L^{2}}\Phi_{K}),
	\end{equation}
for all $k, 1 \leqslant k \leqslant K$.
\end{Lemma}
{\bf Proof.}
Recall that \begin{equation}
	\|u^{k}\|_{H^{4}}\leq \widetilde{C}\|F^{k}\|_{H^{2}}\Phi_{K}^{3}=\widetilde{C}(\|F^{k}\|_{W^{1,2}}+\|\nabla ^{2}F^{k}\|_{L^{2}})\Phi_{K}^{3},
\end{equation}
where $F^{k}=\rho^{k}f-\rho^{k}u^{k}_{t}+\rho^{k}u^{k-1}\partial_{x}u^{k}+\rho^{k}\int_{0}^{y}\partial_{x}u^{k-1}(s)ds\partial_{y}u^{k}$.\\

We only give the estimate of the term 
\begin{equation}
	\begin{aligned}
	\widetilde{C}	\|\nabla^{2}F^{k}\|_{L^{2}}\Phi_{K}^{3}, \end{aligned}\end{equation}where
	\begin{equation}\label{3.77}
	\begin{aligned} 
&\|\nabla^{2}F^{k}\|_{L^{2}}\\
&=\|\nabla^{2}\rho^{k}\partial_{t}u^{k}\|_{L^{2}}+\|\nabla\rho^{k}\nabla \partial_{t}u
^{k}\|_{L^{2}}+\|\nabla\rho^{k}\nabla \partial_{t}u
^{k}\|_{L^{2}}+\|\rho^{k}\nabla^{2}\partial_{t}u^{k}\|_{L^{2}}\\&\quad+\|\nabla^{2}\rho^{k}u^{k-1}\partial_{x}u^{k}\|_{L^{2}}+\|\nabla\rho^{k}\nabla u^{k-1}\partial_{x}u^{k}\|_{L^{2}}+\|\nabla \rho^{k}u^{k-1}\nabla\partial_{x}u^{k}\|_{L^{2}}\\&\quad
+\|\nabla^{2}u^{k-1}\partial_{x}u^{k}\rho^{k}\|_{L^{2}}+\|\nabla u^{k-1}\nabla \partial_{x}u^{k}\rho^{k}\|_{L^{2}}+\|\nabla u^{k-1}\partial_{x}u^{k}\nabla \rho^{k}\|_{L^{2}}\\&\quad+\|\nabla u^{k-1}\nabla \partial_{x}u^{k}\rho^{k}\|_{L^{2}}+\|u^{k-1}\nabla^{2}\partial_{x}u^{k}\rho^{k}\|_{L^{2}}+\|u^{k-1}\nabla \partial_{x}u^{k}\nabla \rho^{k}\|_{L^{2}}\\&\quad+\|\nabla^{2}\rho^{k}\int_{0}^{y}
\partial_{x}u^{k-1}(s)ds\partial_{y}u^{k}\|_{L^{2}}+\|\nabla \rho^{k}\nabla \int_{0}^{y}\partial_{x}u^{k-1}(s)ds\partial_{y}u^{k}\|_{L^{2}}\\&\quad+\|\nabla \rho^{k} \int_{0}^{y}\partial_{x}u^{k-1}(s)ds\nabla\partial_{y}u^{k}\|_{L^{2}}+\|\nabla \rho^{k}\int_{0}^{y}\partial_{x x}u^{k-1}(s)ds\partial_{y}u^{k}\|_{L^{2}}\\&\quad+\| \rho^{k}\nabla\int_{0}^{y}\partial_{x x}u^{k-1}(s)ds\partial_{y}u^{k}\|_{L^{2}}+\| \rho^{k}\int_{0}^{y}\partial_{x x}u^{k-1}(s)ds\nabla\partial_{y}u^{k}(s)ds\|_{L^{2}}\\&\quad+\|\nabla \rho^{k}\partial_{x}u^{k-1}\partial_{y }u^{k}\|_{L^{2}}+\| \rho^{k}\nabla\partial_{x}u^{k-1}\partial_{y }u^{k}\|_{L^{2}}+\| \rho^{k}\partial_{x}u^{k-1}\nabla\partial_{y }u^{k}\|_{L^{2}}+\nabla \rho^{k}\int_{0}^{y}\partial_{x}u^{k-1}(s)ds\nabla \partial_{y}u^{k}\|_{L^{2}}\\&\quad+\| \rho^{k}\nabla\int_{0}^{y}\partial_{x}u^{k-1}(s)ds\nabla \partial_{y}u^{k}\|_{L^{2}}+\| \rho^{k}\int_{0}^{y}\partial_{x}(s)dsu^{k-1}\nabla^{2} \partial_{y}u^{k}\|_{L^{2}}\\&=\sum_{i=1}^{25}O_{i}.
\end{aligned}
\end{equation}
We now provide estimates for the right-hand side terms of (\ref{3.77}).
By direct estimates, we have   
$$O_{1}\leq
\|\nabla^{2}\rho^{k}\|_{L^{\infty}}\|u_{t}^{k}\|_{L^{2}}.$$
$$O_{2}+O_{3}\leq\|\nabla\rho^{k}\|_{L^{\infty}}\|\nabla u^{k}_{t}\|_{L^{2}}.$$
$$O_{4}\leq\|\rho^{k}\|_{L^{\infty}}\|\nabla^{2}u_{t}^{k}\|_{L^{2}}.$$
$$O_{5}\leq\|\nabla^{2}\rho^{k}\|_{L^{\infty}}\|u^{k-1}\|_{L_{y}^{\infty}L_{x}^{2}}\|\partial_{x}u^{k}\|_{L_{x}^{2}L_{y}^{\infty}}.$$
\begin{equation}
	\begin{aligned}
		O_{6}+O_{10}+O_{20}\leq\|\nabla \rho^{k}\|_{L^{\infty}}\|\nabla u^{k-1}\|_{L_{x}^{2}L_{y}^{\infty}}\|\nabla u^{k}\|_{L_{y}^{2}L_{x}^{\infty}}.
	\end{aligned}
\end{equation}
\begin{equation*}
	\begin{aligned}
		O_{7}+O_{13}+O_{16}+O_{23}\leq \|\nabla \rho^{k}\|_{L^{\infty}}\|u^{k-1}\|_{H^{2}}\|\nabla^{2} u^{k}\|_{L^{2}}.
	\end{aligned}
\end{equation*}
\begin{equation*}
	\begin{aligned}
		O_{8}+O_{21}\leq\|\rho^{k}\|_{L^{\infty}}\|\nabla u^{k}\|_{L^{\infty}}\|\nabla ^{2}u^{k-1}\|_{L^{2}}.
	\end{aligned}
\end{equation*}
\begin{equation*}
\begin{aligned}
	O_{9}+O_{11}+O_{22}\leq\|\rho^{k}\|_{L^{\infty}}\|\nabla u^{k-1}\|_{L_{x}^{2}L_{y}^{\infty}}\|\nabla ^{2}u^{k}\|_{L_{y}^{2}L_{x}^{\infty}}.
\end{aligned}
\end{equation*}
\begin{equation*}
	\begin{aligned}
		O_{12}\leq\|\rho^{k}\|_{L^{\infty}}\|u^{k-1}\|_{L_{x}^{2}L_{y}^{\infty}}\|\nabla ^{2}\partial_{x}u^{k}\|_{L_{y}^{2}L_{x}^{\infty}}.
	\end{aligned}
\end{equation*}
\begin{equation*}
	\begin{aligned}
		O_{14}\leq \|\nabla^{2}\rho^{k}\|_{L^{\infty}}\|\int_{0}^{y}\partial_{x}u^{k-1}(s)ds\|_{L_{x}^{2}L_{y}^{\infty}}\|\partial_{y}u^{k}\|_{L_{x}^{\infty}L_{y}^{2}}.
	\end{aligned}
\end{equation*}
\begin{equation*}
	\begin{aligned}
	O_{15}+O_{17}&\leq \|\nabla \rho^{k}\|_{L^{\infty}}\|\partial_{x}u^{k-1}\|_{L_{x}^{2}L_{y}^{\infty}}\|\partial_{y}u^{k}\|_{L_{y}^{2}L_{x}^{\infty}}\\&\quad+\|\nabla\rho^{k}\|_{L^{\infty}}\|\int_{0}^{y}\partial_{x x}u^{k-1}(s)ds\|_{L_{x}^{2}L_{y}^{\infty}}\|\partial_{y}u^{k}\|_{L_{y}^{2}L_{x}^{\infty}}.
	\end{aligned}
\end{equation*}
\begin{equation*}
	\begin{aligned}
	O_{18}+O_{19}+O_{24}&\leq\|\rho^{k}\|_{L^{\infty}}\|\int_{0}^{y}\partial_{x xx}u^{k-1}(s)ds\|_{L_{x}^{2}L_{y}^{\infty}}\|\partial_{y}u^{k}\|_{L_{x}^{\infty}L_{y}^{2}}\\&\quad+\|\rho^{k}\|_{L^{\infty}}\|\int_{0}^{y}\partial_{ xx}u^{k-1}(s)ds\|_{L^{\infty}}\|\nabla\partial_{y}u^{k}\|_{L^{2}}\\&\leq C
	 \|u^{k-1}\|_{H^{3}}\|u^{k}\|_{H^{2}}.
	\end{aligned}
\end{equation*}
\begin{equation*}
	\begin{aligned}
	O_{25}&\leq \|\rho^{k}\|_{L^{\infty}}\|\int_{0}^{y}\partial_{x}u^{k-1}(s)ds\|_{L_{x}^{2}L_{y}^{\infty}}\|\nabla^{2}\partial_{y}u^{k}\|_{L_{y}^{2}L_{x}^{\infty}}\\&\leq\Phi_{K}(\|u^{k}\|_{H^{4}}^{\frac{1}{2}}\|u^{k}\|_{H^{3}}^{\frac{1}{2}}+\|u^{k}\|_{H^{3}}).
	\end{aligned}
\end{equation*}
Summing the above inequalities yields that
\begin{equation*}
	\begin{aligned}
\|u^{k}\|_{H{^4}}&=\widetilde{C}(\|F^{k}\|_{W^{1,2}}\Phi_{K}^{3}+\|\nabla^{2}F_{k}\|_{L^{2}}\Phi_{K}^{3})\\&\leq \widetilde{C}[\|\nabla u^{k-1}_{t}\|^{2}_{L^{2}}+ \|\nabla u^{k}_{t}\|^{2}_{L^{2}}+\|\rho^{k}\|_{L^{\infty}}\|\nabla^{2}u^{k}_{t}\|_{L^{2}}+\|f\|_{H^{1}} +\Phi_{K}^{28}]\\&\leq \widetilde{C}\Phi_{K}^{28}+\widetilde{C}\|\nabla^{2}u^{k}_{t}\|_{L^{2}}, \end{aligned}
	\end{equation*}
where we use the facts that 
$		\|u\|_{H^{2}}\leq C(1+\phi_{K}^{5})$  and $ \|u\|_{H^{3}}\leq C(1+\phi_{K}^{20}).$

Similar to (\ref{{3.44}}), we consider hydrostatic stokes equation: $$-\mathrm{div}(\nabla \mu(\rho^{k})\nabla u^{k}_{t})-\partial_{x }P_{t}^{k}=F^{k}_{t},$$
where\begin{equation*}
	\begin{aligned}F^{k}_{t}&=-\rho^{k}_{t}u^{k}_{t}-\rho^{k} u^{k}_{tt}-\rho^{k}_{t}u^{k-1}\partial_{x}u^{k}-\rho^{k} u^{k-1} \partial_{x}u_{t}^{k}-\rho^{k} u^{k-1}_{t}\partial_{x}u^{k}+\rho^{k}_{t}\int_{0}^{y}\partial_{x }u^{k-1}(s)ds\partial_{y}u^{k}\\&\quad+\rho^{k}\int_{0}^{y}\partial_{x }u_{t}^{k-1}(s)ds\partial_{y}u^{k}+\rho^{k}\int_{0}^{y}\partial_{x }u^{k-1}(s)ds\partial_{y}u_{t}^{k}+\mu(\rho^{k})_{t}\nabla u^{k}+\mu(\rho^{k})_{t}\Delta u^{k}\\&\quad+\rho^{k}_{t} f+\rho^{k} f_{t}.
	\end{aligned}
\end{equation*}
From Lemma
$\ref{Lemma-41}$ and (\ref{311}), we obtain
$$\|\nabla ^{2}u^{k}_{t}\|_{L^{2}}\leq C\Phi_{K}^{12}+C\|\sqrt{\rho^{k}}u^{k}_{tt}\|_{L^{2}}\Phi_{K}.$$
Moreover,
\begin{equation}
	\|u\|_{H^{4}}\leq \widetilde{C} (\Phi_{K}^{28}+\|\sqrt{\rho^{k}}u^{k}_{tt}\|_{L^{2}}\Phi_{K}).
\end{equation}
$\hfill\Box$\\
Now, we give the proof of Lemma \ref{Pro-3.6}.\\
\textbf{Proof of Lemma $\ref{Pro-3.6}$} It follows from $(\ref{2.6})$ and $(\ref{2.7})$ that
\begin{equation}
	\begin{aligned}
&	\| \nabla ^{2}\rho(t)\|_{L^{\infty
		}}+	\| \nabla ^{3}\rho(t)\|_{L^{2}}\\&\leq(\| \nabla ^{2}\rho_{0}\|_{W^{1,2}}+\|\nabla\rho_{0}\|_{W^{1,\infty}})
	\exp(C\int_{0}^{t}\|u(s)\|_{H^{4}}ds)
\\&\leq C	\exp(C\int_{0}^{t}\|u(s)\|_{H^{4}}ds)\\&\leq C\exp[\int_{0}^{t} \widetilde{C} (\Phi_{K}(s)^{28}+\|\sqrt{\rho^{k}}u^{k}_{tt}(s)\|_{L^{2}}^{2}) ds]\\&\leq C \exp[\widetilde{C}_{\delta}(1+\|V_{1}\|_{L^{2}})^{2}\exp(C\int_{0}^{t}\Phi_{K}(s)^{42}ds)]\exp \left[\widetilde{C}\int_{0}^{t} \Phi_{K}(s)^{28} \mathrm{~d} s\right]\\&\leq C\exp[\widetilde{C}_{\delta}(1+\|V_{1}\|_{L^{2}})^{2}\exp(\widetilde{C}\int_{0}^{t}\Phi_{K}(s)^{42}ds)],
\end{aligned}
\end{equation}
where $\widetilde{C}_{\delta}$ depending on $\widetilde{C}$ and $\delta$.
The proof of Lemma \ref{Pro-3.6} is complete.$\hfill\Box$

To get the uniform local boundness of $\Phi_{K}$ with respect to $K$, we combine $(\ref{42})$, $(\ref{433l})$, $(\ref{431})$ and  $(\ref{3.65})$-$(\ref{3.67})$ to  obtain
\begin{equation*}
	\Phi_{K}(t)  \leq C\exp[\widetilde{C}_{\delta}(1+\|V_{1}\|_{L^{2}})^{2}\exp(\widetilde{C}\int_{0}^{t}\Phi_{K}(s)^{42}ds)].
\end{equation*}

Define $\Psi_{K}(t)=\log(\widetilde{C}_{\delta}^{-1}\log[C^{-1}\Phi_{K}(t)]).$
Then it follows that
$$\Psi_{K}(t)\leq \log (1+\|V_{1}\|_{L^{2}})+\widetilde{C}\int_{0}^{t} e^{\widetilde{C}_{\delta}e^{\Psi_{K}(s)}}ds.$$

Thanks to Lemma $\ref{Lemma-2}$, there exists $T_{0}\in (0,T)$ depending only on $\widetilde{C}_{\delta}$ such that
$$\sup _{0 \leqslant t \leqslant T_{0}} \Phi_{K}(t) \leqslant C \exp \left[\widetilde{C}_{\delta}(1+\|V_{1}\|_{L^{2}})\right].$$

Recall that $K$ is a fixed large integer. Thanks to Lemma $\ref{Pro-1}$-$\ref{pro-8}$, we deduce the following uniform bounds of the approximate solutions:
\begin{equation}\label{541}
	\begin{aligned}
		&\sup_{0 \leqslant t \leqslant T_{0}}\left(\left\|u^{k}\right\|_{H^{2}}+\left\|P^{k}\right\|_{H^{1}}+\left\| \sqrt{\rho^{k} } u^{k}_{t}\right\|_{L^{2}}+\left\| \rho^{k}\right \|_{W^{2, \infty}}+\left\| \rho^{k}\right \|_{W^{3, 2}}+\|u^{k}\|_{H^{3}}+\|\nabla u^{k}_{t}\|_{L^{2}}\right) \\
		&+\int_{0}^{T_{0}}\left(\left\|u_{t}^{k}\right\|_{H^{2}}^{2}+\left\|u^{k}\right\|_{H^{4}}^{2}+\left\|P^{k}\right\|_{H^{2}}^{2}\right) \mathrm{d} t \leqslant C\exp \{C \exp(\widetilde{C}_{\delta}(1+\|V_{1}\|_{L^{2}}))  \}
	\end{aligned}
\end{equation}
for all $k\geqslant1$.
\subsubsection{Convergence}
Based on uniform estimates of (\ref{541}),
we prove that the approximate solution $(\rho^{k},u^{k}, P^{k})$ convergences to a solution to (\ref{3.1}) in a sufficiently strong sense.
Define $$\sigma^{k+1}=\rho^{k+1}-\rho^{k},~\eta^{k+1}=u^{k+1}-u^{k} ~\text{and}  ~-\int_{0}^{y}\partial_{x}\eta^{k+1}(s)ds=-\int_{0}^{y}\partial_{x}u^{k+1}(s)ds+\int_{0}^{y}\partial_{x}u^{k}(s)ds.$$
Then it follows from the linearized momentum equation $(\ref{37})_2$ that
\begin{equation}\label{f}
	\begin{aligned}
	&\rho^{k+1}\partial_{t}\eta^{k+1}+\rho^{k+1}u^{k}\partial_{x}\eta^{k+1}-\rho^{k+1}\int_{0}^{y}\partial_{x}u^{k}(s)ds\partial_{y}\eta^{k+1}+\partial_{x}(P^{k+1}-P^{k})-\operatorname{div}[\mu(\rho^{k+1})\nabla \eta^{k+1}]\\
	&=\sigma^{k+1}(f-\partial_{t}u^{k}-u^{k}\partial_{x}u^{k}+\int_{0}^{y}\partial_{x}u^{k}(s)ds\partial_{y}u^{k})-\rho^{k}\eta^{k}\partial_{x}u^{k}+\rho^{k}\int_{0}^{y}\partial_{x}\eta^{k}(s)ds\partial_{y}u^{k}\\&+\partial_{x}[(\mu(\rho^{k+1})-\mu(\rho^{k}))\partial_{x}u^{k}]+\partial_{y}[(\mu(\rho^{k+1})-\mu(\rho^{k}))\partial_{y}u^{k}].
	\end{aligned}
\end{equation}
Hence multiplying (\ref{f}) by $\eta^{k+1}$, integrating over $\Omega$ and using  $(\ref{37})_1$, we obtain
\begin{equation}\label{g}
	\frac{d}{dt}\int_{\Omega}\rho^{k+1}\frac{|\eta^{k+1}|^{2}}{2}+\int_{\Omega}\mu^{k+1}|\nabla \eta^{k+1}|^{2}=\int_{\Omega}[\sigma^{k+1}(f-\partial_{t}u^{k}-u^{k}\partial_{x}u^{k}+\int_{0}^{y}\partial_{x}u^{k}(s)ds\partial_{y}u^{k})\eta^{k+1}$$$$-\rho^{k}\eta^{k}\partial_{x}u^{k}\eta^{k+1}+\rho^{k}\int_{0}^{y}
	\partial_{x }\eta^{k}(s)ds\partial_{y}u^{k}\eta^{k+1}+\operatorname{div}[(\mu^{k
		+1}-\mu^{k})\nabla u^{k}]\eta^{k+1}+\lambda\rho^{k}_{xx}|\eta^{k+1}|^{2}].
\end{equation}

Using H$\ddot{o}$lder inequalities, Gagliardo-Nirenberg inequalities,  Young inequalities and  the fact that $C^{-1} \leqslant \mu^{k} \leqslant C$ and \textcolor{red}{$\lambda$~}is a small enough number, we deduce that

\begin{equation}\label{441}
	\begin{aligned}
		\frac{d}{dt}\int_{\Omega}\rho^{k+1}|\eta^{k+1}|^{2}+\frac{1}{2C}\int_{\Omega}|\nabla \eta^{k+1}|^{2}&\leq C\|\sigma^{k+1}\|_{L^{2}}^{2}(\|f-\partial_{t}u^{k}-u^{k}\partial_{x}u^{k}\\&\quad+\int_{0}^{y}\partial_{x}u^{k}(s)ds\partial_{y}u^{k}\|_{L^{\infty}_{y}L^{2}_{x}}^{2}+\|\nabla u^{k}\|_{L^{\infty}}^{2})+C\|\sqrt{\rho^{k}}\eta^{k}\|_{L^{2}}^{2}\\&=:C\|\sigma^{k+1}\|_{L^{2}}^{2}B_{k}(t)+C\|\sqrt{\rho^{k}}\eta^{k}\|_{L^{2}}^{2}
	\end{aligned}
\end{equation}
where $ B_{k}(t)=\|f-\partial_{t}u^{k}-u^{k}\partial_{x}u^{k}+\int_{0}^{y}\partial_{x}u^{k}(s)ds\partial_{y}u^{k}\|_{L^{\infty}_{y}L^{2}_{x}}^{2}+\|\nabla u^{k}\|_{L^{\infty}}^{2}$ and we have uesd $$\|\eta^{k+1}\|_{L^{2}_{y}L^{\infty}_{x}}\leq \|\partial_{x}\eta^{k
	+1}\|_{L_{x}^{2}L_{y}^{2}}^{\frac{1}{2}}\|\eta^{k
	+1}\|_{L_{x}^{2}L_{y}^{2}}^{\frac{1}{2}}+\|\eta^{k
	+1}\|_{L_{x}^{2}L_{y}^{2}}\leq \|\nabla\eta^{k+1}\|_{L^{2}}.$$

Now, we deal with the term $ \|\sigma^{k+1}\|_{L^{2}}^{2}$.
It follows from  $(\ref{37})_1$ that
\begin{equation}\label{h}
\sigma^{k+1}_{t}+u^{k}\partial_{x}\sigma^{k+1}+\eta^{k}\partial_{x}\rho^{k}-\int_{0}^{y}\partial_{x}u^{k}(s)ds\partial_{y}\sigma^{k+1}-\int_{0}^{y}\partial_{x}\eta^{k}(s)ds\partial_{x}\rho^{k}-\lambda\sigma^{k+1}_{xx}=0.
\end{equation}
Hence multiplying (\ref{h}) by $\sigma^{k+1}$, integrating over $\Omega$, we obtain
\begin{equation*}
	\begin{aligned}
	&	\frac{d}{dt}\int_{\Omega}(\sigma^{k+1})^{2}+\lambda\|\sigma^{k+1}_{x}\|^{2}_{L^{2}}\\&\leq\|\nabla \rho^{k}\|_{W^{1,\infty}}\|\eta^{k}\|_{L^{2}}\|\sigma^{k+1}\|_{L^{2}}+\|\nabla \rho^{k}\|_{L^{\infty}}\| \eta^{k}\|_{L^{2}}\|\sigma_{x}^{k+1}\|_{L^{2}}\\&\leq C \|\sqrt{\rho^{k}}\eta^{k}\|_{L^{2}}\|\sigma^{k+1}\|_{L^{2}}+C\|\sqrt{\rho^{k}} \eta^{k}\|_{L^{2}}\|\sigma_{x}^{k+1}\|_{L^{2}}\\&\leq C \|\sqrt{\rho^{k}}\eta^{k}\|_{L^{2}}^{2}+C	\|\sigma^{k+1}\|_{L^{2}}^{2}+\epsilon	\|\sigma_{x}^{k+1}\|_{L^{2}}^{2}.
	\end{aligned}
\end{equation*}
\\
Since $\sigma^{k+1}(0)=0$, Gronwall's inequality yields
\begin{equation}\label{4422}
	\|\sigma^{k+1}\|_{L^{2}}^{2}\leq  C \int_{0}^{t}\|\sqrt{\rho^{k}}\eta^{k}(s)\|_{L^{2}}^{2}ds.
\end{equation}

Substituting (\ref{4422})
into (\ref{441}) and integrating (\ref{441}) over $(0,t)$, we obtain\begin{equation}\label{d}
\begin{aligned}
&	\|\sqrt{\rho^{k+1}}\eta^{k+1}(t)\|_{L^{2}}^{2}+\int_{0}^{t}\int_{\Omega}|\nabla \eta^{k+1}(s)|^{2}ds\\&\leq \int_{0}^{t}[C \int_{0}^{s}\|\sqrt{\rho^{k}}\eta^{k}(\tau)\|_{L^{2}}^{2} d\tau B_{k}(s)]ds+C\int_{0}^{t}\|\sqrt{\rho^{k}}\eta^{k}(s)\|_{L^{2}}^{2}ds
		\\&\leq C \int_{0}^{t}\|\sqrt{\rho^{k}}\eta^{k}(\tau)\|_{L^{2}}^{2}d\tau \int_{0}^{t}B_{k}(s)ds+C\int_{0}^{t}\|\sqrt{\rho^{k}}\eta^{k}(s)\|_{L^{2}}^{2}ds
		\\&\leq C\int_{0}^{t}\|\sqrt{\rho^{k}}\eta^{k}(s)\|_{L^{2}}^{2}ds ,
\end{aligned}
\end{equation}
where we have used the fact that $\int_{0}^{t}B_{k}(s)ds\leq C$ for all $k$.
Now, we only check the nonlinear terms of $B_{k}(t)$ and the boundness of other terms directly follows from the regularity (\ref{541}).
For the nonlinear term, one has
$$\|u^{k}\partial_{x} u^{k}\|_{L^{\infty}_{y}L^{2}_{x}}^{2}
\leq C\|u^{k}\|_{H^{2}}^{4},
$$
and
$$\|\int_{0}^{y}\partial_{x}u^{k}(s)ds\partial_{y}u^{k}\|_{L^{\infty}_{y}L^{2}_{x}}^{2}\leq C \|u^{k}\|_{H^{2}}^{4}.$$
Hence, it follows from the regularity (\ref{541}) that $\int_{0}^{t}B_{k}(s)ds\leq C$ for all k.

Finally, a simple recursive argument shows that\begin{equation}
\sup_{0 \leqslant t \leqslant T_{0}}\|\sqrt{\rho^{k+1}}\eta^{k+1}(t)\|_{L^{2}}^{2}\leq C\frac{{T_{0}}^{k-1}}{(k-1)!}\int_{0}^{T_{0}}\|\sqrt{\rho^{1}}\eta^{1}(t)\|_{L^{2}}^{2}.
\end{equation}

Next, thanks to (\ref{541}) and Aubin-Lions Lemma, it is easy to show the existence of subsequences, denoted again by $(u^{k+1},u^{k},\rho^{k+1},P^{k+1})$, such that
$$u^{k+1}, u^{k} \rightarrow u ~\text{in}~ L^{2}(0,T_{0}; H^{1}(\Omega),$$
$$\rho^{k+1} \rightarrow \rho ~\text{in}~ L^{\infty}(0,T_{0};L^{2}(\Omega)),$$
$$P^{k+1}\rightarrow P ~\text{weakly}~ ~\text{in}~L^{2}(0,T_{0};L^{2}(\Omega)).$$
Thus it is a simple matter to check that the limit $(\rho, u, P)$ ia a weak solution to regularization problems $(\ref{3.1})$.

 Combining with  $(\ref{541})$, we prove that $(\rho,u,P)$  is a strong solution to regularization problems $(\ref{3.1})$, satisfying
 \begin{equation}\label{546}
\begin{aligned}
&\sup_{0 \leqslant t \leqslant T_{0}}\left(\left\|u\right\|_{H^{2}}+\left\|P\right\|_{H^{1}}+\left\| \sqrt{\rho} u_{t}\right\|_{L^{2}}+\left\| \rho\right \|_{W^{2, \infty}}+\left\| \rho\right \|_{H^{3}}+\|u\|_{H^{3}}+\|\nabla u_{t}\|_{L^{2}}\right) \\
&+\int_{0}^{T_{0}}\left(\left\|u_{t}\right\|_{H^{1}}^{2}+\left\|u\right\|_{H^{4}}^{2}+\left\|P\right\|_{H^{2}}^{2}\right) \mathrm{d} t \leqslant C\exp \{C \exp(\widetilde{C}_{\delta}+\widetilde{C}_{\delta}\|V_{1}\|_{L^{2}})  \}.
\end{aligned}
\end{equation}
\subsubsection{Uniqueness}
Next, we prove the uniqueness of the strong solution to the system (\ref{3.1}).

Suppose that there exists  solutions $(\bar{u}, \bar{\rho})$ and $(u, \rho)$ satisfying $(\ref{3.1})$. Then we have
 \begin{equation}
 \begin{aligned}
 &\bar{\rho}\partial_{t}(\bar{u}-u)+\bar{\rho}\bar{u}\partial_{x}(\bar{u}-u)-\bar{\rho}\int_{0}^{y}\partial_{x}\bar{u}(s)ds\partial_{y}(\bar{u}-u)\\
 & +\partial_{x}(\bar{p}-p)-\partial_{x}\big[\mu(\bar{\rho})\partial_{x}(\bar{u}-u)\big]-\partial_{y}\big[\mu(\bar{\rho})\partial_{y}(\ \bar{u}-u)\big]\\
 &=(\bar{\rho}-\rho)(f-u_{t})+(\rho-\bar{\rho})u\partial_{x}u-(\rho-\bar{\rho})\partial_{y}u\int_{0}^{y}\partial_{x}u(s)ds\\
 &\quad+\partial_{x}\big[(\mu(\bar{\rho})-\mu(\rho))\partial_{x}u\big]+\partial_{y}\big[(\mu(\bar{\rho})-\mu(\rho))\partial_{y}u\big]
 -\bar{\rho}(\bar{u}-u)\partial_{x}u-\bar{\rho}\int_{0}^{y}\partial_{x}(\bar{u}-u)ds\partial_{y}u.
 \end{aligned}
 \end{equation}

Multiplying $\bar{u}-u$ and using H$\ddot{o}$lder and Young's inequality, we have
\begin{equation}
\begin{aligned}
&\frac{d}{dt}\int_{\Omega}\bar{\rho}(\bar{u}-u)^{2}+\int_{\Omega}\mu(\bar{\rho})|\nabla(\bar{u}-u)|^{2}\\
&\leq\int_{\Omega}[(\bar{\rho}-\rho)(f-u_{t})(\bar{u}-u)+(\rho-\bar{\rho})u\partial_{x}u(\bar{u}-u)-(\rho-\bar{\rho})\partial_{y}u\int_{0}^{y}\partial_{x}u(s)ds(\bar{u}-u)\\
&\quad+\partial_{x}(\bar{u}-u)(\mu(\bar{\rho})-\mu(\rho))\partial_{x}u+\partial_{y}(\bar{u}-u)(\mu(\bar{\rho})-\mu(\rho))\partial_{y}u\\
&\quad-\bar{\rho}(\bar{u}-u)^{2}\partial_{x}u-\bar{\rho}\int_{0}^{y}\partial_{x}(\bar{u}-u)ds\partial_{y}u(\bar{u}-u)]\\
&\leq\left\|\bar{\rho}-\rho\right\|_{L^{2}}\left\|\bar{u}-u\right\|_{L^{\infty}_{y}L^{2}_{x}}\left\|f-u_{t}\right\|_{L^{2}_{y}L^{\infty}_{x}}\quad+
\left\|\bar{\rho}-\rho\right\|_{L^{2}}\left\|\bar{u}-u\right\|_{L^{2}}\left\|\partial_{x}u\right\|_{L^{\infty}}\left\|u\right\|_{L^{\infty}}\\
&\quad+\left\|\bar{\rho}-\rho\right\|_{L^{2}}\left\|\bar{u}-u\right\|_{L^{2}_{y}L^{\infty}_{x}}\left\|\partial_{y}u\right\|_{L^{\infty}}\left\|\int_{0}^{y}
\partial_{x}u(s)ds\right\|_{L^{\infty}_{y}L^{2}_{x}}+\left\|\bar{\rho}-\rho\right\|_{L^{2}}\left\|\nabla u\right\|_{L^{\infty}}
\left\|\nabla (\bar{u}-u)\right\|_{L^{2}}\\
&\quad+\left\|\sqrt{\bar{\rho}}(\bar{u}-u)\right\|_{L^{2}}\left\|\nabla (\bar{u}-u)\right\|_{L^{2}}\left\|u\right\|_{H^{2}}\\
&\leq C\left\|\bar{\rho}-\rho\right\|_{L^{2}}^{2}\left\|\nabla f-\nabla u_{t}\right\|_{L^{2}}^{2}+\epsilon\left\|\nabla (\bar{u}-u)\right\|_{L^{2}}^{2}+C\left\|\bar{\rho}-\rho\right\|_{L^{2}}^{2}\left\|u\right\|_{H^{3}}^{2}+\epsilon\left\|\nabla (\bar{u}-u)\right\|_{L^{2}}^{2}\\
&\quad+C\left\|\sqrt{\bar{\rho}}(\bar{u}-u)\right\|_{L^{2}}^{2}+\epsilon\left\|\nabla (\bar{u}-u)\right\|_{L^{2}}^{2}.
\end{aligned}
\end{equation}
Here and below, $\epsilon$ is some small number.\\
 We use the fact that $\underline {\mu}\leq\mu$ and $0<\epsilon<1$, we have
\begin{equation}\label{449}
\begin{aligned}
&\frac{d}{dt}\int_{\Omega}\bar{\rho}(\bar{u}-u)^{2}+\int_{\Omega}|\nabla(\bar{u}-u)|^{2}\\
&\leq C\left\|\bar{\rho}-\rho\right\|_{L^{2}}^{2}(\left\|\nabla f-\nabla
u_{t}\right\|_{L^{2}}^{2}+\left\|u\right\|_{H^{3}}^{2})+C\left\|\sqrt{\bar{\rho}}(\bar{u}-u)\right\|_{L^{2}}^{2}.
\end{aligned}
\end{equation}

Next, it follows from the first equation of ($\ref{a1}$) that
\begin{equation*}
	\begin{aligned}
		&\partial_{t}(\rho-\bar{\rho})+(u-\bar{u})\partial_{x}\rho-\int_{0}^{y}\partial_{x}(u-\bar{u})(s)ds\partial_{y}\rho+\bar{u}(\partial_{x}\rho-\partial_{x}\bar{\rho})\\&
		-\int_{0}^{y}\partial_{x}\bar{u}(s)ds(\partial_{y}\rho-\partial_{y}\bar{\rho})-\lambda\partial_{ xx}(\rho-\bar{\rho})=0.
	\end{aligned}
\end{equation*}

Mutiplying $\rho-\bar{\rho}$, we have $$
\frac{1}{2}\frac{d}{dt}\int_{\Omega}(\rho-\bar{\rho})^{2}+\int_{\Omega}\lambda|\partial_{x}(\rho-\bar{\rho})|^{2}+\int_{\Omega}(u-\bar{u})\partial_{x}\rho(\rho-\bar{\rho})-\int_{\Omega}\int_{0}^{y}\partial_{x}(u-\bar{u})(s)ds\partial_{y}\rho(\rho-\bar{\rho})=0
.$$

Using H$\ddot{o}$lder inequality, one has
\begin{equation}\label{451}
\begin{aligned}
\frac{1}{2}\frac{d}{dt}\int_{\Omega}(\rho-\bar{\rho})^{2}&\leq\left\|u-\bar{u}\right\|_{L^{2}}\left\|\rho-\bar{\rho}\right\|_{L^{2}}+\left\|\nabla(\bar{u}-u)\right\|_{L^{2}}\left\|\rho-\bar{\rho}\right\|_{L^{2}}\\
&\leq \epsilon\left\|\nabla(u-\bar{u})\right\|_{L^{2}}^{2}+C\left\|\rho-\bar{\rho}\right\|_{L^{2}}^{2}.
\end{aligned}
\end{equation}

Therefore, we deduce
\begin{equation}\label{452}
\begin{aligned}
&\frac{d}{dt}\int_{\Omega}\bar{\rho}(\bar{u}-u)^{2}+\frac{1}{2}\frac{d}{dt}\int_{\Omega}(\rho-\bar{\rho})^{2}+\int_{\Omega}|\nabla(\bar{u}-u)|^{2}\\
&\leq C\left\|\bar{\rho}-\rho\right\|_{L^{2}}^{2}(\left\|\nabla f-\nabla
u_{t}\right\|_{L^{2}}^{2}+\left\|u\right\|_{H^{3}}^{2}+1)+C\left\|\sqrt{\bar{\rho}}(\bar{u}-u)\right\|_{L^{2}}^{2}\\
&\leq C(\left\|\bar{\rho}-\rho\right\|_{L^{2}}^{2}+\left\|\sqrt{\bar{\rho}}(\bar{u}-u)\right\|_{L^{2}}^{2})(\left\|\nabla f-\nabla u_{t}\right\|_{L^{2}}^{2}+\left\|u\right\|_{H^{3}}^{2}+1).
\end{aligned}
\end{equation}

Combining $(\ref{546})$ and Gronwall inequality, one obtains $$
\left\|\rho-\bar{\rho}\right\|_{L^{2}}^{2}+\left\|\sqrt{\bar{\rho}}(\bar{u}-u)\right\|_{L^{2}}^{2}\leq0.$$

Hence, we have $\rho=\bar{\rho}$, and $u=\bar{u}$.
This completes the existence and uniqueness for nonnegative initial density.

\section{Proof of Theorem \ref{the-1}}
\subsection{Existence and uniqueness}
\  \  \
First, let $\lambda\rightarrow 0$ in systems (\ref{3.1}). We can get a unique strong solution of original problems $(\ref{a1})$-$(\ref{14})$ satisfying the initial assumption in Proposition \ref{pro-1}.

Now we prove the existence and uniqueness of Theorem \ref{the-1} for the general case, let $(\rho_{0},u_{0})$ be an initial data satisfying the hypotheses of Theorem \ref{the-1}. 	For each $\delta\in(0,1)$, choose $\rho^{\delta}_{0} \in W^{3,2}, \nabla \rho^{\delta}_{0} \in W^{1,\infty}, u_{0}^{\delta}\in H^{3} $ and $\mu^{\delta}\in C^{3}[0,\infty)$ such that 
	\begin{equation*}\begin{aligned}
			&	0<\delta \leq \rho^{\delta}_{0}\leq \rho_{0}+1,~\rho_{0}^{\delta}\rightarrow \rho_{0}~ \text{in}~ W^{2,2},\nabla \rho_{0}^{\delta}\xrightarrow{*}\nabla\rho_{0}~\text{in}~L^{\infty
			},\\&~u^{\delta}_{0}~\rightarrow u_{0}~ \text{in}~ H^{2} ~\text{and}~\mu^{\delta}\rightarrow\mu ~\text{in}~ C^{2}[0,\infty),
		\end{aligned}	
	\end{equation*} and denote by $(u_{0}^{\delta},P_{0}^{\delta})\in H^{1}_{0}\times L^{2}$ a solution to problem $(\ref{1.5})$. Then, according to Proposition \ref{pro-1}, there is a unique strong solution $(\rho^{\delta},u^{\delta}, P^{\delta}).$ We will show that a subsequence of the approximate solutions $u^{\delta}$ converges to a solution of the original problem. To show this, we need to derive some uniform bounds.  The regularity of solutions follows with similar arguments as above. We will only give the main estimates.
We introduce a functional $\Phi_{}(t)$ defined by
\begin{equation}
	\Phi(t)=1+\left\|\nabla \rho^{\delta}(t)\right\|_{L^{\infty}}+\left\|\nabla	^{2}\rho^{\delta}(t)\right\|_{L^{2}}+\left\|\nabla u^{\delta}(t)\right\|_{L^{2}},
\end{equation}

\begin{Lemma}\label{4.1}
	Assume that $(\rho^{\delta},u^{\delta},P^{\delta})$ is the unique local strong solution to the regularization  problem $(\ref{3.1})$  on $[0,T_{0}]$. Then, for any $t \in [0,T_{0}]$, it holds that
	\begin{equation}\label{42}
		\int_{0}^{t}\left\|\sqrt{\rho^{\delta}} u_{t}^{\delta
		}(s)\right\|_{L^{2}}^{2} \mathrm{~d} s+\left\|\nabla u^{\delta}(t)\right\|_{L^{2}}^{2} \leqslant C+C\int_{0}^{t} \Phi(s)^{8} \mathrm{~d} s.
	\end{equation}

\end{Lemma}
\begin{Lemma}\label{4.2}
	Assume that $(\rho^{\delta},u^{\delta},P^{\delta})$ is the unique local strong solution to the regularization  problem $(\ref{3.1})$  on $[0,T_{0}]$. Then, for any $t \in [0,T_{0}]$, it holds that
	\begin{equation}
		\left\|\sqrt{\rho^{\delta}} u_{t}^{\delta}\right\|_{L^{2}}^{2}+\int_{0}^{t}\left\|\nabla u_{t}^{\delta}(s)\right\|_{L^{2}}^{2} \mathrm{~d} s \leqslant C(1+\|V_{1}\|_{L^{2}})\exp \left[C\int_{0}^{t} \Phi(s)^{16} \mathrm{~d} s\right]
	\end{equation}

\end{Lemma}
\begin{Lemma}\label{4.3}
		Assume that $(\rho^{\delta},u^{\delta},P^{\delta})$ is the unique local strong solution to the regularization  problem $(\ref{3.1})$  on $[0,T_{0}]$. Then, for any $t \in [0,T_{0}]$, it holds that
	\begin{equation}\label{433l}
		\left\|\nabla \rho^{\delta}(t)\right\|_{L^{\infty}} \leqslant C \exp \left[C\left(1+\|V_{1}\|_{L^{2}}\right) \exp \left(C \int_{0}^{t} \Phi(s)^{20} \mathrm{~d} s\right)\right],
	\end{equation}
	and
	\begin{equation}\label{431}
		\left\|\nabla^{2} \rho^{\delta}(t)\right\|_{L^{2}} \leqslant C\exp \left[C\left(1+\|V_{1}\|_{L^{2}}\right) \exp \left(C\int_{0}^{t} \Phi(s)^{20} \mathrm{~d} s\right)\right].
	\end{equation}
\end{Lemma}
Now, combining Lemma $\ref{4.1}$-$\ref{4.3}$, we have
$$\Phi(t) \leqslant C \exp \left[C\left(1+\|V_{1}\|_{L^{2}}\right) \exp \left(\int_{0}^{t} \Phi(s)^{20} \mathrm{~d} s\right)\right].$$

Define $\Psi(t)=\log(C^{-1}\log[C^{-1}\Phi(t)]).$
Then it follows that
$$\Psi(t)\leq \log (1+\|V_{1}\|_{L^{2}})+C\int_{0}^{t} Ce^{Ce^{\Psi(s)}}ds.$$

Thanks to this integral inequality, there exists $T_{*}\in (0,T_{0})$ depending only on $\|V_{1}\|_{L^{2}}$ and $C$ such that
$$\sup _{0 \leqslant t \leqslant T_{*}} \Phi(t) \leqslant C \exp \left(C+ C\|V_{1}\|_{L^{2}}\right).$$
See  the proof of Lemma 6 in \cite{simon}.
We get the uniform local boundness of $\Phi(t)$ with respect to $\delta$. Thanks to Lemma \ref{Lemma-41}, we deduce the following uniform bounds of the approximate solutions:
\begin{equation}\label{4.6}
	\begin{aligned}
		&\sup_{0 \leqslant t \leqslant T_{*}}\left(\left\|u^{\delta}\right\|_{H^{2}}+\left\|P^{\delta}\right\|_{H^{1}}+\left\| \sqrt{\rho^{\delta}} u^{\delta}_{t}\right\|_{L^{2}}+\left\| \rho^{\delta}\right \|_{W^{1, \infty}}+\left\| \rho^{\delta}\right \|_{W^{2, 2}}\right) \\&+\int_{0}^{T_{0}}\left(\left\|u_{t}^{\delta}\right\|_{H^{1}}^{2}+\left\|u^{\delta}\right\|_{H^{3}}^{2}+\left\|P^{\delta}\right\|_{H^{2}}^{2}\right) \mathrm{d} t \leqslant C\exp \left[C \exp \left(C+ C\|V_{1}\|_{L^{2}}\right)\right].	\end{aligned}
\end{equation}

 Since the corresponding solutions $(\rho^{\delta},u^{\delta}, P^{\delta})$ satisfy the bound $(\ref{4.6})$ with $\|V_{1}\|_{L^{2}} $ and the constants $T_{*}$ and $C$ are independent of $\delta$, we can choose a subsequence of solutions $(\rho^{\delta},u^{\delta},P^{\delta})$ which converges to a limit $(\rho,u, P)$ in a weak sense. Obviously, it is a strong solution to the original problem satisfying the regularity estimate $(\ref{4.6})$. The uniqueness follows from the above proof similarly. Hence, we prove the existence and uniqueness of strong solutions of the general case to the system $(\ref{a1})$-$(\ref{14})$.
\subsection{Continuity}
\ \ \
In this subsection, we prove the continuity in time of $(\rho, u, P)$.
Firstly, we consider the continuity of $\rho$. By  the first equation of $(\ref{a1})$ and $(\ref{4.6})$, we deduce that
$$\rho \in L^{\infty}(0,T_{*};W^{1,\infty}) ,\quad \rho_{t}\in L^{\infty}(0,T_{*};L^{q}),$$
where $1\leq q<\infty$, which implies that
 $\rho \in C([0,T_{*}];W^{1,\infty}$-weak).  On the other hand, using $(\ref{2.4})$, one has
$$\lim \sup_{t \rightarrow+0}\|\rho(t)\|_{W^{1, \infty}} \leqslant\|\rho_{0}\|_{W^{1, \infty}}.$$

 It follows that $\lim_{t\rightarrow+0}\left\|\rho(t)-\rho_{0}\right\|_{W^{1,\infty}}=0$, that is, $\rho$ is left-continuous at $t=0$ (see $\cite{galdi}$). Because the continuity equation is time reversible, it yields that $\rho\in C([0,T_{*}];W^{1,\infty})$.
Similarly, we can obtain $\rho\in C([0,T_{*}];W^{2,2})$.

Secondly, we show the continuity in time of $(u,P)$. It is usual to prove that (see \cite{temam})
$$ u\in C([0,T_{*}];W^{1,p})\cap C([0,T_{*}];H^{2}-weak)$$
for $1\leq p<\infty$.

We now show that $(u,P)\in C([0,T_{*}];H^{2}\times H^{1})$. Observe that for each $t\in[0,T_{*}]$, $u=u(t)\in H^{1}_{0,per}\cap H^{2}$ is a solution of the hydrostatic Stokes equations
\begin{equation}\label{453}
-\partial_{x}\left(\mu(\rho)\partial_{x}u\right)-\partial_{y}\left(\mu(\rho)\partial_{y}u\right)+\partial_{x}P=F \quad \text{and}\quad \partial_{y}P =0\quad \text  { in}
\quad \Omega,
\end{equation}
where $F=\rho f-\rho\partial_{t}u+\rho u\partial_{x}u+\rho\int_{0}^{y}\partial_{x}u(s)ds\partial_{y}u .$

It is required to show that $F\in C([0,T_{*}],L^{2})$. From $(\ref{444})$ in Remark $\ref{remark-1}$, we deduce that the function $t \mapsto\left\|\sqrt{\rho} u_{t}(t)\right\|_{L^{2}}^{2}$ is continuous on $[0,T_{*}]$. Therefore, recalling that $\rho\in C([0,T_{*}];W^{1,\infty})$, we conclude that $\rho u_{t}\in C([0,T_{*}],L^{2})$.

It follows that
\begin{equation}
\begin{aligned}
&\left\|\rho(t)\int_{0}^{y}\partial_{x}u(t,\tau)d\tau\partial_{y}u(t)-\rho(s)\int_{0}^{y}\partial_{x}u(s,\tau)d\tau\partial_{y}u(s)\right\|_{L^{2}}\\
&\leq\left\|\rho(t)-\rho(s)\right\|_{L^{\infty}}\left\|\int_{0}^{y}\partial_{x}u(t,\tau)d\tau\right\|_{L^{\infty}_{y}L^{2}_{x}}\left\|\partial_{y}u(t)\right\|_{L^{\infty}_{x}L^{2}_{y}}\\
&\quad+\left\|\rho(s)\right\|_{L^{\infty}}\left\|\int_{0}^{y}\partial_{x}u(t,\tau)d\tau-\int_{0}^{y}\partial_{x}u(s,\tau)d\tau\right\|_{L^{\infty}_{y}L^{2}_{x}}
\left\|\partial_{y}u(t)\right\|_{L^{\infty}_{x}L^{2}_{y}}\\&\quad+\left\|\rho(s)\right\|_{L^{\infty}}\left\|\int_{0}^{y}\partial_{x}u(t,\tau)d\tau\right\|_{L^{\infty}}
\left\|\partial_{y}u(t)-\partial_{y}u(s)\right\|_{L^{2}} \rightarrow0\quad\text{ as}\quad t\rightarrow s.
\end{aligned}
\end{equation}

The term $\rho u\partial_{x}u$ has the similar argument. Hence, we obtain that $F \in C([0,T_{*}];L^{2})$.

Similar to the estimate $(\ref{3333})$, we easily deduce that for any $t,s\in[0,T_{*}]$,
$$\|P(t)-P(s)\|_{L^{2}} \leqslant C\left(\|\mu(\rho(t)) \nabla u(t)-\mu(\rho(s)) \nabla u(s)\|_{L^{2}}+\|F(t)-F(s)\|_{L^{2}}\right),$$
which implies, in particular, that $P\in C([0,T_{*}];L^{2})$.

Rewrite $(\ref{453})$ as
$$-\Delta u+\partial_{x}\tilde{P}=\mu^{-1}(F+ \nabla \mu \nabla u+\tilde{P} \partial_{x} \mu(\rho)),$$
where $\tilde{P}=P/\mu$, $\mu=\mu(\rho)$.

Similar to the estimate $(\ref{34})$, we can show that
\begin{equation}
\begin{aligned}
&\left\|u(t)-u(s)\right\|_{H^{2}}+\left\|\tilde{P}(t)-\tilde{P}(s)\right\|_{H^{1}}\\&\leq C \left\|F(t)-F(s)\right\|_{L^{2}}+\|\nabla\mu(\rho(t))\nabla u(t)-\nabla\mu(\rho(s))\nabla u(s)\|_{L^{2}}\\&\quad+\|\tilde{P}(t)\partial_{x}\mu(\rho(t))-\tilde{P}(s)\partial_{x}\mu(\rho)\|_{L^{2}}\\
&\leq C \left\|F(t)-F(s)\right\|_{L^{2}}+\|[\nabla\mu(\rho(t))-\nabla \mu(\rho(s))]\nabla u(t)\|_{L^{2}}+\|\nabla \mu(\rho(s))(\nabla u(t)-\nabla u(s))\|_{L^{2}}\\&\quad
+\|[\tilde{P}(t)-\tilde{P}(s)]\partial_{x}\mu(\rho(t))\|_{L^{2}}+\|\tilde{P}(s)(\partial_{x}\mu(\rho(t))-\partial_{x}\mu(\rho(s)))\|_{L^{2}}
\rightarrow0\quad\text{ as}\quad t\rightarrow s,
\end{aligned}
\end{equation}
which implies that $(u,P)\in C([0,T_{*}];H^{1}\times H^{2})$.

\subsection{Blow-up criterion}
\ \ \
Now we show the blow-up criterion $(\ref{16})$ in Theorem \ref{the-1}.
Introduce the functionals $\Phi(t)$ and $J(t)$ as
$$\Phi(t)=1+\left\|\nabla \rho(t)\right\|_{L^{\infty}}+\left\|\nabla^{2}\rho(t)\right\|_{L^{2}}+\left\|\nabla u(t)\right\|_{L^{2}},$$
and
$$J(t)=1+\left\|u(t)\right\|_{H^{2}}+\left\|P(t)\right\|_{H^{1}}+\left\| \sqrt{\rho } u_{t}(t)\right\|_{L^{2}}+\left\| \rho(t)\right \|_{W^{1, \infty}}+\left\| \rho(t)\right \|_{W^{2, 2}}$$$$
+\int_{\tau}^{t}\left(\left\|u_{t}(s)\right\|_{H^{1}}^{2}+\left\|u(s)\right\|_{H^{3}}^{2}+\left\|P(s)\right\|_{H^{2}}^{2}\right) \mathrm{d}s.$$

Thanks to  the results in Lemma $\ref{Lemma-41}$ and Lemma $\ref{4.1}$-$\ref{4.3}$, we can prove that for any $t\in(\tau,T^{*})$,
$$
\begin{aligned}
&\left\|\nabla u(t)\right\|_{L^{2}}^{2}\leq C(1+\left\|\nabla u(\tau)\right\|_{L^{2}})+C \int_{\tau}^{t}\Phi(s)^{N}ds,\\
&\left\|u(t)\right\|_{H^{2}}+\left\|P(t)\right\|_{H^{1}}\leq C (1+\left\|\sqrt{\rho}u_{t}(t)\right\|_{L^{2}})\Phi(t)^{N},\\
&\left\|u(t)\right\|_{H^{3}}+\left\|P(t)\right\|_{H^{2}}\leq C (\left\|f\right\|_{H^{1}}^{2}+\left\|\nabla u_{t}\right\|_{L^{2}}^{2}+\Phi(t)^{N}),\\
&\left\|\sqrt{\rho}u_{t}(t)\right\|_{L^{2}}^{2}+\int_{\tau}^{t}\left\|\nabla u_{t}(t)\right\|_{L^{2}}^{2}ds\leq C(1+\left\|\sqrt{\rho}u_{t}(\tau)\right\|_{L^{2}})\exp(C \int_{\tau}^{t}\Phi(s)^{N}ds),\\
&\left\|\rho(t)\right\|_{W^{1,\infty}}\leq \left\|\rho(\tau)\right\|_{W^{1,\infty}}\exp(\int_{\tau}^{t}\left\|u\right\|_{H^{3}}ds),\\
&\left\|\rho(t)\right\|_{W^{2,2}}\leq \left\|\rho(\tau)\right\|_{W^{2,2}}\exp(C\int_{\tau}^{t}\left\|u\right\|_{H^{3}}ds)
\end{aligned}
$$
for some $N>0$.

In view of the above estimates, we deduce
$$J(t)\leq \exp\big[C J(\tau)\exp\big(C J(\tau)\int_{\tau}^{t}\Phi(s)^{N}ds\big)\big]$$
for $\tau<t<T^{*}$.

 By the definition of the maximality of $T^{*}$ in Theorem \ref{the-1},  it follows that $J(t)\rightarrow\infty$ as $t\rightarrow T^{*} $. Hence, we can  show that
$$\int_{\tau}^{T^{*}}(\left\|\nabla \rho(t)\right\|_{L^{\infty}}+\left\|\nabla^{2}\rho(t)\right\|_{L^{2}}+\left\|\nabla u\right\|_{L^{2}})^{N}dt=\infty,$$ which implies the blow up criterion $(\ref{16})$. Up to now, we complete the proof of Theorem \ref{the-1}.

{\bf Acknowledgements.} Q. Jiu is partially
supported by National Natural Sciences Foundation of China (No. 11931010). F.Wang is partially supported by the National Natural Science Foundation of China (Grant No.12201028), Chinese Postdoctoral Science Foundation (Grant No.2022M720383)

\end{document}